\documentclass[1p, times]{elsarticle}

\usepackage{amssymb}

\usepackage{amsmath, amssymb, amsgen, amsthm, amscd, xspace, color,tikz,mathrsfs}

\newcommand{\Hom}{\operatorname{Hom}}

\journal{Journal of Functional Analysis}

\begin{document}

\begin{frontmatter}



\title{The explicit theta correspondence for reductive dual pairs $(Sp(p,q),O^{*}(4))$}


\author{Yixin Bao}
\ead{mabaoyixin1984@amss.ac.cn}

\address{Academy of Mathematics and Systems Science, Chinese Academy of Sciences, Beijing, China}

\begin{abstract}
For every irreducible Harish-Chandra module of $O^{*}(4)$, we determine its theta lift to $Sp(p,q)$ in terms of the Langlands parameter, for all non-negative integers $p$ and $q$. Our strategy is to determine the desired theta lifts by their infinitesimal characters and lowest $K$-types.
\end{abstract}


\end{frontmatter}


\section{Introduction}

In representation theory of classical Lie groups, many interesting representations have been constructed via the local theta correspondence for reductive dual pairs (see \cite{Li89LowRank} and \cite{Li89Singular}). To get the Langlands parameters of the representations thus obtained, and for other applications to automorphic forms, it is of great interest to understand the local theta correspondence as explicitly as possible.

For a real reductive group $H$ with complexified Lie algebra $\mathfrak{h}$ and a fixed maximal compact subgroup $K_{H}$, we denote by $\mathcal{R}(H)$ the set of all isomorphism classes of irreducible admissible $(\mathfrak{h},K_{H})$-modules. By Corollary 4.2.4 of \cite{Wallach88ReductiveI}, every admissible  $(\mathfrak{h},K_{H})$-module is isomorphic to the subspace of $K_{H}$-finite vectors of an admissible representation of $H$. By this result, also let $\mathcal{R}(H)$ denote the admissible dual of $H$, i.e., the set of all infinitesimal equivalence classes of irreducible continuous admissible representation of $H$. By abuse of notation, for a class $\pi\in \mathcal{R}(H)$, also let $\pi$ denote a Harish-Chandra module (resp. an admissible representation) in this class.

Let $p,q$ be non-negative integers such that $p+q>0$ and let $n$ be a positive integer. We fix a maximal compact subgroup $K=Sp(p)\times Sp(q)$ of $G=Sp(p,q)$ and a maximal compact subgroup $K^{'}=U(n)$ of $G^{'}=O^{*}(2n)$. The group pair  $(G,G^{'})$ can be embedded into a symplectic group $Sp=Sp(4n(p+q),\mathbb{R})$ such that $G$ and $G^{'}$ are centralizer of each other (as algebraic groups) in this symplectic group. We call $(G,G^{'})$ a reductive dual pair in $Sp$. Furthermore, we fix a maximal compact subgroup $U$ of $Sp$ as in $\S$3 such that $U\cap G=K$ and $U\cap G^{'}=K^{'}$.

Fix a nontrivial unitary character $\psi(t)=e^{2\pi \textbf{i}t}$ of the additive group $\mathbb{R}$ ($\textbf{i}=\sqrt{-1}\in \mathbb C$). A unitary representation $\omega_{\psi}$ of the two-fold metaplectic cover $\widetilde{Sp}$ of $Sp$ was constructed by Shale \cite{Shale62Symmetries} and Weil \cite{Weil64SurCertains}. We call this representation the oscillator representation or the Weil representation. Let $\widetilde{U}$ denote the inverse image of $U$ in $\widetilde{Sp}$. We omit the subscript $\psi$ and denote this Weil representation by $\omega_{p,q,n}$ and its Fock space (the subspace of $\widetilde{U}$-finite vectors of $\omega_{p,q,n}$) by $\mathscr{F}_{p,q,n}$. Since the two-fold metaplectic covers of both $G$ and $G^{'}$ uniquely split, the Weil representation $\omega_{p,q,n}$ can be regarded as a representation of $Sp(p,q)\times O^{*}(2n)$. Let $\mathcal{R}(Sp(p,q),\omega_{p,q,n})$ (resp. $\mathcal{R}(O^{*}(2n),\omega_{p,q,n})$) denote the set of elements in $\mathcal{R}(Sp(p,q))$ (resp. $\mathcal{R}(O^{*}(2n))$) which can be realized in the form $\mathscr{F}_{p,q,n}/\mathscr{N}$, where $\mathscr{N}\subseteq \mathscr{F}_{p,q,n}$ is a $(\mathfrak{sp}(2(p+q),\mathbb{C}),Sp(p)\times Sp(q))$-invariant (resp. $(\mathfrak{o}(2n,\mathbb{C}),U(n))$-invariant) subspace. Similarly, we define the set $\mathcal{R}(Sp(p,q)\times O^{*}(2n),\omega_{p,q,n})$. Howe \cite{Howe89Transcend} proved that $\mathcal{R}(Sp(p,q)\times O^{*}(2n),\omega_{p,q,n})$ is the graph of a bijection between $\mathcal{R}(Sp(p,q),\omega_{p,q,n})$ and $\mathcal{R}(O^{*}(2n),\omega_{p,q,n})$. This bijection is called theta correspondence for the reductive dual pair $(Sp(p,q),O^{*}(2n))$. For $\pi\in \mathcal{R}(Sp(p,q))$ and $\pi^{'}\in \mathcal{R}(O^{*}(2n))$, $\pi$ corresponds to $\pi^{'}$ in this bijection if and only if $\pi\otimes \pi^{'}\in \mathcal{R}(Sp(p,q)\times O^{*}(2n),\omega_{p,q,n})$. If they correspond, we write $\theta_{p,q}(\pi^{'})=\pi$ and $\theta_{n}(\pi)=\pi^{'}$. Furthermore, we write $\theta_{p,q}(\pi^{'})=0$ (resp. $\theta_{n}(\pi)=0$) if $\pi^{'}$ (resp. $\pi$) does not belong to $\mathcal{R}(O^{*}(2n),\omega_{p,q,n})$ (resp. $\mathcal{R}(Sp(p,q),\omega_{p,q,n})$).  Formally, we also define the theta correspondence for $(Sp(0,0),O^{*}(2n))$ and $(Sp(p,q),O^{*}(0))$.

 Fruitful achievements have been made in \cite{AB95Complex}, \cite{Moeglin89Correspondence}, \cite{Paul1998Unitary} and \cite{Paul2005Symplectic} on the explicit theta correspondence for reductive dual pairs of type I. For  $(Sp(p,q),O^{*}(2n))$, the authors of \cite{LPTZ03Correspondence} described the local theta correspondence explicitly for $p+q\leq n$ as well as some of the cases for $p+q>n$. In \cite{LTZ01TensorProduct}, the authors described the local theta correspondence for $(O(p,q),SL(2,\mathbb{R}))$, $(U(p,q),U(1,1))$ and $(Sp(p,q),O^{*}(4))$ in terms of constituents of degenerate principal series. The explicit theta correspondence in terms of Langlands parameters remains unknown. R. Howe and J.-S. Li determined the theta correspondence for $(O(p,q),SL(2,\mathbb{R}))$ in terms of Langlands parameters in their unpublished preprint. In this paper, we will follow their strategy and determine the theta correspondence for $(Sp(p,q),O^{*}(4))$ in terms of Langlands parameters for all non-negative integers $p$ and $q$.

The paper is arranged as follows: In $\S$2.1 and 2.2, we introduce the Langlands parametrization for irreducible admissible representations of $Sp(p,q)$ and $O^{*}(2n)$. Every irreducible admissible representation of $Sp(p,q)$ (resp. $O^{*}(2n)$) is infinitesimal equivalent to a representation $\pi(r,\lambda,\Psi,\mu,\nu)$, where $r$ is a non-negative integer, $\mu\in (\mathbb{Z}_{\geq 1})^{r}$, $\nu\in \mathbb{C}^{r}$, $\lambda$ is the Harish-Chandra parameter for a limit of discrete series $\pi(\lambda,\Psi)$ of $Sp(p-r,q-r)$ (resp. $O^{*}(2(n-2r))$) and $\Psi$ is a system of positive roots of $Sp(p-r,q-r)$ (resp. $O^{*}(2(n-2r))$) satisfying condition (A) of $\S$2.1. In $\S$2.3, we introduce the algorithm for calculating the lowest $Sp(p)\times Sp(q)$-types (resp. $U(n)$-types) of $\pi(r,\lambda,\Psi,\mu,\nu)$.

In $\S$3, we review the correspondence between $Sp(p)\times Sp(q)$-types and $U(n)$-types in the joint harmonics $\mathscr{H}_{p,q,n}$ and the correspondence of infinitesimal characters.  Moreover, we also introduce the following two results for reducing the problem.

\smallskip

\noindent{\bf Theorem A.}(Lemma 3.32 of \cite{LPTZ03Correspondence}){\it~For $\pi^{'}\in \mathcal{R}(O^{*}(2n))$, $(\theta_{p,q}(\pi^{'}))^{\vee}=\theta_{q,p}((\pi^{'})^{\vee})$.}

For a representation $\pi$ of $Sp(p,q)$ (resp. $O^{*}(2n)$), we denote its contragredient representation by $\pi^{\vee}$.

\smallskip

\noindent{\bf Theorem B.}(Theorem 1 of \cite{PrP08FirstOcuurrence}){\it~For $\pi^{'}\in \mathcal{R}(O^{*}(2n))$, $\theta_{p,q}(\pi^{'})$ is nonzero if $\min\{p,q\}\geq n$.}

In $\S$4, we introduce the induction principle and prove the following theorem.

\smallskip

\noindent{\bf Theorem C.}{\it ~For $\pi^{'}\in \mathcal{R}(O^{*}(2n))$, if $\theta_{p,q}(\pi^{'})=\pi(r,\lambda,\Psi,\mu,\nu)$, then
\begin{equation*}
  \theta_{p+s,q+s}(\pi^{'})=\pi(r+s,\lambda,\Psi,\mu^{s},\nu^{s}),
\end{equation*}
where $\mu^{s}=(\mu,\underbrace{1,\cdots,1}_{s})$ and
\begin{equation*}
  \nu^{s}=(\nu,2p+2q-2n+3,2p+2q-2n+7,\cdots,2p+2q-2n+4s-1).
\end{equation*}
}

\smallskip

By Theorem A, Theorem B and Theorem C, we can reduce the problem to explicit description of $\theta_{p,q}$ for $p\geq q$ and $q\leq n$. In $\S$5, we determine $\theta_{p,q}$ explicitly for reductive dual pairs $(Sp(p,q),O^{*}(2))$. This work is done by R. Howe and J.-S. Li in an unpublished preprint. For an integer $k$, let $\chi_{k}$ be the character $t\mapsto t^{k}$ of $O^{*}(2)\cong U(1)$. Since $O^{*}(2)$ is isomorphic to the compact group $U(1)$, every irreducible admissible representation of $O^{*}(2)$ is of the form $\chi_{k}$. We define the occurrence set $\mathcal{O}_{k}=\{(p,q)|\theta_{p,q}(\chi_{k})\neq 0\}$. We give the occurrence pictures of $\chi_{k}$ for $k=0$, $k>0$ and $k<0$ from left to right. The representation $\chi_{k}$ occurs in the theta correspondence for $(Sp(p,q),O^{*}(2))$ if and only if the integral point $(p,q)$ lies in the shadow region or on the bold boundary line. The contragredient representation of $\chi_{k}$ is the representation $\chi_{-k}$. By Lemma A, the occurrence picture of $\chi_{k}$ is just the reflection of the occurrence picture of $\chi_{-k}$ over the line $p=q$.

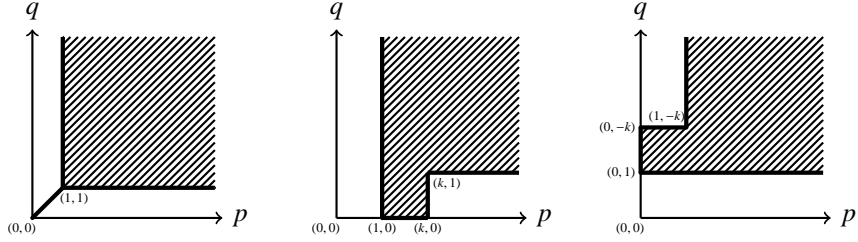
\begin{tikzpicture}[scale=0.2]
\draw[thick, ->](0,0)--(12.5,0) node[right]{$p$};
\draw[thick, ->](0,0) node[scale=0.5, anchor=45]{$(0,0)$}-- (0,12.5) node[above]{$q$};
\draw[ultra thick](0,0)-- (2,2) node[scale=0.5, anchor=135]{$(1,1)$};
\draw[ultra thick](2,2)-- (2,12);
\draw[ultra thick](2,2)-- (12,2);
\fill [black] (0,0) circle(4pt);
\fill [black] (2,2) circle(4pt);
\foreach \x in {2.5,3,...,12}
\draw[thick](2,\x)--(14-\x,12);
\foreach \x in {2,2.5,...,12}
\draw[thick](\x,2)--(12,14-\x);
\draw[thick, ->](20,0)-- (32.5,0) node[right]{$p$};
\draw[thick, ->](20,0) node[scale=0.5,anchor=45]{$(0,0)$}-- (20,12.5) node[above]{$q$};
\draw[ultra thick](23,0)node[scale=0.5,below]{$(1,0)$}-- (23,12);
\draw[ultra thick](26,3)-- (32,3);
\draw[ultra thick](23,0)-- (26,0)node[scale=0.5,below]{$(k,0)$};
\draw[ultra thick](26,0)-- (26,3)node[scale=0.5,anchor=150]{$(k,1)$};
\fill [black] (23,3) circle(4pt);
\fill [black] (23,0) circle(4pt);
\foreach \x in {0.5,1,...,2.5}
\draw[thick](23+\x,0)--(26,3-\x);
\foreach \x in {3,3.5,...,11.5}
\draw[thick](23,\x)--(35-\x,12);
\foreach \x in {0,0.5,1,...,2.5}
\draw[thick](23,\x)--(32,9+\x);
\foreach \x in {0.5,1,...,5.5}
\draw[thick](26+\x,3)--(32,9-\x);
\draw[thick, ->](40,0)-- (52.5,0) node[right]{$p$};
\draw[thick, ->](40,0) node[scale=0.5, anchor=45]{$(0,0)$}-- (40,12.5) node[above]{$q$};
\draw[ultra thick](40,3)node[scale=0.5,left]{$(0,1)$}--(52,3);
\draw[ultra thick](40,3)--(40,6)node[scale=0.5,left]{$(0,-k)$};
\draw[ultra thick](40,6)--(43,6)node[scale=0.5, anchor=-30]{$(1,-k)$};
\draw[ultra thick](43,6)--(43,12);
\fill [black](40,3) circle(4pt);
\fill [black](40,6) circle(4pt);
\fill [black](43,6) circle(4pt);
\foreach \x in {0,0.5,...,3}
\draw[thick](40+\x,3)--(49+\x,12);
\foreach \x in {3.5,4,...,11.5}
\draw[thick](40+\x,3)--(52,15-\x);
\foreach \x in {6.5,7,...,11.5}
\draw[thick](43,\x)--(55-\x,12);
\foreach \x in {3.5,4,...,5.5}
\draw[thick](40,\x)--(46-\x,6);
\end{tikzpicture}

In $\S$6, we determine $\theta_{p,q}$ explicitly for reductive dual pairs $(Sp(p,q),O^{*}(4))$. The following is a list of irreducible admissible representations (up to infinitesimal equvialence) of $O^{*}(4)$:
\begin{enumerate}[(1)]
  \item Irreducible principal series $P_{\lambda_{1},\lambda_{2}}$ with $\lambda_{1},\lambda_{2}\notin \mathbb{Z}$, $\mathrm{Re}(\lambda_{1}+\lambda_{2})\geq 0$ and $\lambda_{1}-\lambda_{2}\in \mathbb{Z}_{\geq 1}$;
  \item Limits of lowest weight discrete series $D_{\lambda_{1},\lambda_{2}}$ with $\lambda_{1},\lambda_{2}\in \mathbb{Z}$, $\lambda_{1}+\lambda_{2}\in \mathbb{Z}_{\geq 0}$ and $\lambda_{1}-\lambda_{2}\in \mathbb{Z}_{\geq 1}$;
  \item Limits of highest weight discrete series $\overline{D}_{\lambda_{1},\lambda_{2}}$ with $\lambda_{1},\lambda_{2}\in \mathbb{Z}$, $\lambda_{1}+\lambda_{2}\in \mathbb{Z}_{\leq 0}$ and $\lambda_{1}-\lambda_{2}\in \mathbb{Z}_{\geq 1}$;
  \item Finite dimensional representations $F_{\lambda_{1},\lambda_{2}}$ with $\lambda_{1},\lambda_{2}\in \mathbb{Z}$, $\lambda_{1}+\lambda_{2}\in \mathbb{Z}_{\geq 1}$ and $\lambda_{1}-\lambda_{2}\in \mathbb{Z}_{\geq 1}$.
\end{enumerate}

The representations in each case are determined by their infinitesimal characters $(\lambda_{1},\lambda_{2})$. We define the occurrence set for them as for those in $\mathcal{R}(O^{*}(2))$. The occurrence sets of all irreducible principle series are the same set $\mathcal{O}_{P}=\{(p,q)|\mathrm{min}\{p,q\}\geq 1\}$. We denote the occurrence set of $F_{\lambda_{1},\lambda_{2}
}$ by $\mathcal{O}_{\lambda_{1},\lambda_{2}
}$, the occurrence set of $D_{\lambda_{1},\lambda_{2}}$ by $\mathcal{O}_{\lambda_{1},\lambda_{2}}^{+}$ and the occurrence set of $\overline{D}_{\lambda_{1},\lambda_{2}}$ by $\mathcal{O}_{\lambda_{1},\lambda_{2}}^{-}$.  By results in $\S$6,  we know that $\mathcal{O}_{\lambda_{1},\lambda_{2}}^{+}=(1,0)+\mathcal{O}_{\lambda_{2}}$
and $\mathcal{O}_{\lambda_{1},\lambda_{2}}^{-}=(0,1)+\mathcal{O}_{\lambda_{1}}$. Then the occurrence picture of $D_{\lambda_{1},\lambda_{2}}$ is just the occurrence picture of $\chi_{\lambda_{2}}$ moved right by one unit and the occurrence picture of $\overline{D}_{\lambda_{1},\lambda_{2}}$ is just the occurrence picture of $\chi_{\lambda_{1}}$ moved up by one unit. The following is a list of occurrence pictures of finite dimensional representations $F_{\lambda_{1},\lambda_{2}}$.

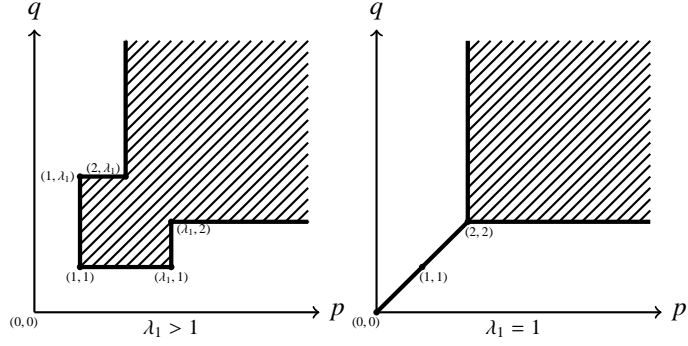
\begin{tikzpicture}[scale=0.3]
\draw[thick, ->](0,0)--(6,0) node[scale=0.8,below]{$\lambda_{1}>1$}-- (12.5,0) node[right]{$p$};
\draw[thick, ->](0,0) node[scale=0.5,anchor=45]{$(0,0)$}-- (0,12.5) node[above]{$q$};
\draw[ultra thick](2,2)node[scale=0.5,below]{$(1,1)$}--(6,2)node[scale=0.5,below]{$(\lambda_{1},1)$};
\draw[ultra thick](6,2)--(6,4)node[scale=0.5,anchor=160]{$(\lambda_{1},2)$};
\draw[ultra thick](6,4)--(12,4);
\draw[ultra thick](2,2)--(2,6)node[scale=0.5,left]{$(1,\lambda_{1})$};
\draw[ultra thick](2,6)--(4,6)node[scale=0.5, anchor=-20]{$(2,\lambda_{1})$};
\draw[ultra thick](4,6)--(4,12);
\fill [black] (2,2) circle(4pt);
\fill [black] (6,2) circle(4pt);
\fill [black] (6,4) circle(4pt);
\fill [black] (2,6) circle(4pt);
\fill [black] (4,6) circle(4pt);
\foreach \x in {2,2.5,...,4}
\draw[thick](\x,2)--(12,14-\x);
\foreach \x in {2.5,3,...,4}
\draw[thick](2,\x)--(14-\x,12);
\foreach \x in {6.5,7,...,12}
\draw[thick](\x,4)--(12,16-\x);
\foreach \x in {6.5,7,...,12}
\draw[thick](4,\x)--(16-\x,12);
\foreach \x in {4.5,5,5.5}
\draw[thick](\x,2)--(6,8-\x);
\foreach \x in {4.5,5,5.5}
\draw[thick](2,\x)--(8-\x,6);
\draw[thick, ->](15,0)--(21,0) node[scale=0.8,below]{$\lambda_{1}=1$}-- (27.5,0) node[right]{$p$};
\draw[thick, ->](15,0) node[scale=0.5,anchor=45]{$(0,0)$}-- (15,12.5) node[above]{$q$};
\draw[ultra thick](15,0)--(17,2)node[scale=0.5,anchor=135]{$(1,1)$}--(19,4)node[scale=0.5,anchor=135]{$(2,2)$};
\draw[ultra thick](19,4)--(19,12);
\draw[ultra thick](19,4)--(27,4);
\fill [black] (15,0) circle(4pt);
\fill [black] (17,2) circle(4pt);
\fill [black] (19,4) circle(4pt);
\foreach \x in {0,0.5,...,7.5}
\draw[thick](19+\x,4)--(27,12-\x);
\foreach \x in {0.5,1,...,7.5}
\draw[thick](19,4+\x)--(27-\x,12);
\end{tikzpicture}
\newline

Let $\Psi_{1}$,$\Psi_{2}$,$\Psi_{3}$,$\Psi_{4}$,$\Psi_{5}$ and $\Psi_{6}$ be the systems of positive roots mentioned in $\S$5 and 6. For irreducible principal series, by Theorem A, Theorem C and the occurrence set $\mathcal{O}_{P}$, we only need to calculate theta $(p,1)$-lifts for $p\geq 1$.

\smallskip

\noindent{\bf Theorem D.}{\it ~The following is a list of theta $(p,1)$-lifts of irreducible principal series for $p\geq 1$:
\begin{equation*}
  \theta_{p,1}(P_{\lambda_{1},\lambda_{2}})=\pi(1,(p-1,p-2,\cdots,1),\Psi_{1},\lambda_{1}-\lambda_{2},\lambda_{1}+\lambda_{2}).
\end{equation*}}

By the occurrence set $\mathcal{O}_{\lambda_{1},\lambda_{2}}^{-}$, we know that all theta $(p,0)$-lifts of limits of highest weight discrete series are zero. For limits of highest weight discrete series, we only need to calculate theta $(p,1)$-lifts and theta $(p,2)$-lifts. We list all theta $(p,1)$-lifts in the following theorem.
\smallskip

\noindent{\bf Theorem E.}{\it~The following is a list of theta $(p,1)$-lifts of limits of highest weight discrete series for $p\geq 1$:
\begin{enumerate}[(1)]
  \item If $\lambda_{1}\geq p$, then
\begin{equation*}
  \theta_{p,1}(\overline{D}_{\lambda_{1},\lambda_{2}})=\pi((\lambda_{1},p-1,p-2,\cdots,1;-\lambda_{2}),\Psi_{2}).
\end{equation*}
  \item If $\lambda_{1}<p$, then
\begin{equation*}
  \theta_{p,1}(\overline{D}_{\lambda_{1},\lambda_{2}})=0.
\end{equation*}
\end{enumerate}}

We list all theta $(p,2)$-lifts of limits of highest weight discrete series in the following theorem.
\smallskip

\noindent{\bf Theorem F.}{\it~The following is a list of theta $(p,2)$-lifts of limits of highest weight discrete series for $p\geq 2$:
\begin{enumerate}[(1)]
  \item If $\lambda_{1}\leq -p$, then
\begin{equation*}
  \theta_{p,2}(\overline{D}_{\lambda_{1},\lambda_{2}})=\pi((p,p-1,\cdots,1;-\lambda_{2},-\lambda_{1}),\Psi_{5}).
\end{equation*}
  \item If $\lambda_{2}\leq -p+1\leq \lambda_{1}<p-1$, then
\begin{equation*}
  \theta_{p,2}(\overline{D}_{\lambda_{1},\lambda_{2}})=\pi(1,(p-1,p-2,\cdots,1;-\lambda_{2}),\Psi_{2},p-\lambda_{1},p+\lambda_{1}).
\end{equation*}
  \item If $\lambda_{2}> -p+1$ and $\lambda_{1}<p-1$, then
\begin{equation*}
  \theta_{p,2}(\overline{D}_{\lambda_{1},\lambda_{2}})=\pi(2,(p-2,p-3,\cdots,1),\Psi_{1},(\mu_{1},\mu_{2}),(\nu_{1},\nu_{2})),
\end{equation*}
with $(\mu_{1},\mu_{2})=(p-1-\lambda_{2},p-\lambda_{1})$ and $(\nu_{1},\nu_{2})=(p-1+\lambda_{2},p+\lambda_{1})$.
  \item If $\lambda_{1}\geq p-1$, then $\theta_{p,2}(\overline{D}_{\lambda_{1},\lambda_{2}})$ can be determined by Theorem C.
\end{enumerate}}

For finite dimensional irreducible representations, we have the following theorem.
\smallskip

\noindent{\bf Theorem G.}{\it~The following is a list of nonzero theta $(p,q)$-lifts of finite dimensional irreducible representations for $p\geq q$ and $q\leq 2$:
\begin{enumerate}[(1)]
  \item $\theta_{0,0}(F_{1,0})=\pi(0,\emptyset)$, where $\pi(0,\emptyset)$ is the trivial representation of the trivial group.
  \item If $\lambda_{1}\geq p$, then
\begin{equation*}
  \theta_{p,1}(F_{\lambda_{1},\lambda_{2}})=\pi(1,(p-1,p-2,\cdots,1),\Psi_{1},\lambda_{1}-\lambda_{2},\lambda_{1}+\lambda_{2}).
\end{equation*}
  \item If $\lambda_{1}<p-1$, then
\begin{equation*}
  \theta_{p,2}(F_{\lambda_{1},\lambda_{2}})=\pi(2,(p-2,p-3,\cdots,1),\Psi_{1},(\mu_{1},\mu_{2}),(\nu_{1},\nu_{2})),
\end{equation*}
with $(\mu_{1},\mu_{2})=(p-1-\lambda_{2},p-\lambda_{1})$ and $(\nu_{1},\nu_{2})=(p-1+\lambda_{2},p+\lambda_{1}))$.
\item If $\lambda_{1}\geq p-1$, then $\theta_{p,2}(F_{\lambda_{1},\lambda_{2}})$ can be determined by Theorem C.
\end{enumerate}}

For limits of lowest weight discrete series, first we list all nonzero theta $(p,0)$-lifts.
\smallskip

\noindent{\bf Theorem H.}{\it~The following is a list of nonzero theta $(p,0)$-lifts of limits of lowest weight discrete series:
\begin{enumerate}[(1)]
  \item $\theta_{1,0}(D_{\lambda_{1},0})=\pi((\lambda_{1}),\Psi_{1})$.
  \item If $p\geq 2$ and $\lambda_{2}\geq p-1$,  then
\begin{equation*}
  \theta_{p,0}(D_{\lambda_{1},\lambda_{2}})=\pi((\lambda_{1},\lambda_{2},p-2,p-3,\cdots,1),\Psi_{1}).
\end{equation*}
\end{enumerate}}

Next we list all theta $(p,1)$-lifts for $p\geq 1$.

\smallskip

\noindent{\bf Theorem I.}{\it~The following is a list of theta $(p,1)$-lifts of limits of lowest weight discrete series for $p\geq 1$:
\begin{enumerate}[(1)]
  \item If $\lambda_{1}\leq p-1$ and $\lambda_{1}+\lambda_{2}>0$, then
\begin{equation*}
  \theta_{p,1}(D_{\lambda_{1},\lambda_{2}})=\pi(1,(p-1,p-2,\cdots,1),\Psi_{1},\lambda_{1}-\lambda_{2},\lambda_{1}+\lambda_{2}).
\end{equation*}
  \item If $\lambda_{1}\leq p-1$ and $\lambda_{1}+\lambda_{2}=0$,
  then
\begin{equation*}
  \theta_{p,1}(D_{\lambda_{1},-\lambda_{1}})=\pi((p-1,\cdots,\lambda_{1}+1,\lambda_{1},\lambda_{1},\lambda_{1}-1,\cdots,1;\lambda_{1}),\Psi)
\end{equation*}
with $\Psi$ uniquely determined by condition $(A)$ of $\S2.1$.
  \item If $\lambda_{1}\geq p$ and $\lambda_{2}\leq -p$, then
\begin{equation*}
  \theta_{p,1}(D_{\lambda_{1},\lambda_{2}})=\pi((\lambda_{1},p-1,p-2,\cdots,1;-\lambda_{2}),\Psi_{4}).
\end{equation*}
  \item If $p\geq 2$, $1-p<\lambda_{2}\leq p-2$ and $\lambda_{1}\geq p$, then
\begin{equation*}
  \theta_{p,1}(D_{\lambda_{1},\lambda_{2}})=\pi(1,(\lambda_{1},p-2,p-3,\cdots,1),\Psi_{1},p-1-\lambda_{2},p-1+\lambda_{2}).
\end{equation*}
  \item If $p\geq 2$, $\lambda_{1}\geq p$ and $\lambda_{2}=1-p$, then
\begin{equation*}
  \theta_{p,1}(D_{\lambda_{1},1-p})=\pi((\lambda_{1},p-1,p-2,\cdots,1;p-1),\Psi_{4}).
\end{equation*}
  \item If $p\geq 3$ and $\lambda_{2}\geq p-1$, then
\begin{equation*}
  \theta_{p,1}(D_{\lambda_{1},\lambda_{2}})=\pi(1,(\lambda_{1},\lambda_{2},p-3,p-4,\cdots,1),\Psi_{1},1,2p-3).
\end{equation*}
  \item If $p=2$ and $\lambda_{2}\geq 1$, then
\begin{equation*}
  \theta_{2,1}(D_{\lambda_{1},\lambda_{2}})=\pi((\lambda_{1},\lambda_{2};1),\Psi_{3}).
\end{equation*}
  \item If $p=1$ and $\lambda_{2}\geq 0$, then
\begin{equation*}
  \theta_{1,1}(D_{\lambda_{1},\lambda_{2}})=0.
\end{equation*}
\end{enumerate}}

Finally, we calculate theta $(p,2)$-lifts of limits of lowest weight discrete series for $p\geq 2$. By Theorem C and Theorem I, we only need to calculate theta $(2,2)$-lifts of $D_{\lambda_{1},\lambda_{2}}$ for $\lambda_{2}\geq 0$.

\noindent{\bf Theorem J.}{\it~The following is a list of theta $(2,2)$-lifts of $D_{\lambda_{1},\lambda_{2}}$ for $\lambda_{2}\geq 0$:
\begin{enumerate}[(1)]
  \item If $\lambda_{2}=0$, then
\begin{equation*}
  \theta_{2,2}(D_{\lambda_{1},0})=\pi(1,(\lambda_{1};1),\Psi_{3},2,2).
\end{equation*}
  \item If $\lambda_{2}=1$, then
\begin{equation*}
  \theta_{2,2}(D_{\lambda_{1},1})=\pi(1,(\lambda_{1};1),\Psi_{3},3,1).
\end{equation*}
  \item If $\lambda_{2}\geq 2$, then
\begin{equation*}
  \theta_{2,2}(D_{\lambda_{1},\lambda_{2}})=\pi((\lambda_{1},\lambda_{2};2,1),\Psi_{6}).
\end{equation*}
\end{enumerate}}

Notice that $\pi=\pi^{\vee}$ if $\pi$ is an irreducible admissible representation of $Sp(p,q)$. On the other hand, $(P_{\lambda_{1},\lambda_{2}})^{\vee}=P_{\lambda_{1},\lambda_{2}}$, $(F_{\lambda_{1},\lambda_{2}})^{\vee}=F_{\lambda_{1},\lambda_{2}}$ and $(D_{\lambda_{1},\lambda_{2}})^{\vee}=\overline{D}_{-\lambda_{2},-\lambda_{1}}$. By the ten theorems and the occurrence sets given in this introduction, we  determine theta $(p,q)$-lifts of all irreducible Harish-Chandra module of $O^{*}(4)$ explicitly, for all non-negative integers $p$ and $q$.

\section{Langlands Parameters and Lowest $K$-types}

For a compact Lie group $K$, a $K$-type means an irreducible finite dimensional representation (up to isomorphism) of $K$. Let $H$ be a real reductive group with a maximal compact subgroup $K_{H}$. We sometimes refer to $K_{H}$-types as $K$-types for $H$, or simply as $K$-types if the group $H$ is clearly understood. We denote by $\mathcal{R}(H)$ the set of all isomorphism classes of irreducible Harish-Chandra modules of $H$.  Since every Harish-Chandra module of $H$ can be realized as the subspace of $K_{H}$-finite vectors of an admissible representation of $H$ (see Corollary 4.2.4 of \cite{Wallach88ReductiveI}), $\mathcal{R}(H)$ can also be regarded as the set of infinitesimal equivalence classes of irreducible admissible representations of $H$. By abuse of notation, for an isomorphism class (resp. an infinitesimal equivalence class) $\pi$ in $\mathcal{R}(H)$, also let $\pi$ denote an irreducible Harish-Chandra module (resp. an irreducible admissible representation) in this class. For every $\pi\in \mathcal{R}(H)$, we denote by $\mathcal{K}(\pi)$ the set of $K_{H}$-types occurring in $\pi$ and by $\mathcal{A}(\pi)$ the set of lowest $K_{H}$-types of $\pi$ (see Definition 5.1 of \cite{Vogan79AlgebraicStructure}). For $\pi\in \mathcal{R}(H)$ and $\sigma\in \mathcal{K}(\pi)$,  $m(\sigma,\pi)=\dim\Hom_{K_{H}}(\sigma,\pi)$ is a non-negative integer and we call it the multiplicity of $\sigma$ in $\pi$.

Let $\mathbb{H}$ be a quaternion algebra such that $\mathbb{H}=\mathbb{R}\oplus \mathbb{R}\textbf{i}\oplus \mathbb{R}\textbf{j}\oplus \mathbb{R}\textbf{k}$ with $\textbf{i}^{2}=\textbf{j}^{2}=\textbf{k}^{2}=-1$, $\textbf{ij}=\textbf{k}$, $\textbf{jk}=\textbf{i}$ and $\textbf{ki}=\textbf{j}$. We define a involution on $\mathbb{H}$ such that $\overline{a+b\textbf{i}+c\textbf{j}+d\textbf{k}}=a-b\textbf{i}-c\textbf{j}-d\textbf{k}$.

In $\S$2.1 and $\S$2.2, we briefly review the Langlands parametrization for admissible dual of the groups $Sp(p,q)$ and $O^{*}(2n)$. We mainly follow the content in $\S$2 of \cite{LPTZ03Correspondence} and a more detailed discussion can be found in $\S$3 of \cite{Paul2005Symplectic}. In $\S$2.3, we introduce the algorithm for calculating the lowest $K$-types of all irreducible admissible representations of $Sp(p,q)$ and $O^{*}(2n)$. In $\S$2.4, we state a more precise classification of admissible dual of the group $O^{*}(4)$. This work plays a key role in the explicit calculation of the local theta correspondence.

\subsection{Parametrization for $\mathcal{R}(Sp(p,q))$}

 For every positive integer $n$, we denote by $I_{n}$ the $n\times n$ identity matrix. For fixed non-negative integers $p$ and $q$, let $G=Sp(p,q)$ be the isometry group of the Hermitian form $(,)$ on right $\mathbb{H}$ column vector space $\mathbb{H}^{p+q}$ given by
\begin{equation}\label{Hermitian form}
  (v,w)=\sum_{i=1}^{p}\bar{v_{i}}w_{i}-\sum_{i=p+1}^{p+q}\bar{v_{i}}w_{i}
\end{equation}
for $v,w\in \mathbb{H}^{p+q}$. Then
\begin{equation*}
  Sp(p,q)=\{g\in GL(p+q,\mathbb{H)}|g^{*}K_{p,q}g=K_{p,q}\}
\end{equation*}
with
\begin{equation*}
  K_{p,q}=\left(
                                               \begin{array}{cc}
                                                 I_{p} & 0 \\
                                                 0 & -I_{q} \\
                                               \end{array}
                                             \right).
\end{equation*}
Here $g^{*}=\overline{g}^{T}$. The group $Sp(p,q)$ has a natural left action on the right $\mathbb{H}$ column vector space $\mathbb{H}^{p+q}$.
Let
\begin{equation}\label{Maximal Compact Subgroup1}
  K=\{(g_{1},g_{2})\in GL(p,H)\times GL(q,H)|g_{1}^{*}I_{p}g_{1}=I_{p},g_{2}^{*}I_{q}g_{2}=I_{q}\}.
\end{equation}
We know that $K=Sp(p)\times Sp(q)$ is a maximal compact subgroup of $Sp(p,q)$ (see $\S$\uppercase\expandafter{\romannumeral1}.1 of \cite{Knapp86Overview}). We denote by $S^{1}$ the set of complex numbers whose absolute values are one. Let
\begin{equation*}
  T=\left\{\left(
        \begin{array}{ccc}
          a_{1} &  &  \\
           & \ddots &  \\
           &  & a_{p+q} \\
        \end{array}
      \right)
|a_{i}\in S^{1}\right\}.
\end{equation*}
Then $T$ is a maximal Cartan subgroup of $K$. The real Lie algebra $\mathfrak{g}_{0}$ of $G$ is
\begin{equation}\label{Real Lie Algebra1}
  \left\{\left(\begin{array}{cc}
       A & B \\
       C & D
     \end{array}\right)
  \mid A+A^{*}=D+D^{*}=0,C=B^{*}\right\}
\end{equation}
with $A\in M_{p\times p}(\mathbb{H})$, $B\in M_{p\times q}(\mathbb{H})$, $C\in M_{q\times p}(\mathbb{H})$ and $D\in M_{q\times q}(\mathbb{H})$.

For every positive integer $n$, there is an injective homomorphism from $GL(n,\mathbb{H})$ (resp. $M_{n\times n}(\mathbb{H})$) into $GL(2n,\mathbb{C})$ (resp. $M_{2n\times 2n}(\mathbb{C})$):
\begin{equation}\label{Injective Homomorphism}
  A+B\textbf{j}\mapsto \left(
                 \begin{array}{cc}
                   A & B \\
                   -\overline{B} & \overline{A} \\
                 \end{array}
               \right)
\end{equation}
with $A,B\in M_{n\times n}(\mathbb{C})$. By this homomorphism, the group $Sp(p,q)$ is embedded into $GL(2(p+q),\mathbb{C})$ as a closed subgroup and the real Lie algebra $\mathfrak{g}_{0}$ is isomorphic to
\begin{equation*}
  \left\{ \left(
    \begin{array}{cc}
      A & B \\
      -\overline{B} & \overline{A} \\
    \end{array}
  \right)
| A^{*}K_{p,q}+K_{p,q}A=0, \, B^{t}K_{p,q}=K_{p,q}B \right\}
\end{equation*}
with $A,B\in M(p+q,\mathbb{C})$. Here
\begin{equation*}
  A=\left(
          \begin{array}{cc}
            A_{1} & A_{2} \\
            A_{3} & A_{4} \\
          \end{array}
        \right)
\end{equation*}
with $A_{1}\in M_{p}(\mathbb{C})$, $A_{4}\in M_{q}(\mathbb{C})$, $A_{2}\in M_{p\times q}(\mathbb{C})$ and $A_{3}\in M_{q\times p}(\mathbb{C})$. We know that
\begin{equation*}
  A_{1}+A_{1}^{*}=0, A_{4}+A_{4}^{*}=0, A_{3}=A_{2}^{*}.
\end{equation*}
Similarly,
\begin{equation*}
  B=\left(
                \begin{array}{cc}
                  B_{1} & B_{2} \\
                  B_{3} & B_{4} \\
                \end{array}
              \right)
\end{equation*}
with $B_{1}\in M_{p}(\mathbb{C})$, $B_{4}\in M_{q}(\mathbb{C})$, $B_{2}\in M_{p\times q}(\mathbb{C})$ and $B_{3}\in M_{q\times p}(\mathbb{C})$. We know that
\begin{equation*}
  B_{1}=B_{1}^{t}, B_{4}=B_{4}^{t}, B_{3}=-B_{2}^{t}.
\end{equation*}
The real Lie algebra $\mathfrak{k}_{0}$ of $K$ is the set of matrices in $\mathfrak{g}_{0}$ such that $A_{2}=B_{2}=0$ and the real Lie algebra $\mathfrak{t}_{0}$ of $T$
 is the set of diagonal matrices in $\mathfrak{g}_{0}$. Let $\mathfrak{g}_{0}=\mathfrak{k}_{0}\oplus {\mathfrak{p}}_{0}$ be the Cartan decomposition of $\mathfrak{g}_{0}$. We omit the subscript $0$ to denote the complexified Lie algebras. We choose a basis $\{e_{1},\cdots,e_{p};f_{1},\cdots,f_{q}\}$ of $\textbf{i}\mathfrak{t}_{0}^{*}$ such that the set of compact roots is $\Delta_{c}=\{\pm 2e_{i},\pm 2f_{i},\pm(e_{i}\pm e_{j}),\pm(f_{i}\pm f_{j})\}$ and the set of non-compact roots is $\Delta_{n}=\{\pm (e_{i}\pm f_{j})\}$. We choose a system of positive compact roots $\Delta_{c}^{+}=\{2e_{i},2f_{i},e_{i}\pm e_{j},f_{i}\pm f_{j}|i<j\}$. One-half the sum of positive compact roots is
\begin{equation*}
  \rho_{c}=(p,p-1,\cdots,1;q,q-1,\cdots,1).
\end{equation*}

The discrete series of $G$ is parameterized by Harish-Chandra parameters $\lambda$ as follows:
\begin{equation*}
  (a_{1},a_{2},\cdots,a_{p};b_{1},b_{2},\cdots,b_{q})
\end{equation*}
with $a_{i},b_{j}\in \mathbb{Z}$,$a_{1}>a_{2}>\cdots >a_{p}\geq 1$, $b_{1}>b_{2}>\cdots >b_{q}\geq 1$ and $a_{i}\neq b_{j}$ for all $i,j$.

A limit of discrete series of $G$ is parameterized by a pair $(\lambda,\Psi)$, where $\lambda$ is the Harish-Chandra parameter and $\Psi$ is a system of positive roots containing $\Delta_{c}^{+}$. The parameter $\lambda=(\lambda_{1};\lambda_{2})$ is of the form
\begin{equation}\label{Harish-Chandra Parameter1}
  (\overbrace{a_{1},\cdots,a_{1}}^{m_{1}},\overbrace{a_{2},\cdots,a_{2}}^{m_{2}},\cdots,\overbrace{a_{k},\cdots,a_{k}}^{m_{k}};\overbrace{a_{1},\cdots,a_{1}}^{n_{1}},\cdots,\overbrace{a_{k},\cdots,a_{k}}^{n_{k}}),
\end{equation}
where $a_{i}\in \mathbb{Z}$, $a_{1}>a_{2}>\cdots >a_{k}>0$ and $|m_{i}-n_{i}|\leq 1$ for all $i$. Here $\Psi$ is a system of positive roots satisfying the following condition (also see condition F-1 of \cite{Vogan84Unitarzability}):
\newline
(A)$\left\{
  \begin{array}{ll}
     \hbox{The system of compact positive roots $\Delta_{c}^{+}$ is contained in $\Psi$;} \\
     \hbox{The Harish-Chandra parameter $\lambda$ is dominant with respect to $\Psi$;} \\
     \hbox{If $<\lambda,\alpha>=0$ for a simple root $\alpha$ in $\Psi$, then $\alpha$ is non-compact.}
  \end{array}
\right.$
\newline
Consequently, there are $2^{t}$ different systems of positive roots satisfying condition (A), where $t$ is the number of $i$ such that $m_{i}=n_{i}>0$ (c.f. $\S$3.1 of \cite{Paul2005Symplectic}). We say  $\lambda$ satisfies condition (B) if $t=0$. If $\lambda$ satisfies condition (B), there is a unique system of positive roots $\Psi$ satisfying condition (A). We denote by $\pi(\lambda,\Psi)$ the limit of discrete series parameterized by $(\lambda,\Psi)$. Its unique lowest $K$-type is $\Lambda=\lambda+\rho_{n}-\rho_{c}$, where $\rho_{n}$ is one-half the sum of non-compact roots in $\Psi$ (We identify $K$-types with their highest weight with respect to $\Delta_{c}^{+}$.).

Cuspidal parabolic subgroups of $Sp(p,q)$ are of the form $P=MAN=LN$ with
\begin{equation}\label{Cuspidal Parabolic Subgroup1}
  L\cong Sp(p-r,q-r)\times GL(1,\mathbb{H})^{r}
\end{equation}
and $0\leq r\leq \min\{p,q\}$. In fact, set $R=\min\{p,q\}$. Then $(\mathbb{H}^{p+q},(,))$ admits a polar decomposition
\begin{equation}\label{Polar Decomposition1}
  \mathbb{H}^{p+q}=V_{R}^{+}\oplus V^{0}\oplus V_{R}^{-},
\end{equation}
where $V_{R}^{+}=\mathrm{span}_{\mathbb{H}}\{v_{1},v_{2},\cdots,v_{R}\}$ and $V_{R}^{-}=\mathrm{span}_{\mathbb{H}}\{v_{1}^{'},v_{2}^{'},\cdots,v_{R}^{'}\}$ are two maximal isotropic subspaces in $\mathbb{H}^{p+q}$ and are dual to each other with respect to $(,)$. Here $v_{i}$, $v_{i}^{'}$ can be chosen such that $(v_{i},v_{j}^{'})=\delta_{ij}$.
For every integer $r$ such that $0\leq r\leq R$, we define
\begin{equation}\label{Isotropic Subspace1}
  V_{r}^{+}=\mathrm{span}_{\mathbb{H}}\{v_{1},v_{2},\cdots,v_{r}\},
\end{equation}
\begin{equation}\label{Isotropic Subspace2}
  V_{r}^{-}=\mathrm{span}_{\mathbb{H}}\{v_{1}^{'},v_{2}^{'},\cdots,v_{r}^{'}\},
\end{equation}
and define $V^{r}$ to be
\begin{equation}\label{Complementary Subspace}
\begin{split}
   \mathrm{span}_{\mathbb{H}}\{v_{R-r+1},\cdots,v_{R-1},v_{R}\}\oplus V^{0}\oplus \mathrm{span}_{\mathbb{H}}\{v_{R-r+1}^{'},\cdots,v_{R-1}^{'},v_{R}^{'}\}.
\end{split}
\end{equation}
Then
\begin{equation}\label{Isotropic Flag1}
  \{0\}\subset V_{1}^{+}\subset V_{2}^{+}\subset\cdots \subset V_{r}^{+}
\end{equation}
is an isotropic flag in $V_{R}^{+}$. The stabilizer $P_{r}$ of this flag in $Sp(p,q)$ is a cuspidal parabolic subgroup of $Sp(p,q)$. Furthermore, we embed the group $Sp(p-r,q-r)\times GL(1,\mathbb{H})^{r}$ into $P_{r}$. The group $Sp(p-r,q-r)$ is regarded as the isometry group of $V^{R-r}$. The group $GL(1,\mathbb{H})^{r}$ is equipped with an left action $\circ$ on $\mathbb{H}^{p+q}$ which is right $\mathbb{H}$-linear as follows. For $h=(h_{1},\cdots,h_{r})\in GL(1,\mathbb{H})^{r}$, the action of $h$ on $\mathbb{H}^{p+q}$ is
\begin{equation}\label{Left Action}
\begin{split}
 h\circ v_{i}=v_{i}\cdot h_{i}, \quad \forall 1\leq i\leq r,\\
 h\circ v_{i}^{'}=v_{i}^{'}\cdot (\overline{h_{i}})^{-1}, \quad \forall 1\leq i\leq r,\\
 h\circ u=u, \quad \forall u\in V^{R-r},
\end{split}
\end{equation}
where $\cdot$ is the natural right multiplication on the right $\mathbb{H}$ column vector space $\mathbb{H}^{p+q}$. It is easy to check that
\begin{equation*}
  (h\circ v_{i},h\circ v_{i}^{'})=(v_{i},v_{i}^{'})
\end{equation*}
for each $h\in GL(1,\mathbb{H})^{r}$ and each $1\leq i\leq r$. Then we embed $Sp(p-r,q-r)\times GL(1,\mathbb{H})^{r}$ into the parabolic subgroup $P_{r}$. In general, let
\begin{equation}\label{Isotropic Flag2}
  \{0\}\subset V_{r_{1}}^{+}\subset V_{r_{2}}^{+}\subseteq \cdots\subset V_{r_{s}}^{+}
\end{equation}
be an isotropic flag in $V_{R}^{+}$, where $1\leq r_{1}<r_{2}<\cdots <r_{s}\leq R$. The stabilizer $P_{\{r_{1},\cdots,r_{s}\}}$ of this flag is a parabolic subgroup of $Sp(p,q)$. Set $r_{0}=0$ and $d_{i}=r_{i}-r_{i-1}$ for $1\leq i\leq s$. The Levi subgroup of $P_{\{r_{1},\cdots,r_{s}\}}$ is isomorphic to
\begin{equation}\label{Cuspidal Parabolic Subgroup2}
  Sp(p-r_{s},q-r_{s})\times \Pi_{i=1}^{s}GL(d_{i},\mathbb{H}).
\end{equation}
Up to conjugation, each proper parabolic subgroup of $Sp(p,q)$ is the stabilizer of an isotropic flag.

By the injective homomorphism (\ref{Injective Homomorphism}), the group $GL(1,\mathbb{H})$ is
\begin{equation*}
  \left\{\left(
     \begin{array}{cc}
       a & b \\
       -\bar{b} & \bar{a} \\
     \end{array}
   \right)|a,b\in \mathbb{C}\right\}\cap GL(2,\mathbb{C}).
\end{equation*}
We define the determinant of an element in $GL(1,\mathbb{H})$ to be $|a|^{2}+|b|^{2}$. For every non-negative integer $m$,  we denote the $m+1$-dimensional irreducible representation of $SU(2)$ by $\tau_{m}$. Since $GL(1,\mathbb{H})\cong SU(2)\times \mathbb{R}$, every irreducible admissible representation of  $GL(1,\mathbb{H})$ is of the form $\tau(\mu,\nu)=\tau_{\mu-1}\otimes \det(\cdot)^{\frac{\nu}{2}}$ with $\mu\in {\mathbb{Z}}_{\geq 1}$ and $\nu\in \mathbb{C}$.  The infinitesimal character of $\tau(\mu,\nu)$ is $(\frac{\nu+ \mu}{2},\frac{\nu- \mu}{2})$.

Given a cuspidal parabolic subgroups $P=LN$ of $Sp(p,q)$ with the Levi component as in (\ref{Cuspidal Parabolic Subgroup1}), a limit of discrete series $\varpi=\pi(\lambda,\Psi)$ of $Sp(p-r,q-r)$ and a representation $\tau=\otimes _{i=1}^{r}\tau(\mu_{i},\nu_{i})$ of $GL(1,\mathbb{H})^{r}$, we call the normalized induced representation
\begin{equation}\label{Standard Module}
  Ind_{P}^{Sp(p,q)}(\varpi\otimes \tau\otimes 1{\kern-0.36em} 1)
\end{equation}
a standard module and denote it by $X(P,r,\lambda,\Psi,\mu,\nu)$ with $\mu=(\mu_{1},\cdots,\mu_{r})\in (\mathbb{Z}_{\geq 1})^{r}$ and $\nu=(\nu_{1},\cdots,\nu_{r})\in \mathbb{C}^{r}$. Since there is a unique cuspidal parabolic subgroup up to conjugation in $Sp(p,q)$ whose Levi component is isomorphic to $Sp(p-r,q-r)\times GL(1,\mathbb{H})^{r}$, for convenience, we usually omit the notation $P$ and denote the standard module by $X(r,\lambda,\Psi,\mu,\nu)$. By abuse of notation, we also denote the underlying Harish-Chandra module of the standard module by $X(r,\lambda,\Psi,\mu,\nu)$. The non-parity condition F-2 in \cite{Vogan84Unitarzability} amounts to the requirement that $\mu_{i}$ is odd if $\nu_{i}=0$. If the non-parity condition is satisfied, the standard module $X(r,\lambda,\Psi,\mu,\nu)$ has a unique constituent containing all lowest $Sp(p)\times Sp(q)$-types. We denote this constituent (infinitesimal equivalence class) by $\pi(r,\lambda,\Psi,\mu,\nu)$.

We identify the infinitesimal characters of Harish-Chandra modules (resp. admissible representations) of $G$ with elements of the dual of a Cartan subalgebra of $\mathfrak{g}$ (module the action of Weyl group), via the Harish-Chandra map. For $Sp(p,q)$, we choose the maximal compact Cartan subalgebra $\mathfrak{t}$ for our Cartan subalgebra. Since the infinitesimal characters are preserved by normalized induction, the infinitesimal character of the standard module $X(r,\lambda,\Psi,\mu,\nu)$ and $\pi(r,\lambda,\Psi,\mu,\nu)$ is
\begin{equation}\label{Infinitesimal Characters}
  (\lambda,\frac{\nu_{1}+\mu_{1}}{2},\frac{\nu_{1}-\mu_{1}}{2},\frac{\nu_{2}+\mu_{2}}{2},\frac{\nu_{2}-\mu_{2}}{2},\cdots,\frac{\nu_{r}+\mu_{r}}{2},\frac{\nu_{r}-\mu_{r}}{2})\in \mathfrak{t}^{*}/\sim,
\end{equation}
where $\sim$ means the equivalence under the action of Weyl group $W$ of $\triangle(\mathfrak{g}:\mathfrak{t})$.

For the classification of admissible dual of $Sp(p,q)$, we have the following theorem (see $\S$2.2 of \cite{LPTZ03Correspondence}).

\newtheorem{thm}{Theorem}[section]
\newdefinition{rmk}[thm]{Remark}

\begin{thm}\label{Langlands Classification}
Every irreducible admissible representation of $Sp(p,q)$ is infinitesimal equivalent to some $\pi(r,\lambda,\Psi,\mu,\nu)$ as mentioned above. Furthermore, $\pi(r,\lambda,\Psi,\mu,\nu)$ is infinitesimal equivalent to $\pi(r^{'},\lambda^{'},\Psi^{'},\mu^{'},\nu^{'})\Leftrightarrow r=r^{'}$, $\lambda=\lambda^{'}$, $\Psi=\Psi^{'}$ and $(\mu^{'},\nu^{'})$ may be obtained from $(\mu,\nu)$ by simultaneous permutation of the coordinates, and by replacing some of the $\nu_{i}$ by $-\nu_{i}$.
\end{thm}

Thus we obtain a parametrization for $\mathcal{R}(Sp(p,q))$.

\subsection{Parametrization for $\mathcal{R}(O^{*}(2n))$}

For a fixed positive integer $n$, let $G=O^{*}(2n)$ be the isometry group of the skew-Hermitian form $<,>$ on left $\mathbb{H}$ row vector space $\mathbb{H}^{n}$ given by
\begin{equation}\label{Skew-Hermitian Form}
  <v,w>=\sum_{i=1}^{n}v_{i}\textbf{i}\bar{w_{i}}
\end{equation}
for $v,w\in \mathbb{H}^{n}$. The group $GL(n,\mathbb{H})$ has a natural right action on $\mathbb{H}^{n}$.  Then
\begin{equation}\label{Isometry Group of Skew-Hermitian form}
  O^{*}(2n)=\{g\in GL(n,\mathbb{H)}|g\textbf{i}I_{n}g^{*}=\textbf{i}I_{n}\}.
\end{equation}
By (\ref{Isometry Group of Skew-Hermitian form}), it is not difficult to see that $K=U(n)$ is naturally embedded into $O^{*}(2n)$ as a maximal compact group (see $\S$\uppercase\expandafter{\romannumeral1}.1 of \cite{KnVo95Cohomological}). Furthermore, $T=U(1)^{n}$ is naturally embedded into $O^{*}(2n)$ and is regarded as a maximal Cartan subgroup of $K$. By the homomorphism (\ref{Injective Homomorphism}), $O^{*}(2n)$ can be embedded into $GL(2n,\mathbb{C})$ as a closed subgroup, while
\begin{equation}\label{Maximal Compact Subgroup2}
  K=\left\{\left(
              \begin{array}{cc}
                A & 0 \\
                0 & \overline{A} \\
              \end{array}
            \right)|A\in U(n)
  \right\}.
\end{equation}
The real Lie algebra $\mathfrak{g}_{0}=\mathfrak{o}^{*}(2n)$ of $O^{*}(2n)$ is
\begin{equation*}
  \{g\in M_{n\times n}(\mathbb{H})|g\textbf{i}I_{n}+\textbf{i}I_{n}g^{*}=0\}.
\end{equation*}
By the homomorphism (\ref{Injective Homomorphism}), the real Lie algebra $\mathfrak{g}_{0}$ is isomorphic to
\begin{equation}\label{Real Lie Algebra2}
  \left\{\left(
     \begin{array}{cc}
       A & B \\
       -\overline{B} & \overline{A} \\
     \end{array}
   \right)
| A+A^{*}=0, B+B^{t}=0\right\}
\end{equation}
with $A,B\in M_{n\times n}(\mathbb{C})$. The real Lie algebra $\mathfrak{k}_{0}$ of $K$ is the set of matrices in $\mathfrak{g}_{0}$ such that $B=0$ and the real Lie algebra $\mathfrak{t}_{0}$ of $T$
 is the set of diagonal matrices in $\mathfrak{g}_{0}$. Let $\mathfrak{g}_{0}=\mathfrak{k}_{0}\oplus {\mathfrak{p}}_{0}$ be the Cartan decomposition. We omit the subscript $0$ to denote the complexified Lie algebras.

We choose a basis $\{e_{1},\cdots,e_{n}\}$ of $i\mathfrak{t}_{0}^{*}$ such that the set of compact roots is $\Delta_{c}=\{\pm(e_{i}-e_{j})\mid i<j\}$ and the set of non-compact roots is $\Delta_{n}=\{\pm(e_{i}+e_{j})\mid i<j\}$.  We choose a system of positive compact roots $\Delta_{c}^{+}=\{e_{i}-e_{j}\mid i<j\}$. One-half the sum of positive compact roots is
\begin{equation*}
  \rho_{c}=(\frac{n-1}{2},\frac{n-3}{2},\cdots,\frac{1-n}{2}).
\end{equation*}

The discrete series of $G$ is parameterized by Harish-Chandra parameters $\lambda$ as follows:
\begin{equation*}
  (a_{1},a_{2},\cdots,a_{n})
\end{equation*}
with $a_{i}\in \mathbb{Z}$,$a_{1}>a_{2}>\cdots >a_{n}$ and $a_{i}+a_{j}\neq 0$.

A limit of discrete series of $G$ is parameterized by a pair $(\lambda,\Psi)$, where $\lambda$ is the Harish-Chandra parameter and $\Psi$ is a positive root system containing $\Delta_{c}^{+}$. The parameter $\lambda
$ is of the form
\begin{equation}\label{Harish-Chandra Parameter2}
  (\overbrace{a_{1},\cdots,a_{1}}^{m_{1}},\overbrace{a_{2},\cdots,a_{2}}^{m_{2}},\cdots,\overbrace{a_{k},\cdots,a_{k}}^{m_{k}},\overbrace{-a_{k},\cdots,-a_{k}}^{n_{k}},\cdots,,\overbrace{-a_{1},\cdots,-a_{1}}^{n_{1}})
\end{equation}
with $a_{i}$ as in (\ref{Harish-Chandra Parameter1}); or with one zero between the last $a_{k}$ and the first $-a_{k}$. Here $\Psi$ is a system of positive roots satisfying condition (A) as for $Sp(p,q)$.
Consequently, there are $2^{t}$ different systems of positive roots satisfying condition (A), where $t$ is the number of $i$ such that $0<m_{i}=n_{i}$. If $\lambda$ satisfies condition (B) as in $\S$2.1, there is a unique system of positive roots $\Psi$ satisfying condition (A). As in $\S$2.1, we denote by $\pi(\lambda,\Psi)$ the limit of discrete series parameterized by $(\lambda,\Psi)$. Its unique lowest $K$-type is $\Lambda=\lambda+\rho_{n}-\rho_{c}$, where $\rho_{n}$ is one-half the sum of non-compact roots in $\Psi$.

Cuspidal parabolic subgroups of $O^{*}(2n)$ are of the form $P=MAN=LN$ with
\begin{equation}\label{Cuspidal Parabolic Subgroup3}
  L\cong O^{*}(2(n-2r))\times GL(1,\mathbb{H})^{r}
\end{equation}
and $0\leq r\leq \lfloor\frac{n}{2}\rfloor$. Set $R=\lfloor\frac{n}{2}\rfloor$. Then $(\mathbb{H}^{n},<,>)$ admits a polar decomposition
\begin{equation}\label{Polar Decomposition2}
  \mathbb{H}^{n}=V_{R}^{+}\oplus V^{0}\oplus V_{R}^{-},
\end{equation}
where $V_{R}^{+}=\mathrm{span}_{\mathbb{H}}\{v_{1},v_{2},\cdots,v_{R}\}$ and $V_{R}^{-}=\mathrm{span}_{\mathbb{H}}\{v_{1}^{'},v_{2}^{'},\cdots,v_{R}^{'}\}$ are two maximal isotropic subspaces in $\mathbb{H}^{n}$ and are dual to each other with respect to $<,>$. Here $v_{i}$, $v_{i}^{'}$ can be chosen such that $<v_{i},v_{j}^{'}>=\delta_{ij}$.
For every integer $r$ such that $0\leq r\leq R$, we define $V_{r}$, $W_{r}$ and $X_{r}$ as in (\ref{Isotropic Subspace1}), (\ref{Isotropic Subspace2}) and (\ref{Complementary Subspace}). Then
\begin{equation}
  \{0\}\subset V_{1}^{+}\subset V_{2}^{+}\subset\cdots \subset V_{r}^{+}
\end{equation}
is an isotropic flag in $V_{R}^{+}$. The stabilizer $P_{r}$ of this flag in $O^{*}(2n)$ is a cuspidal parabolic subgroup of $O^{*}(2n)$. Furthermore, the group $O^{*}(2(n-2r))\times GL(1,\mathbb{H})^{r}$ can be embedded into $P_{r}$. The group $O^{*}(2(n-2r))$ is regarded as the isometry group of $V^{R-r}$. The group $GL(1,\mathbb{H})^{r}$ is equipped with an right action $\circ$ on $\mathbb{H}^{n}$ which is left $\mathbb{H}$-linear as follows. For $h=(h_{1},\cdots,h_{r})\in GL(1,\mathbb{H})^{r}$, the action of $h$ on $\mathbb{H}^{n}$ is
\begin{equation}\label{Right Action}
\begin{split}
 v_{i}\circ h=(\overline{h_{i}})^{-1}\cdot v_{i}, \quad \forall 1\leq i\leq r,\\
 v_{i}^{'}\circ h=h_{i}\cdot v_{i}^{'}, \quad \forall 1\leq i\leq r,\\
 u\circ h=u, \quad \forall u\in V^{R-r},
\end{split}
\end{equation}
where $\cdot$ is the natural left multiplication on the left $\mathbb{H}$ row vector space $\mathbb{H}^{n}$.
It is easy to check that
\begin{equation*}
  (v_{i}\circ h,v_{i}^{'}\circ h)=(v_{i},v_{i}^{'})
\end{equation*}
for each $h\in GL(1,\mathbb{H})^{r}$ and each $1\leq i\leq r$. Then we embed $O^{*}(2(n-2r))\times GL(1,\mathbb{H})^{r}$ into the parabolic subgroup $P_{r}$. In general, let
\begin{equation}
  \{0\}\subset V_{r_{1}}^{+}\subset V_{r_{2}}^{+}\subseteq \cdots\subset V_{r_{s}}^{+}
\end{equation}
be an isotropic flag in $V$, where $1\leq r_{1}<r_{2}<\cdots <r_{s}\leq R$. The stabilizer $P_{\{r_{1},\cdots,r_{s}\}}$ of this flag is a parabolic subgroup of $O^{*}(2n)$. Set $r_{0}=0$ and $d_{i}=r_{i}-r_{i-1}$ for $1\leq i\leq s$. Its Levi subgroup is isomorphic to
\begin{equation}\label{Cuspidal Parabolic Subgroup4}
  O^{*}(2(n-2r_{s}))\times \Pi_{i=1}^{s}GL(d_{i},\mathbb{H}).
\end{equation}
Here the group $GL(d_{i},\mathbb{H})$ has a right action on the left row $\mathbb{H}$ space
\begin{equation*}
  \mathrm{span}_{\mathbb{H}}\{v_{r_{i-1}+1}^{'},v_{r_{i-1}+2}^{'},\cdots,v_{r_{i}}^{'}\} \end{equation*}
as the natural right multiplication and a right action on the left $\mathbb{H}$ space
\begin{equation*}
  \mathrm{span}_{\mathbb{H}}\{v_{r_{i-1}+1},v_{r_{i-1}+2},\cdots,v_{r_{i}}\} \end{equation*}
such that $<v_{k}\circ h,v_{l}^{'}\circ h>=<v_{k},v_{l}^{'}>$ for $r_{i-1}<k, l\leq r_{i}$ and $h\in GL(d_{i},\mathbb{H})$. Up to conjugation, each proper parabolic subgroup of $O^{*}(2n)$ is the stabilizer of an isotropic flag.

 Substituting $Sp(p-r,q-r)$ by $O^{*}(2(n-2r))$, we define $X(r,\lambda,\Psi,\mu,\nu)$, $\pi(r,\lambda,\Psi,\mu,\nu)$ and obtain a parametrization for $\mathcal{R}(O^{*}(2n))$ as in $\S$2.1. For $O^{*}(2n)$, the non-parity condition F-2 also amounts to the requirement that $\mu_{i}$ is odd if $\nu_{i}=0$. The infinitesimal character of $X(r,\lambda,\Psi,\mu,\nu)$ and $\pi(r,\lambda,\Psi,\mu,\nu)$ is similar to that for $Sp(p,q)$.

\subsection{The Lowest $K$-types}

Using the standard theory of \cite{KnVo95Cohomological} and \cite{Vogan81representations}, we compute the lowest $K$-types of irreducible admissible representations of the groups $Sp(p,q)$ and $O^{*}(2n)$ from the Langlands parameters. The algorithm can also be founded in $\S$2.3 of \cite{LPTZ03Correspondence} and $\S$3 of \cite{Paul2005Symplectic}.

Let $G=Sp(p,q)$ and let $K=Sp(p)\times Sp(q)$ be the maximal compact subgroup of $G$ as in (\ref{Maximal Compact Subgroup1}). Let $T$ and $\Delta_{c}^{+}$ be those in $\S$2.1. We identify a $K$-type $\sigma$ with its highest weight with respect to $\Delta_{c}^{+}$, also denoted by $\sigma$. By definition 5.1 of \cite{Vogan79AlgebraicStructure}, let $\|\cdot\|_{K}$ be the norm of $K$-types as follows: for a $K$-type $\sigma$, $\|\sigma\|_{K}=\langle \sigma+2\rho_{c},\sigma+2\rho_{c} \rangle$, where $\langle , \rangle$ is the Killing form and $\rho_{c}$ is that in $\S$2.1. A $K$-types $\sigma$ is of the form
\begin{equation*}
  (a_{1},a_{2},\cdots,a_{p};b_{1},b_{2},\cdots,b_{q})
\end{equation*}
with $a_{1}\geq a_{2}\geq \cdots\geq a_{p}\geq 0$ and $b_{1}\geq b_{2}\geq \cdots\geq b_{p}\geq 0$. The norm of $\sigma$ is
\begin{equation*}
  \|\sigma\|_{K}=\sum_{i=1}^{p}(a_{i}+2p+2-2i)^{2}+\sum_{j=1}^{q}(b_{j}+2q+2-2j)^{2}.
\end{equation*}
If a $K$-type $\sigma$ occurs with minimal norm in an admissible representation $\pi$ of $G$, we call $\sigma$ a lowest $K$-type of $\pi$.

Let $\pi(r,\lambda,\Psi,\mu,\nu)$ be an irreducible admissible representation of $Sp(p,q)$, with
\begin{equation*}
  \lambda=(\lambda_{1};\lambda_{2})=(\overbrace{a_{1},\cdots,a_{1}}^{m_{1}},\overbrace{a_{2},\cdots,a_{2}}^{m_{2}},\cdots,\overbrace{a_{k},\cdots,a_{k}}^{m_{k}};\overbrace{a_{1},\cdots,a_{1}}^{n_{1}},\cdots,\overbrace{a_{k},\cdots,a_{k}}^{n_{k}})
\end{equation*}
as in (\ref{Harish-Chandra Parameter1}). For $\pi(r,\lambda,\Psi,\mu,\nu)$, we write $\lambda^{'}$ as the dominant weight which is $W$-conjugate to $(\lambda_{1},\frac{\mu}{2};\lambda_{2},\frac{\mu}{2})$. Write
\begin{equation}\label{Twisted Harish-Chandra Parameter1}
  \lambda^{'}=(\overbrace{\alpha_{1},\cdots,\alpha_{1}}^{M_{1}},\overbrace{\alpha_{2},\cdots,\alpha_{2}}^{M_{2}},\cdots,\overbrace{\alpha_{s},\cdots,\alpha_{s}}^{M_{s}};\overbrace{\alpha_{1},\cdots,\alpha_{1}}^{N_{1}},\cdots,\overbrace{\alpha_{s},\cdots,\alpha_{s}}^{N_{s}})
\end{equation}
with $\alpha_{1}>o_{2}>\cdots >\alpha_{s}>0$. Notice that $\alpha_{i}\in \frac{1}{2}\mathbb{Z}$ and $|M_{i}-N_{i}|\leq 1$. If $M_{i}\neq N_{i}$, $\alpha_{i}$ is an integer.

Let $\mathfrak{q}=\mathfrak{l}\oplus \mathfrak{u}$ be the $\theta$-stable parabolic subalgebra of $\mathfrak{g}$ defined by $\lambda^{'}$ (see Definition 5.2.1 of \cite{Vogan81representations}). The lowest $K$-types of $\pi(r,\lambda,\Psi,\mu,\nu)$ are of the form
\begin{equation}\label{Lowest Type1}
  \Lambda=\lambda^{'}+\rho(\mathfrak{u}\cap \mathfrak{p})-\rho(\mathfrak{u}\cap \mathfrak{k})+\delta_{L}.
\end{equation}
Here $\rho(\mathfrak{u}\cap \mathfrak{p})$ and $\rho(\mathfrak{u}\cap \mathfrak{k})$ are one-half of the sums of roots associated to $\mathfrak{u}\cap \mathfrak{p}$ and $\mathfrak{u}\cap \mathfrak{k}$, respectively, and $\delta_{L}$ will be given explicitly below.

Set $\widetilde{M}_{i}=M_{1}+\cdots +M_{i}$ and $\widetilde{N}_{i}=N_{1}+\cdots +N_{i}$. Then $\widetilde{M}_{s}=p$, $\widetilde{N}_{s}=q$ and
we write $\lambda^{'}+\rho(\mathfrak{u}\cap \mathfrak{p})-\rho(\mathfrak{u}\cap \mathfrak{k})$ from (\ref{Lowest Type1}) as
\begin{equation}\label{Form of Lowest K-types1}
  (\overbrace{\beta_{1},\cdots,\beta_{1}}^{M_{1}},\cdots,\overbrace{\beta_{s},\cdots,\beta_{s}}^{M_{s}};\overbrace{\gamma_{1},\cdots,\gamma_{1}}^{N_{1}},\cdots,\overbrace{\gamma_{s},\cdots,\gamma_{s}}^{N_{s}}),
\end{equation}
where
\begin{equation}\label{Form of Entries1}
\begin{split}
 \beta_{i}=\alpha_{i}+\widetilde{M}_{i}-\widetilde{N}_{i}-\frac{1}{2}(M_{i}-N_{i}+1)+q-p,\\
 \gamma_{i}=\alpha_{i}-\widetilde{M}_{i}+\widetilde{N}_{i}-\frac{1}{2}(N_{i}-M_{i}+1)+p-q.
\end{split}
\end{equation}

\begin{thm}\label{Algorithm1}
Let $\pi(r,\lambda,\Psi,\mu,\nu)$ be an irreducible admissible representation of $Sp(p,q)$. Write $\lambda$ as in $(\ref{Harish-Chandra Parameter1})$ with $k$, $a_{i}$, $m_{i}$ and $n_{i}$ as defined there. Write $\lambda^{'}+\rho(\mathfrak{u}\cap \mathfrak{p})-\rho(\mathfrak{u}\cap \mathfrak{k})$ from $(\ref{Lowest Type1})$ as in $(\ref{Form of Lowest K-types1})$ and $(\ref{Form of Entries1})$. Let $\widetilde{m}_{j}=\sum _{c=1}^{j}m_{c}$ and let $\widetilde{n}_{j}=\sum _{c=1}^{j}n_{c}$ for $1\leq j\leq k$. Let $\widetilde{m}_{0}=\widetilde{n}_{0}=0$. The lowest $K$-types of  $\pi(r,\lambda,\Psi,\mu,\nu)$ are precisely those of the form
\begin{equation*}
  \lambda^{'}+\rho(\mathfrak{u}\cap \mathfrak{p})-\rho(\mathfrak{u}\cap \mathfrak{k})+\delta_{L}
\end{equation*}
with
\begin{equation*}
  \delta_{L}=(\overbrace{\delta_{1},\cdots,\delta_{1}}^{M_{1}},\cdots,\overbrace{\delta_{s},\cdots,\delta_{s}}^{M_{s}};\overbrace{-\delta_{1},\cdots,-\delta_{1}}^{N_{1}},\cdots,\overbrace{-\delta_{s},\cdots,-\delta_{s}}^{N_{s}})
\end{equation*}
satisfying the following condition:
\newline
$(1)$ If $\beta_{i}$ is an integer, then $\delta_{i}=0$.
\newline
$(2)$ Suppose $\beta_{i}\in \mathbb{Z}+\frac{1}{2}$, then $\delta_{i}=\frac{1}{2}$ or $-\frac{1}{2}$. If $\alpha_{i}$ does not occur as an entry of in $\lambda$, then both choices occur. If $\alpha_{i}=a_{j}$, then $\delta_{i}=\frac{1}{2}$ if $e_{\widetilde{m}_{j-1}+1}-f_{\widetilde{n}_{j-1}+1}\in \Psi$, and $\delta_{i}=-\frac{1}{2}$ otherwise.
\end{thm}

\begin{rmk}\label{Multiplicity}
By the definition of $\pi(r,\lambda,\Psi,\mu,\nu)$, its lowest $Sp(p)\times Sp(q)$-types coincide with those of the standard module $X(r,\lambda,\Psi,\mu,\nu)$. In fact, every lowest $Sp(p)\times Sp(q)$-type of the standard module occurs with multiplicity one (see $\S$3.2 of \cite{Paul1998Unitary}).
\end{rmk}

Let $G=O^{*}(2n)$ and let $K=U(n)$ be the maximal compact subgroup as in (\ref{Maximal Compact Subgroup2}). Let $T$, $\Delta_{c}^{+}$ and $\rho_{c}$ be those in $\S$2.2. Let $\mathfrak{t}_{0}$ be the real Lie algebra of $T$ as in $\S$2.2. For $O^{*}(2n)$, a $K$-type $\sigma$ is of the form
\begin{equation*}
  (a_{1},a_{2},\cdots,a_{n})
\end{equation*}
with $a_{1}\geq a_{2}\geq \cdots\geq a_{n}$. The norm of $\sigma$ is
\begin{equation*}
  \|\sigma\|_{K}=\langle \sigma+2\rho_{c},\sigma+2\rho_{c} \rangle=\sum_{i=1}^{n}(a_{i}+n+1-2i)^{2}.
\end{equation*}
We define the lowest $K$-types for admissible representations of $O^{*}(2n)$ as for those of $Sp(p,q)$.

Let $\pi(r,\lambda,\Psi,\mu,\nu)$ be an irreducible admissible representation of $O^{*}(2n)$, with
\begin{equation*}
  \lambda=(\overbrace{a_{1},\cdots,a_{1}}^{m_{1}},\overbrace{a_{2},\cdots,a_{2}}^{m_{2}},\cdots,\overbrace{a_{k},\cdots,a_{k}}^{m_{k}},\overbrace{-a_{k},\cdots,-a_{k}}^{n_{k}},\cdots,\overbrace{-a_{1},\cdots,-a_{1}}^{n_{1}});
\end{equation*}
as in (\ref{Harish-Chandra Parameter2}); or with one zero between the last $a_{k}$ and the first $-a_{k}$. We write $\lambda^{'}$ as the dominant weight which is $W$-conjugate to $(\lambda,\frac{\mu}{2},-\frac{\mu}{2})$. Write
\begin{equation}
  \lambda^{'}=(\overbrace{\alpha_{1},\cdots,\alpha_{1}}^{M_{1}},\overbrace{\alpha_{2},\cdots,\alpha_{2}}^{M_{2}},\cdots,\overbrace{\alpha_{s},\cdots,\alpha_{s}}^{M_{s}},\overbrace{-\alpha_{s},\cdots,-\alpha_{s}}^{N_{s}},\cdots,\overbrace{-\alpha_{1},\cdots,-\alpha_{1}}^{N_{1}})
\end{equation}
with $\alpha_{i}$ as in (\ref{Twisted Harish-Chandra Parameter1}); or with one zero between the last $\alpha_{s}$ and the first $-\alpha_{s}$.

Let $\mathfrak{q}=\mathfrak{l}\oplus \mathfrak{u}$ be the $\theta$-stable parabolic subalgebra of $\mathfrak{g}$ defined by $\lambda^{'}$. The lowest $K$-types of $\pi(r,\lambda,\Psi,\mu,\nu)$ are of the form
\begin{equation}\label{Lowest Type2}
  \Lambda=\lambda^{'}+\rho(\mathfrak{u}\cap \mathfrak{p})-\rho(\mathfrak{u}\cap \mathfrak{k})+\delta_{L}.
\end{equation}
Here $\rho(\mathfrak{u}\cap \mathfrak{p})$ and $\rho(\mathfrak{u}\cap \mathfrak{k})$ are one-half of the sums of roots associated to $\mathfrak{u}\cap \mathfrak{p}$ and $\mathfrak{u}\cap \mathfrak{k}$, respectively, and $\delta_{L}$ will be given explicitly below.

Set $\widetilde{M}_{i}=M_{1}+\cdots +M_{i}$ and $\widetilde{N}_{i}=N_{1}+\cdots +N_{i}$. Then $\widetilde{M}_{s}+\widetilde{N}_{s}=n$ or $n-1$ and we write $\lambda^{'}+\rho(\mathfrak{u}\cap \mathfrak{p})-\rho(\mathfrak{u}\cap \mathfrak{k})$ from (\ref{Lowest Type2}) as
\begin{equation}\label{Form of Lowest K-types2}
  \gamma^{'}=(\overbrace{\beta_{1},\cdots,\beta_{1}}^{M_{1}},\cdots,\overbrace{\beta_{s},\cdots,\beta_{s}}^{M_{s}},\overbrace{-\gamma_{s},\cdots,-\gamma_{s}}^{N_{s}},\cdots,\overbrace{-\gamma_{1},\cdots,-\gamma_{1}}^{N_{1}})
\end{equation}
if $\widetilde{M}_{s}+\widetilde{N}_{s}=n$; or with $\widetilde{M}_{s}-\widetilde{N}_{s}$ between the last $\beta_{s}$ and the first $-\gamma_{s}$ if $\widetilde{M}_{s}+\widetilde{N}_{s}=n-1$. Here
\begin{equation}\label{Form of Entries2}
\begin{split}
  \beta_{i}=\alpha_{i}+\widetilde{M}_{i}-\widetilde{N}_{i}-\frac{1}{2}(M_{i}-N_{i}+1),\\
  \gamma_{i}=\alpha_{i}-\widetilde{M}_{i}+\widetilde{N}_{i}-\frac{1}{2}(N_{i}-M_{i}+1).
\end{split}
\end{equation}

\begin{thm}\label{Algorithm2}
 Let $\pi(r,\lambda,\Psi,\mu,\nu)$ be an irreducible admissible representation of $O^{*}(2n)$. Write $\lambda$ as in $(\ref{Harish-Chandra Parameter2})$ with $k$, $a_{i}$, $m_{i}$ and $n_{i}$ as defined there. Write $\lambda^{'}+\rho(\mathfrak{u}\cap \mathfrak{p})-\rho(\mathfrak{u}\cap \mathfrak{k})$ from $(\ref{Lowest Type2})$ as in $(\ref{Form of Lowest K-types2})$ and $(\ref{Form of Entries2})$. Let $\widetilde{m}_{j}=\sum _{c=1}^{j} m_{j}$ and let $\widetilde{n}_{j}=n-2r-\sum _{c=1}^{j} n_{j}$ for $1\leq j\leq k$. Let $\widetilde{m}_{0}=0$ and let $\widetilde{n}_{0}=n-2r$. The lowest $K$-types of  $\pi(r,\lambda,\Psi,\mu,\nu)$ are precisely those of the form
\begin{equation*}
  \Lambda=\gamma+\rho(\mathfrak{u}\cap \mathfrak{p})-\rho(\mathfrak{u}\cap \mathfrak{k})+\delta_{L}.
\end{equation*}
with
\begin{equation*}
  \delta_{L}=(\overbrace{\delta_{1},\cdots,\delta_{1}}^{M_{1}},\cdots,\overbrace{\delta_{s},\cdots,\delta_{s}}^{M_{s}},\overbrace{\delta_{s},\cdots,\delta_{s}}^{N_{s}},\cdots,\overbrace{\delta_{1},\cdots,\delta_{1}}^{N_{1}})
\end{equation*}
if $\widetilde{M}_{s}+\widetilde{N}_{s}=n$; or with one zero as the $(\widetilde{M}_{s}+1)$-th entry if $\widetilde{M}_{s}+\widetilde{N}_{s}=n-1$. Here $\delta_{L}$ is any element in $\texttt{i}\mathfrak{t}_{0}^{*}$ satisfying the following condition:
\newline
$(1)$ If $\beta_{i}$ is an integer, then $\delta_{i}=0$.
\newline
$(2)$ Suppose $\beta_{i}\in \mathbb{Z}+\frac{1}{2}$, then $\delta_{i}=\frac{1}{2}$ or $-\frac{1}{2}$. If $\alpha_{i}$ does not occur as an entry of in $\lambda$, then both choices occur. If $\alpha_{i}=a_{j}$, then $\delta_{i}=\frac{1}{2}$ if $e_{\widetilde{m}_{j-1}+1}+e_{\widetilde{n}_{j}+1}\in \Psi$, and $\delta_{i}=-\frac{1}{2}$ otherwise.
\end{thm}

For $O^{*}(2n)$, we have the similar conclusions as in Remark \ref{Multiplicity}.

\subsection{Admissible Dual of $O^{*}(4)$}

The classification of the admissible dual of the groups $O^{*}(2n)$ has been given in $\S$2.2. In this section, we will give a more accurate description of the admissible dual of $O^{*}(4)$. As in (\ref{Real Lie Algebra2}), the real Lie algebra $\mathfrak{o}^{*}(4)$ is
\begin{equation*}
  \left\{\left(
            \begin{array}{cc}
              A & B \\
              -\overline{B} & \overline{A} \\
            \end{array}
          \right)
  |A+A^{*}=0,B+B^{t}=0\right\}
\end{equation*}
with $A,B\in M_{2\times 2}(\mathbb{C})$. Let $\varphi_{1}$ be a injective homomorphism from $\mathfrak{su}(2)$ into $\mathfrak{o}^{*}(4)$:
\begin{equation*}
  A\mapsto \left(\begin{array}{cc}
              A & 0 \\
              0 & \overline{A}
            \end{array}
  \right).
\end{equation*}
Let $\varphi_{2}$ be a injective homomorphism from $\mathfrak{sl}(2,\mathbb{R})$ into $\mathfrak{o}^{*}(4)$:
\begin{equation*}
\begin{split}
\left(
     \begin{array}{cc}
       1 & 0 \\
       0 & -1 \\
     \end{array}
   \right)
  \mapsto \left(
            \begin{array}{cc}
              0 & H \\
              -H & 0 \\
            \end{array}
          \right)
  , H=\left(
          \begin{array}{cc}
            0 & 1 \\
            -1 & 0 \\
          \end{array}
        \right)
  ,\\
  \left(
     \begin{array}{cc}
       0 & 1 \\
       1 & 0 \\
     \end{array}
   \right)
  \mapsto \left(
            \begin{array}{cc}
              0 & X \\
              X & 0 \\
            \end{array}
          \right)
  , X=\left(
          \begin{array}{cc}
            0 & -\textbf{i} \\
            \textbf{i} & 0 \\
          \end{array}
        \right),\\
 \left(
     \begin{array}{cc}
       0 & 1 \\
       -1 & 0 \\
     \end{array}
   \right)
  \mapsto \left(
            \begin{array}{cc}
              Y & 0 \\
              0 & -Y \\
            \end{array}
          \right)
  , Y=\left(
          \begin{array}{cc}
            \textbf{i} & 0 \\
            0 & \textbf{i} \\
          \end{array}
        \right).
\end{split}
\end{equation*}
Then $\mathfrak{o}^{*}(4)$ is isomorphic to $\mathfrak{su}(2)\oplus \mathfrak{sl}(2,\mathbb{R})$. As in (\ref{Maximal Compact Subgroup2}), $U(2)$ can be regarded as a maximal compact subgroup of $O^{*}(4)$. Then the homomorphism $\varphi_{1}$ can be lifted to an embedding $\Phi_{1}$ of $SU(2)$ into $O^{*}(4)$ and the homomorphism $\varphi_{2}$ can be lifted to an embedding $\Phi_{2}$ of $SL(2,\mathbb{R})$ into $O^{*}(4)$ (The homomorphism  $\varphi_{2}$ restricted on the maximal compact subalgebra $\mathfrak{so}(2)$ of $\mathfrak{sl}(2,\mathbb{R})$ can be lifted to an isomorphism from $SO(2)$ to $U(1)$.).  So we have a projection map $\Phi=\Phi_{1}\times \Phi_{2}$ from $SU(2)\times SL(2,\mathbb{R})$ onto $O^{*}(4)$. As in $\S$5 of \cite{LTZ01TensorProduct}, we identify $G=O^{*}(4)$ with $(SU(2)\times SL(2,\mathbb{R}))/\{\pm I\}$. To understand the representations of $O^{*}(4)$, it is necessary for us to understand the representations of $SU(2)$ and $SL(2,\mathbb{R})$.

Let $G_{1}=SU(2)$ and let $G_{2}=SL(2,\mathbb{R})$. All irreducible admissible representations of $SU(2)$ are of the form $\tau_{k}$, where $\tau_{k}$ is the $k+1$-dimensional irreducible representation of $SU(2)$. Then we concentrate on the admissible dual of $SL(2,\mathbb{R})$. The classification can be founded in \cite{Knapp86Overview} and \cite{Lang86} and we state the results in this section.

 The unique proper parabolic subgroup of $G_{2}$ up to conjugation is the upper triangular group $P_{2}=M_{2}A_{2}N_{2}$, where $M_{2}=\{\pm I\}$,
\begin{equation*}
  A_{2}=\{\left(
                           \begin{array}{cc}
                             a & 0 \\
                             0 & a^{-1} \\
                           \end{array}
                         \right)
|a>0 \},
\end{equation*}
and
\begin{equation*}
  N_{2}=\{\left(\begin{array}{cc}
                                                                    1 & n \\
                                                                    0 & 1
                                                                  \end{array}\right)
|n\in \mathbb{R}\}.
\end{equation*}
For $k\in \mathbb{Z}_{\geq 0}$ and $\nu\in \mathbb{C}$, we denote by $I^{k}(\nu)$ the induced representation $Ind_{P_{2}}^{SL(2,\mathbb{R})}(\mathrm{sgn}^{k}\otimes \nu\otimes 1{\kern-0.36em} 1)$. For $\nu\in \mathbb{Z}_{\geq 1}$, we denote by $D_{\nu}$ the irreducible lowest weight $SL(2,\mathbb{R})$-module of lowest weight $\nu$, and by $\overline{D}_{-\nu}$ the irreducible highest weight $SL(2,\mathbb{R})$-module of highest weight weight $-\nu$. For $\nu\in \mathbb{Z}_{\geq 0}$, we denote by $F_{\nu}$ the irreducible finite dimensional representation of $SL(2,\mathbb{R})$ of highest weight $\nu$ (or dimension $\nu+1$). We have the following theorems (see Chapter \textrm{VI} of \cite{Lang86}).

\begin{thm}\label{Admissible Dual1}
$(1)$ If $k\not \equiv \nu-1~(\mathrm{mod}~2)$, then $I^{k}(\nu)$ is an irreducible representation. Furthermore, $I^{k}(\nu)$ is infinitesimal equivalent to $I^{k}(-\nu)$.
\newline
$(2)$ If $k\equiv\nu-1~(\mathrm{mod}~2)$ and $\nu>0$, then $I^{k}(\nu)$ has two irreducible submodule: $D_{\nu+1}$ and $\overline{D}_{-(\nu+1)}$. The unique irreducible quotient of $I^{k}(\nu)$ is $F_{\nu-1}$.
\newline
$(3)$ If $k\equiv\nu-1~(\mathrm{mod}~2)$ and $\nu>0$, then $I^{k}(-\nu)$ has a unique submodule $F_{\nu-1}$, and its quotient by $F_{\nu-1}$ is isomorphic to $D_{\nu+1}\oplus \overline{D}_{-(\nu+1)}$.
\newline
$(4)$ If $k$ is odd, then $I^{k}(0)$ is isomorphic to $D_{1}\oplus \overline{D}_{-1}$.
\end{thm}

\begin{thm}
The following is a list of infinitesimal equivalence classes of irreducible admissible representations of $SL(2,\mathbb{R})$:
\begin{enumerate}[(1)]
  \item Irreducible principal series $I^{k}(\nu)$, $\mathrm{Re}\nu\geq 0$;
  \item Lowest weight modules $D_{\nu+1}$, $\nu\in \mathbb{Z}_{\geq 0}$;
  \item Highest weight modules $\overline{D}_{-(\nu+1)}$, $\nu\in \mathbb{Z}_{\geq 0}$;
  \item Finite dimensional representations $F_{\nu-1}$, $\nu\in \mathbb{Z}_{\geq 1}$.
\end{enumerate}
\end{thm}

Let $P=(G_{1}\times P_{2})/\{\pm I\}$. It is not difficult to check that $P$ is the stabilizer of the isotropic row vector space $\mathbb{H}(1,-\textbf{j})$ in $O^{*}(4)$ as in $\S$2.2, where $\mathbb{H}^{2}=\mathbb{H}(1,-\textbf{j})\oplus \mathbb{H}(\frac{\textbf{i}}{2},\frac{\textbf{k}}{2})$ such that $<(1,-\textbf{j}),(\frac{\textbf{i}}{2},\frac{\textbf{k}}{2})>=1$. For $k\in \mathbb{Z}_{\geq 0}$ and $\nu\in \mathbb{C}$, it is not difficult to see that the representation $\tau_{k}\otimes I^{k}(\nu)$ is infinitesimal equivalent to the standard module $X(1,0,\emptyset,k+1,\nu)$ as a representation of $SU(2)\times SL(2,\mathbb{R})$ (also see $\S$5 of \cite{LTZ01TensorProduct}). Furthermore, if $\nu$ is a non-negative integer such that $\nu\equiv k+1~(\mathrm{mod}~2)$, $\tau_{k}\otimes D_{\nu+1}$ is infinitesimal equivalent to the limit of discrete series $\pi((\frac{\nu+k+1}{2},\frac{\nu-k-1}{2}),\{e_{1}\pm e_{2}\})$ while $\tau_{k}\otimes \overline{D}_{-(\nu+1)}$ is infinitesimal equivalent to the limit of discrete series $\pi((\frac{-\nu+k+1}{2},\frac{-\nu-k-1}{2}),\{e_{2}\pm e_{1}\})$ (as representations of $SU(2)\times SL(2,\mathbb{R})$). If $\nu\in \mathbb{Z}_{\geq 1}$, by Theorem \ref{Admissible Dual1}, it is not difficult to see that the unique irreducible quotient of $X(1,0,\emptyset,k+1,\nu)$ is infinitesimal equivalent to $\tau_{k}\otimes F_{\nu-1}$ as a representation of $SU(2)\times SL(2,\mathbb{R})$, and this quotient is infinitesimal equivalent to $\pi(1,0,\emptyset,k+1,\nu)$ as a representation of $O^{*}(4)$. The following is a list of infinitesimal equivalence classes of irreducible admissible representations of $O^{*}(4)$:
\begin{enumerate}[(1)]
  \item Irreducible principal series $\pi(1,0,\emptyset,k+1,\nu)$, $\mathrm{Re}\nu\geq 0$, $\nu\not \equiv k+1~(\mathrm{mod}~2)$;
  \item Limits of lowest weight discrete series $\pi((\frac{\nu+k+1}{2},\frac{\nu-k-1}{2}),\Psi)$ with $\Psi=\{e_{1}-e_{2},e_{1}+e_{2}\}$, $\nu\in \mathrm{Z}_{\geq 0}$ and $\nu\equiv k+1~(\mathrm{mod}~2)$;
  \item Limits of highest weight discrete series $\pi((\frac{-\nu+k+1}{2},\frac{-\nu-k-1}{2}),\Psi)$ with $\Psi=\{e_{1}-e_{2},-(e_{1}+e_{2})\}$, $\nu\in \mathrm{Z}_{\geq 0}$ and $\nu\equiv k+1~(\mathrm{mod}~2)$;
  \item Finite dimensional representations $\pi(1,0,\emptyset,k+1,\nu)$, $\nu\in \mathrm{Z}_{\geq 1}$, $\nu\equiv k+1~(\mathrm{mod}~2)$.
\end{enumerate}

It is easy to check the infinitesimal equivalence classes in each case are uniquely determined by their infinitesimal characters $(\lambda_{1},\lambda_{2})=(\frac{\nu+k+1}{2},\frac{\nu-k-1}{2})$. By results in Chapter $\textrm{VI}$ of \cite{Lang86}, it is not difficult to determine their structures of $U(2)$-types. We have the following theorem.

\begin{thm}\label{Admissible Dual2}
The following is a list of infinitesimal equivalence classes of irreducible admissible representations of $O^{*}(4)$:
\begin{enumerate}[(1)]
  \item Irreducible principal series $P_{\lambda_{1},\lambda_{2}}=\pi(1,0,\emptyset,\lambda_{1}-\lambda_{2},\lambda_{1}+\lambda_{2})$ with $\lambda_{1},\lambda_{2}\notin \mathbb{Z}$, $\mathrm{Re}(\lambda_{1}+\lambda_{2})\geq0$ and $\lambda_{1}-\lambda_{2}\in \mathbb{Z}_{\geq 1}$. Here $\mathcal{K}(P_{\lambda_{1},\lambda_{2}})=\{(p+\lambda_{1}-\lambda_{2}-1,p)|\forall p\in \mathbb{Z}\}$;
  \item Limits of lowest weight discrete series $D_{\lambda_{1},\lambda_{2}}=\pi((\lambda_{1},\lambda_{2}),\Psi)$ with $\Psi=\{e_{1}-e_{2},e_{1}+e_{2}\}$.  Here $\mathcal{K}(D_{\lambda_{1},\lambda_{2}})=\{(p+\lambda_{1}-\lambda_{2}-1,p)|\forall p\geq \lambda_{2}+1\}$;
  \item Limits of highest weight discrete series $\overline{D}_{\lambda_{1},\lambda_{2}}=\pi((\lambda_{1},\lambda_{2}),\Psi)$ with $\Psi=\{e_{1}-e_{2},-(e_{1}+e_{2})\}$.  Here $\mathcal{K}(\overline{D}_{\lambda_{1},\lambda_{2}})=\{(p+\lambda_{1}-\lambda_{2}-1,p)|\forall p\leq \lambda_{2}\}$;
  \item Finite dimensional representations $F_{\lambda_{1},\lambda_{2}}=\pi(1,0,\emptyset,\lambda_{1}-\lambda_{2},\lambda_{1}+\lambda_{2})$ with $\lambda_{1},\lambda_{2}\in \mathbb{Z}$, $\lambda_{1}+\lambda_{2}\in \mathbb{Z}_{\geq 1}$ and $\lambda_{1}-\lambda_{2}\in \mathbb{Z}_{\geq 1}$. Here $\mathcal{K}(F_{\lambda_{1},\lambda_{2}})=\{(p+\lambda_{1}-\lambda_{2}-1,p)| 1-\lambda_{1}\leq p\leq \lambda_{2}\}$.
\end{enumerate}
Furthermore, for each $\pi\in \mathcal{R}(O^{*}(4))$, the multiplicity $m(\sigma,\pi)$ of each $\sigma\in \mathcal{K}(\pi)$ equals to one.
\end{thm}

\section{Local Theta Correspondence}

\newtheorem{lem}[thm]{Lemma}
\newproof{pf}{Proof}

For non-negative integers $p,q,n$ such that $p+q\geq 1$ and $n\geq 1$, let $V=\mathbb{H}^{p+q}$ be a right $\mathbb{H}$ column vector space with Hermitian form (\ref{Hermitian form}) and let $V^{'}=\mathbb{H}^{n}$ be a left $\mathbb{H}$ row vector space with skew-Hermitian form (\ref{Skew-Hermitian Form}). Let $W=\mathbb{H}^{p+q}\otimes _{\mathbb{H}}\mathbb{H}^{n}$ be a $4n(p+q)$-dimensional real space with symplectic form $\langle, \rangle=\mathrm{tr}_{\mathbb{H}/\mathbb{R}}((,)\otimes \overline{<,>})$. As in $\S$2, the group $Sp(p,q)$ has a left action on $\mathbb{H}^{p+q}$ and the group $O^{*}(2n)$ has a right action on $\mathbb{H}^{n}$. Then the group $Sp(p,q)\times O^{*}(2n)$ has a left action on $W$ as follows. For $g\in Sp(p,q)$, $h\in O^{*}(2n)$,$v\in \mathbb{H}^{p+q}$ and $w\in \mathbb{H}^{n}$,
\begin{equation}\label{Bilateral Action}
  (g,h)\cdot (v\otimes w)=(g\cdot v)\otimes (w\cdot h^{-1}).
\end{equation}
Since this action keeps the symplectic form $\langle, \rangle$ invariant,  $(Sp(p,q),O^{*}(2n))$ is embedded into $Sp(4n(p+q),\mathbb{R})$. Furthermore, $Sp(p,q)$ and $O^{*}(2n)$ are centralizer of each other (as algebraic groups) in the symplectic group $Sp=Sp(4n(p+q),\mathbb{R})$. We call $(Sp(p,q),O^{*}(2n))$ a reductive dual pair in $Sp$. We endow $W$ with a complex Hermitian form as follows. Let $V_{1}$ be the subspace of vectors of $V$ whose last $q$ entries are zero and let $V_{2}$ be the subspace of vectors of $V$ whose first $p$ entries are zero. We define a complex multiplication $\ast$ on $W$ as follows:
\begin{equation}
\begin{split}
 c\ast (v_{1}\otimes v^{'})=v_{1}\otimes v^{'}\cdot c,~ \forall c\in \mathbb{C}, ~\forall v_{1}\in V_{1}, ~\forall v^{'}\in V^{'},\\
 c\ast (v_{2}\otimes v^{'})=v_{2}\otimes v^{'}\cdot \overline{c},~ \forall c\in \mathbb{C}, ~\forall v_{2}\in V_{2}, ~\forall v^{'}\in V^{'}.
\end{split}
\end{equation}
Here $\cdot$ is the natural right multiplication on $\mathbb{H}^{n}$. We define a norm $\|\cdot \|$ on $W$: $\|w\|^{2}=\langle w,\textbf{i}\ast w \rangle$. The isometry group $U$ of this norm is a maximal compact subgroup of $Sp$ and is isomorphic to $U(2n(p+q))$. By the action (\ref{Bilateral Action}), the maximal compact group $K=Sp(p)\times Sp(q)$ of $G=Sp(p,q)$ as in (\ref{Maximal Compact Subgroup1}) and the maximal compact subgroup $K^{'}=U(n)$ of $G^{'}=O^{*}(2n)$ as in (\ref{Maximal Compact Subgroup2}) act on $W$. It is easy to check that the actions of $K$ and $K^{'}$ on $W$ commute with the complex multiplication and keep the norm invariant. Then $K$ and $K^{'}$ can be regarded as compact subgroups of $U$.

For a fixed additive character $\psi(t)=e^{2\pi\textbf{i}t}$ of $\mathbb{R}$ and the two-fold metaplectic cover group $\widetilde{Sp}$ of $Sp$, a unitary representation of $\widetilde{Sp}$ associated to $\psi$ was constructed by Shale \cite{Shale62Symmetries} and Weil \cite{Weil64SurCertains}. We call this representation the oscillator representation.  We denote by $\omega_{p,q,n}$ the oscillator representation for the reductive dual pair $(Sp(p,q),O^{*}(2n))$. Let $\widetilde{U}$ denote the inverse image of $U$ in the metaplectic cover and let $\mathscr{F}_{p,q,n}$ denote the Fock space (the Harish-Chandra module) of $\omega_{p,q,n}$. Since the two-fold metaplectic groups of both $Sp(p,q)$ and $O^{*}(2n)$ uniquely split, we can state the results in terms of the groups themselves instead of their covering groups. We denote by $\mathcal{R}(Sp(p,q),\omega_{p,q,n})$ (resp. $\mathcal{R}(O^{*}(2n),\omega_{p,q,n})$) the set of elements in $\mathcal{R}(Sp(p,q))$ (resp. $\mathcal{R}(O^{*}(2n))$) which can be realized as quotients of $\mathscr{F}_{p,q,n}$. We define the set $\mathcal{R}(Sp(p,q)\times O^{*}(2n),\omega_{p,q,n})$ in the same way. Howe \cite{Howe89Transcend} proved the following theorem.

\begin{thm}
The set $\mathcal{R}(Sp(p,q)\times O^{*}(2n),\omega_{p,q,n})$ defines a bijection between $\mathcal{R}(Sp(p,q),\omega_{p,q,n})$ and $\mathcal{R}(O^{*}(2n),\omega_{p,q,n})$.For $\pi\in \mathcal{R}(Sp(p,q),\omega_{p,q,n})$ and $\pi^{'}\in \mathcal{R}(O^{*}(2n),\omega_{p,q,n})$, they correspond to each other in this bijection if and only if $\pi\otimes \pi^{'}\in \mathcal{R}(Sp(p,q)\times O^{*}(2n),\omega_{p,q,n})$.
\end{thm}

We call this bijection the local theta correspondence or the Howe correspondence. We denote by $\theta_{p,q}$ the bijection from the set $\mathcal{R}(O^{*}(2n),\omega_{p,q,n})$ to the set $\mathcal{R}(Sp(p,q),\omega_{p,q,n})$ and by $\theta_{n}$ the bijection in the other direction. We say $\theta_{p,q}(\pi^{'})=0$ for $\pi^{'}\in \mathcal{R}(O^{*}(2n))$ if $\pi^{'}\notin \mathcal{R}(O^{*}(2n),\omega_{p,q,n})$. Similarly, we say $\theta_{n}(\pi)=0$ for $\pi\in \mathcal{R}(O^{*}(2n))$ if $\pi\notin \mathcal{R}(Sp(p,q),\omega_{p,q,n})$. Formally, we also define the theta correspondence for reductive dual pairs $(Sp(p,q),O^{*}(0))$ and $(Sp(0,0),O^{*}(2n))$ while only trivial representations occur in the correspondence.

In \cite{CP66Quantum}, the Fock space $\mathscr{F}_{p,q,n}$ can be realized on a complex polynomial space with $2n(p+q)$ variables  and the action of the two-fold metaplectic cover $\widetilde{U}$ keeps the degree of polynomials invariant. This allow us to associate to each $Sp(p)\times Sp(q)$-type $\sigma$ (resp. $U(n)$-type $\sigma^{'}$) occurring in $\mathscr{F}_{p,q,n}$ a degree which is the minimal degree of polynomials in the $\sigma$-isotypic subspace (resp.$\sigma^{'}$-isotypic subspace) of $\mathscr{F}_{p,q,n}$. For $\pi\in \mathcal{R}(Sp(p,q),\omega_{p,q,n})$, a $Sp(p)\times Sp(q)$-type $\sigma$ is called a lowest degree $Sp(p)\times Sp(q)$-type of $\pi$ if and only if it is of minimal degree among all $Sp(p)\times Sp(q)$-types occurring in $\pi$. Similarly, we define the lowest degree $U(n)$-types for $\pi^{'}\in \mathcal{R}(O^{*}(2n),\omega_{p,q,n})$. Howe \cite{Howe89Transcend} proved that there exists a $(Sp(p)\times Sp(q))\times U(n)$-invariant subspace $\mathscr{H}_{p,q,n}$ of $\mathscr{F}_{p,q,n}$ called the space of joint harmonics with the following properties.

\begin{thm}\label{Joint Harmonics}
There is a one-one correspondence of $Sp(p)\times Sp(q)$-types and $U(n)$-types on the joint harmonics $\mathscr{H}_{p,q,n}$ with the following properties. Suppose $\pi\in \mathcal{R}(Sp(p,q),\omega_{p,q,n})$ and $\pi^{'}\in \mathcal{R}(O^{*}(2n),\omega_{p,q,n})$ correspond to each other in the theta correspondence for $(Sp(p,q),O^{*}(2n))$. Let $\sigma$ be a lowest degree $Sp(p)\times Sp(q)$-type of $\pi$. Then $\sigma$ occurs in the joint harmonics $\mathscr{H}_{p,q,n}$ and corresponds to a $U(n)$-type $\sigma^{'}$. Furthermore, $\sigma^{'}$ is a lowest degree $U(n)$-type of $\pi^{'}$. The statement of this theorem is also true with the roles of $\pi$ and $\pi^{'}$ (also of $\sigma$ and $\sigma^{'}$) reversed.
\end{thm}

Let $\sigma$ be a $Sp(p)\times Sp(q)$-type and let $\sigma^{'}$ be a $U(n)$-type. Write
\begin{equation}\label{Form of K-types}
\begin{split}
  \left\{
     \begin{array}{ll}
       \sigma=(a_{1},a_{2},\cdots,a_{r},0,\cdots,0;b_{1},b_{2},\cdots,b_{s},0,\cdots,0) & \hbox{$r\leq p, s\leq q$,} \\
       \sigma^{'}=(a_{1}^{'},a_{2}^{'},\cdots,a_{r^{'}}^{'},0,\cdots,0,-b_{s^{'}}^{'},\cdots,-b_{2}^{'},-b_{1}^{'})+
       \\
       \hspace{0.8cm} (p-q,\cdots,p-q) & \hbox{$r^{'}+s^{'}\leq n$,}
     \end{array}
   \right.
\end{split}
\end{equation}
where $a_{1}\geq \cdots\geq a_{r}>0$, $b_{1}\geq \cdots\geq b_{s}>0$, and similarly for the $a_{i}^{'}$ and $b_{j}^{'}$. We have the following lemma (Lemma 3.4 of \cite{LPTZ03Correspondence}).

\begin{lem}\label{Correspondence in Joint Harmonics}
$(1)$ Let $\sigma$ be a $Sp(p)\times Sp(q)$-type as in $(\ref{Form of K-types})$. Then $\sigma$ occurs in $\mathscr{F}_{p,q,n}$ if and only if $r,s\leq n$. Let $\sigma^{'}$ be a $U(n)$-type as in $(\ref{Form of K-types})$. Then $\sigma^{'}$ occurs in $\mathscr{F}_{p,q,n}$ if and only if $r^{'}\leq 2p, s^{'}\leq 2q$. If the conditions are satisfied, then
\begin{equation}\label{Degree Formula1}
  \mathrm{deg}(\sigma)=a_{1}+\cdots+a_{r}+b_{1}+\cdots+b_{s},
\end{equation}
\begin{equation}\label{Degree Formula2}
  \mathrm{deg}(\sigma^{'})=a_{1}^{'}+\cdots+a_{r^{'}}^{'}+b_{1}^{'}+\cdots+b_{s^{'}}^{'}.
\end{equation}
$(2)$ Let $\sigma$ be a $Sp(p)\times Sp(q)$-type as in $(\ref{Form of K-types})$. Then $\sigma$ occurs in the joint harmonics if and only if $r+s\leq n$. Let $\sigma^{'}$ be a $U(n)$-type as in $(\ref{Form of K-types})$. Then $\sigma^{'}$ occurs in the joint harmonics if and only if $r^{'}\leq p, s^{'}\leq q$. Furthermore, $\sigma$ corresponds to $\sigma^{'}$ in the joint harmonics if and only if $r=r^{'}, s=s^{'}$, $a_{i}=a_{i}^{'}$ and $b_{j}=b_{j}^{'}$.
\end{lem}

Let $H$ be a real reductive group with a maximal compact subgroup $K_{H}$. For an admissible representation $\pi$ of $H$, we denote by $\pi^{\vee}$ its contragredient representation. We have the following lemma (see Lemma 3.32 of \cite{LPTZ03Correspondence}).

\begin{lem}\label{Contragredient Representation}
Let $\pi^{'}$ be an irreducible admissible representation of $O^{*}(2n)$. Then $(\theta_{p,q}(\pi^{'}))^{\vee}=\theta_{q,p}((\pi^{'})^{\vee})$.
\end{lem}

Notice that if $\pi$ is an irreducible admissible representation of $Sp(p,q)$, $\pi=\pi^{\vee}$. By Lemma \ref{Contragredient Representation}, we reduce the problem to calculating $\theta_{p,q}$ ($p\geq q$) explicitly. The following theorem is on the correspondence of infinitesimal characters (see Theorem 1.19 of \cite{Pre96Infinitesimal}).

\begin{thm}\label{Correspondence of Infinitesimal Characters}
Let $p,q,n$ be non-negative integers such that $p+q\geq n$. Suppose $\pi\leftrightarrow \pi^{'}$ for the reductive dual pair $(Sp(p,q),O^{*}(2n))$. If the infinitesimal character of $\pi^{'}$ is $(\lambda_{1}^{'},\lambda_{2}^{'},\cdots,\lambda_{n}^{'})$, then the infinitesimal character of $\pi$ is $(\lambda_{1}^{'},\lambda_{2}^{'},\cdots,\lambda_{n}^{'},p+q-n,\cdots,2,1)$.
\end{thm}

To reduce the problem of explicit calculation, we introduce the stable range theorem. In \cite{PrP08FirstOcuurrence}, the authors proved the stable range theorem.

\begin{thm}\label{Stable Ranger Theorem}
Let $\pi^{'}$ be an irreducible admissible representation of $O^{*}(2n)$. Then $\theta_{p,q}(\pi^{'})$ is nonzero if $\min\{p,q\}\geq n$.
\end{thm}

For $\pi\in \mathcal{R}(Sp(p,q),\omega_{p,q,n})$, we denote by $\mathcal{A}_{p,q,n}(\pi)$ the set of lowest degree $Sp(p)\times Sp(q)$-types occurring in $\pi$ for the local theta correspondence for $(Sp(p,q),O^{*}(2n))$. Similarly, for $\pi^{'}\in \mathcal{R}(O^{*}(2n),\omega_{p,q,n})$, we also define the set $\mathcal{A}_{p,q,n}(\pi^{'})$. We state the following lemma without proof (see Lemma 5.1.1 of \cite{Paul1998Unitary}).

\begin{lem}\label{Form of Lowest Degree K-types}
Let $X(r,\lambda,\Psi,\mu,\nu)$ be the standard module of $Sp(p,q)$ as in $(\ref{Standard Module})$ and let $\sigma$ be a $Sp(p)\times Sp(q)$-type occurring in $X(r,\lambda,\Psi,\mu,\nu)$. Then $\sigma$ is of the form
\begin{equation*}
  \sigma=\Lambda+\Sigma n_{\alpha}\alpha,
\end{equation*}
where $\Lambda$ is a lowest $K$-type of $X(P,r,\lambda,\Psi,\mu,\nu)$, the sum runs over roots in $\Delta(\mathfrak{l}:\mathfrak{t})\cup \Delta(\mathfrak{u}\cap \mathfrak{p})$ and $n_{\alpha}\geq 0$ for all $\alpha$. Here $\mathfrak{l}$ and $\mathfrak{u}$ are those in $\S2.3$.
\end{lem}

This lemma is due to the standard theory in \cite{KnVo95Cohomological}. By this lemma, we prove the following theorem.

\begin{thm}\label{Lowest and Lowest Degree K-types}
Let $\pi=\pi(r,\lambda,\Psi,\mu,\nu)$ be a representation in $\mathcal{R}(Sp(p,q),\omega_{p,q,n})$. Then $\mathcal{A}(\pi)\subseteq \mathcal{A}_{p,q,n}(\pi)$.
\end{thm}
\begin{pf}
By Theorem \ref{Algorithm1} and Lemma \ref{Correspondence in Joint Harmonics}, it is not difficult to check that the degrees of all lowest $Sp(p)\times Sp(q)$-types of $\pi$ equal to each other. By Lemma \ref{Form of Lowest Degree K-types}, any $Sp(p)\times Sp(q)$-type $\sigma$ occurring in $\pi(r,\lambda,\Psi,\mu,\nu)$ is of the form $\Lambda+\Sigma n_{\alpha}\alpha$, where $\Lambda$ is a lowest $Sp(p)\times Sp(q)$-type of $\pi$. The roots in $\Delta(\mathfrak{l}:\mathfrak{t})\cup \Delta(\mathfrak{u}\cap \mathfrak{p})$ are of the form
\begin{equation}\label{Form of Roots}
\begin{split}
  \pm(e_{i}-e_{j})   \quad     1\leq i <j\leq p,\\
  \pm(f_{i}-f_{j})   \quad    1\leq i <j\leq q,\\
  \pm(e_{i}-f_{j})   \quad 1\leq i\leq p,1\leq j\leq q,\\
   e_{i}+f_{j}   \quad  1\leq i\leq p,1\leq j\leq q.
\end{split}
\end{equation}
Let $\Sigma n_{\alpha}\alpha=(x_{1},\cdots,x_{p};y_{1},\cdots,y_{q})$. By (\ref{Degree Formula1}),  $\mathrm{deg}(\sigma)=\mathrm{deg}(\Lambda)+\sum_{i=1}^{p}x_{i}+\sum_{j=1}^{q}y_{j}$. For each $\alpha$ of the form in (\ref{Form of Roots}), the summation of its all coordinates is non-negative. Then $\mathrm{deg}(\sigma)\geq \mathrm{deg}(\Lambda)$ and $\Lambda$ is a lowest degree $Sp(p)\times Sp(q)$-type of $\pi$. Furthermore, all lowest $Sp(p)\times Sp(q)$-types are of minimal degree. We finish the proof.
\end{pf}

For $\pi^{'}\in \mathcal{R}(O^{*}(2n),\omega_{p,q,n})$,  $\mathcal{A}(\pi^{'})$ is not necessarily contained in $\mathcal{A}_{p,q,n}(\pi^{'})$ since the degree formula (\ref{Degree Formula2}) depends on the choice of $(p,q)$.

\section{Parabolic Induction Principle}
\newtheorem{prop}[thm]{Proposition}

In this chapter, we briefly review the mechanism of parabolic induction principle. The induction principle is due to Kudla \cite{Ku1986Local}. The proof of parabolic induction principle is based on the explicit realization of the Schrodinger Model (see \cite{Rao93Explicit}). We mainly follow the content in $\S$4 of \cite{LPTZ03Correspondence}. The detailed proof of induction principle can also be founded in $\S$3 of \cite{AB95Complex} and $\S$4 of \cite{Paul1998Unitary}. Using induction principle, we prove the Going-Up theorem which is useful in explicit calculation.

Let $\mathbf{D}=\mathbb{H}$ with the standard involution as in $\S$2.

Let $V$ (resp. $V^{'}$) be a right (resp. left) $\mathbf{D}$ space equipped with a Hermitian form $(,)$ (resp. skew-Hermitian form $<,>$). We define $W=V\otimes _{\mathbf{D}}V^{'}$ to be a real symplectic space with the symplectic form $\langle, \rangle=\mathrm{tr}_{\mathbf{D}/\mathbb{R}}((,)\otimes \overline{<,>})$, where $\mathrm{tr}$ is the reduced trace on $\mathbb{H}$ such that $\mathrm{tr}(a+b\textbf{i}+c\textbf{j}+d\textbf{k})=2a$. We denote by $G=G(V)$ (resp. $G^{'}=G(V^{'})$) the isometry group of $(V,(,))$ (resp. $(V^{'},<,>)$). As in $\S$3, $(G(V),G(V^{'}))$ can be embedded into $Sp(W)$ and is a reductive dual pair in this symplectic group.

Suppose that
\begin{equation}\label{Polar Decomposition3}
  V=V_{+}\oplus V_{0}\oplus V_{-},
\end{equation}
\begin{equation}\label{Polar Decomposition4}
  V^{'}=V_{+}^{'}\oplus V_{0}^{'}\oplus V_{-}^{'},
\end{equation}
where $V_{+},V_{-}$ (resp. $V_{+}^{'},V_{-}^{'}$) are totally isotropic subspaces and are dual to each other with respect to $(,)$ (resp. $<,>$). Denote $d=\dim_{\mathbb{R}}\mathbb{H}=4$ and $d_{0}=\dim_{\mathbb{R}}\{t\in \mathbb{H}|\overline{t}=-t\}=3$. Set $m_{0}=\dim_{\mathbb{H}}V_{0}$, $n_{0}=\dim_{\mathbb{H}}V_{0}^{'}$, $k=\dim_{\mathbb{H}}V_{+}$ and $l=\dim_{\mathbb{H}}V_{+}^{'}$.

Let $W_{0}=V_{0}\otimes_{\mathbb{H}} V^{'}_{0}$ and let  $W_{0}=X_{0}\oplus Y_{0}$ is a complete polarization. Then $W=X\oplus Y$ is a complete polarization of $W$, where
\begin{equation}\label{Polar Decomposition5}
   X=(V\otimes_{\mathbb{H}} V_{+}^{'})\oplus (V_{+}\otimes_{\mathbb{H}} V_{0}^{'})\oplus (X_{0}),
\end{equation}
\begin{equation}\label{Polar Decomposition6}
  Y=(V\otimes_{\mathbb{H}} V_{-}^{'})\oplus (V_{-}\otimes_{\mathbb{H}} V_{0}^{'})\oplus (Y_{0}).
\end{equation}

Let $P_{V}=P(V_{+})$ be the stabilizer group of $V_{+}$ in $G(V)$ and let $P_{V^{'}}=P(V_{+}^{'})$ be the stabilizer group of  $V_{+}^{'}$ in $G(V^{'})$. Then $P_{V}=L_{V}N_{V}$ and $P_{V^{'}}=L_{V^{'}}N_{V^{'}}$, where $L_{V}\cong G(V_{0})\times GL(k,\mathbb{H})$  as in (\ref{Cuspidal Parabolic Subgroup2}) and $L_{V^{'}}\cong G(V_{0}^{'})\times GL(l,\mathbb{H})$ as in (\ref{Cuspidal Parabolic Subgroup4}).

Let $W_{L}=(V_{-}\otimes_{\mathbb{H}} V_{+}^{'})\oplus (V_{+}\otimes_{\mathbb{H}} V_{-}^{'})\oplus (V_{0}\otimes_{\mathbb{H}} V_{0}^{'})$. Then $W_{L}$ is a real symplectic subspace of $W$. Let $W_{L}=X_{L}\oplus Y_{L}$ be a complete polarization of $W_{L}$, where
\begin{equation}\label{Polar Decomposition7}
  X_{L}=(V_{-}\otimes_{\mathbb{H}} V_{+}^{'})\oplus X_{0},
\end{equation}
\begin{equation}\label{Polar Decomposition8}
  Y_{L}=(V_{+}\otimes_{\mathbb{H}} V_{-}^{'})\oplus Y_{0}.
\end{equation}
Then $(L_{V},L_{V^{'}})$ is a reductive dual pair in $Sp(W_{L})$. In fact, $(G(V_{0}),G(V_{0}^{'}))$ is a reductive dual pair of type I and $(GL(k,\mathbb{H}),GL(l,\mathbb{H}))$ is a reductive dual pair of type II (c.f. $\S$6 of \cite{LPTZ03Correspondence}).

For a fixed additive character $\psi(t)=e^{2\pi\textbf{i}t}$, we denote by $\omega$ the oscillator representation attached to the dual pair $(G(V),G(V^{'}))\subseteq Sp(W)$ and by $\omega_{L}$ the oscillator representation attached to the dual pair $(L_{V},L_{V^{'}})\subseteq Sp(W_{L})$. The oscillator representation $\omega$ may be realized on the Schwartz space $S(Y)$ and the oscillator representation $\omega_{L}$ may be realized on the Schwartz space $S(Y_{L})$ (see \cite{Rao93Explicit}). Let
\begin{equation}\label{Restriction Map}
  \rho: S(Y)\rightarrow S(Y_{L})
\end{equation}
be the restriction map. Obviously $\rho$ is surjective.

By the homomorphism (\ref{Injective Homomorphism}), the group $GL(k,\mathbb{H})$ can be embedded into the group $GL(2k,\mathbb{C})$. For every $h\in GL(k,\mathbb{H})$, let $\mathrm{det}(h)$ be the usual determinant of $h$ realized as an element in $GL(2k,\mathbb{C})$. Similarly, we define the determinant of $GL(l,\mathbb{H})$. Let $\xi$ be a character of $P_{V}\times P_{V^{'}}$:
\begin{equation}\label{Character of Induction Principle}
  \xi((g_{0}hn, g_{0}^{'}h^{'}n^{'}))=\mathrm{det}(h)^{n_{0}+l}\mathrm{det}(h^{'})^{m_{0}+k},
\end{equation}
where $(g_{0},g_{0}^{'})\in (G(V_{0}),G(V_{0}^{'}))$, $(h,h^{'})\in (GL(k,\mathbb{H}),GL(l,\mathbb{H}))$ and $(n,n^{'})\in (N_{V},N_{V^{'}})$.  Since $N_{V}$ (resp. $N_{V^{'}}$) is normal in $P_{V}$ (resp. $P_{V^{'}}$), the character is well-defined. Then we have the following theorem (Proposition 4.14 in \cite{LPTZ03Correspondence}).

\begin{thm}
Let $\xi$ be the character as in $(\ref{Character of Induction Principle})$. The restriction map $\rho$ in $(\ref{Restriction Map})$ is a  $P_{V}\times P_{V^{'}}$-equivalent map
\begin{equation}
  \omega\rightarrow \omega_{L}\otimes \xi,
\end{equation}
where $N_{V}\times N_{V^{'}}$ acts trivially on $\omega_{L}$.
\end{thm}

The detail of proof can be founded in $\S$4 of \cite{Paul1998Unitary}. Let $\chi_{V}$ and $\chi_{V^{'}}$ be the characters of $GL(k,\mathbb{H})$ and $GL(l,\mathbb{H})$ given by
\begin{equation}\label{Form of Characters1}
  \chi_{V}(h)=\det(h)^{\frac{1}{4}((n_{0}+l-m_{0}-k+1)d-2d_{0})}, \quad \forall h\in GL(k,\mathbb{H}),
\end{equation}
\begin{equation}\label{Form of Characters2}
   \chi_{V^{'}}(h^{'})=\det(h^{'})^{\frac{1}{4}((m_{0}+k-n_{0}-l-1)d+2d_{0})}, \quad \forall h^{'}\in GL(l,\mathbb{H}).
\end{equation}
Using Frobenius Reciprocity, we state the following theorem on Induction Principle (Theorem 4.20 of \cite{LPTZ03Correspondence}, also see Theorem 5.10 of \cite{Paul2005Symplectic}).

\begin{thm}\label{Parabolic Induction Principle1}
Let $\pi\in \mathcal{R}(G(V_{0})$, $\pi^{'}\in \mathcal{R}(G(V_{0}^{'}))$, $\tau\in \mathcal{R}(GL(k,\mathbb{H}))$ and $\tau^{'}\in \mathcal{R}(GL(l,\mathbb{H}))$. Suppose that $\pi$ corresponds to $\pi^{'}$ for the real reductive dual pair $(G(V_{0}),G(V_{0}^{'}))$ and $\tau$ corresponds to $\tau^{'}$ for the real reductive dual pair $(GL(k,\mathbb{H}),GL(l,\mathbb{H}))$. There is a nonzero $G(V)\times G(V^{'})$ map (on the level of Harish-Chandra module)
\begin{equation}\label{Projection1}
  \Phi: \omega\rightarrow Ind_{P_{V}}^{G(V)}(\pi\otimes \tau\otimes \chi_{V}\otimes 1{\kern-0.36em} 1)\otimes Ind_{P_{V^{'}}}^{G(V^{'})}(\pi^{'}\otimes \tau^{'}\otimes \chi_{V^{'}}\otimes 1{\kern-0.36em} 1).
\end{equation}
\end{thm}

 To be more precise, we state the following theorem (Theorem 4.24 of \cite{LPTZ03Correspondence}). The details of proof can be founded in $\S$3 of \cite{AB95Complex}.

\begin{thm}\label{Parabolic Induction Principle2}
 Let $\Phi$ be as in Theorem $\ref{Parabolic Induction Principle1}$, Denote $G=G(V)$, $G^{'}=G(V^{'})$, $L=L_{V}$ and $L^{'}=L_{V^{'}}$, and let $K$ and $K^{'}$ be the maximal compact subgroups of $G$ and $G^{'}$, respectively. Suppose $\sigma$ is a $K$-type and $\kappa$ is a $K\cap L$-type such that the following conditions are satisfied:
\newline
$(1)$ $\kappa$ occurs and is of minimal degree in $\pi\otimes \tau$.
\newline
$(2)$ $\sigma$ occurs and is of minimal degree and of multiplicity one in $Ind_{P_{V}}^{G(V)}(\pi\otimes \tau\otimes \chi_{V}\otimes 1{\kern-0.36em} 1)$.
\newline
$(3)$ $\sigma$ and $\kappa$ have the same degree, and the restriction of $\sigma$ to $K\cap L$ contains $\kappa$.
\newline
$(4)$ There exist characters $\alpha$ and $\alpha^{'}$ of $L$ and $L^{'}$ which are trivial on $K\cap L$ and $K^{'}\cap L^{'}$, such that $\pi\otimes \tau\otimes \alpha$ and $\pi^{'}\otimes \tau^{'}\otimes \alpha^{'}$ correspond to each other for the reductive dual pair $(L,L^{'})$, and $Ind_{P_{V}}^{G(V)}(\pi\otimes \tau\otimes \alpha\otimes \chi_{V}\otimes 1{\kern-0.36em} 1)$ is irreducible.

Let $\sigma^{'}$ be the $K^{'}$ type which corresponds to $\sigma$ in the joint harmonics space. Then $\sigma\otimes \sigma^{'}$ is in the image of $\Phi$.

The statement of this theorem is also true if we exchange the role of $V$ and $V^{'}$.
\end{thm}

For further application, we need to know something about the theta correspondence for the reductive dual pairs $(GL(n,\mathbb{H}),GL(m,\mathbb{H}))$. For positive integers $n,m$, $(GL(n,\mathbb{H}),GL(m,\mathbb{H}))$ is a reductive dual pair in $Sp(8mn,\mathbb{R})$ (see $\S$6 of \cite{LPTZ03Correspondence}). The maximal compact subgroups of $GL(n,\mathbb{H})$ and $GL(m,\mathbb{H})$ are isomorphic to $Sp(n)$ and $Sp(m)$, respectively. We state the following two propositions (Proposition 4.25 and 4.27 of \cite{LPTZ03Correspondence}).

\begin{prop}\label{Type II1}
Let $n\leq m$. The correspondence of $Sp(n)$-types and $Sp(m)$-types in the joint harmonics for the reductive dual pair $(GL(n,\mathbb{H}),GL(m,\mathbb{H}))$ is given as follow. If $\sigma=(a_{1},\cdots,a_{n})$ is a $Sp(n)$-type, then $\sigma$ occurs in the joint harmonics and corresponds to the $Sp(m)$-type $\sigma^{'}=(a_{1},\cdots,a_{n},0,\cdots,0)$. The degree of $\sigma$ is $\Sigma_{i=1}^{n}a_{i}$.
\end{prop}

\begin{prop}\label{Type II2}
Let $\tau=\tau(\mu,\nu)$ be an irreducible admissible representation of $GL(1,\mathbb{H})$ as in $\S\mathrm{2.1}$. Then $\tau$ corresponds to $\tau^{\vee}$ for the reductive dual pair $(GL(1,\mathbb{H}),GL(1,\mathbb{H}))$, where $\tau^{\vee}=\tau(\mu,-\nu)$.
\end{prop}

Using Induction Principle, we prove the following theorem.

\begin{thm}\label{Going-Up}
Let $\pi^{'}$ be an irreducible admissible representation of $O^{*}(2n)$. Suppose $\theta_{p,q}(\pi^{'})$ is nonzero and $\theta_{p,q}(\pi^{'})=\pi(r,\lambda,\Psi,\mu,\nu)$. Let $s$ be a non-negative integer. Then
\begin{equation*}
  \theta_{p+s,q+s}(\pi^{'})=\pi(r+s,\lambda,\Psi,\mu^{s},\nu^{s}),
\end{equation*}
where $\mu^{s}=(\mu,\underbrace{1,\cdots,1}_{s})$ and
\begin{equation*}
  \nu^{s}=(\nu,2p+2q-2n+3,2p+2q-2n+7,\cdots,2p+2q-2n+4s-1).
\end{equation*}
\end{thm}
\begin{pf}
We use induction on $s$. Let $\pi_{s}=\pi(r+s,\lambda,\Psi,\mu^{s},\nu^{s})$. By Theorem \ref{Joint Harmonics} and Theorem \ref{Lowest and Lowest Degree K-types}, there exists a lowest $Sp(p)\times Sp(q)$-type $\eta_{0}$ of $\pi_{0}=\theta_{p,q}(\pi^{'})$ such that $\eta_{0}$ corresponds to a lowest degree $U(n)$-type $\eta^{'}$ of $\pi^{'}$ in $\mathscr{H}_{p,q,n}$. Assume $\theta_{p+s-1,q+s-1}(\pi^{'})=\pi_{s-1}=\pi(r+s-1,\lambda,\Psi,\mu^{s-1},\nu^{s-1})$. Take $G(V)=Sp(p+s,q+s)$, $G(V^{'})=O^{*}(2n)$, $L=Sp(p+s-1,q+s-1)\times GL(1,\mathbb{H})$ and $L^{'}=O^{*}(2n)$. By Theorem \ref{Parabolic Induction Principle1}, there is a nonzero $Sp(p+s,q+s)\times O^{*}(2n)$ map
\begin{equation*}
  \omega_{p+s,q+s,n}\mapsto Ind_{P}^{Sp(p+s,q+s)}(\pi_{s-1}\otimes \chi_{s}\otimes 1{\kern-0.36em} 1)\otimes \pi^{'},
\end{equation*}
where $P=LN$ with $L=Sp(p+s-1,q+s-1)\cdot GL(1,\mathbb{H})$, and $\chi_{s}$ is a character of $GL(1,\mathbb{H})$ such that
\begin{equation*}
  \chi_{s}(h)=\mathrm{det}(h)^{n-p-q-2s+\frac{1}{2}}
\end{equation*}
for each $h\in GL(1,\mathbb{H})$. Then $\theta_{p+s,q+s}(\pi^{'})$ is a constituent of the induced representation $Ind_{P}^{Sp(p+s,q+s)}(\pi_{s-1}\otimes \chi_{s}\otimes 1{\kern-0.36em} 1)$. Let $P_{s-1}$ be a cuspidal subgroup of $Sp(p+s-1,q+s-1)$ with Levi component $L_{s-1}\cong Sp(p-r,q-r)\times GL(1,\mathbb{H})^{r+s-1}$. Then $\pi_{s-1}$ is a constituent of the standard module  $X(r+s-1,\lambda,\Psi,\mu^{s-1},\nu^{s-1})$. Since $P_{s-1}\cdot GL(1,\mathbb{H})$ is a parabolic subgroup of $L$, then $P_{s}=P_{s-1}\cdot GL(1,\mathbb{H})\cdot N$ is a parabolic subgroup of $Sp(p+s,q+s)$. By the double induction formula in $\S\mathrm{7.2}$ of \cite{Knapp86Overview}, the induced representation $Ind_{P}^{Sp(p+s,q+s)}(\pi_{s-1}\otimes \chi_{s}\otimes 1{\kern-0.36em} 1)$ is a sub-quotient of the standard module $X(r+s,\lambda,\Psi,\mu_{s},\nu_{s})$, where $\mu_{s}=(\mu^{s-1},1)$ and $\nu_{s}=(\nu^{s-1},2n-2p-2q-4s+1)$. Furthermore, $\theta_{p+s,q+s}(\pi^{'})$ is a constituent of this standard module. Suppose $\eta_{0}=(a_{1},\cdots,a_{p};b_{1},\cdots,b_{q})$. Let $\eta_{s}$ be a $Sp(p+s)\times Sp(q+s)$-type such that
\begin{equation*}
  \eta_{s}=(a_{1},\cdots,a_{p},\overbrace{0,\cdots,0}^{s};b_{1},\cdots,b_{q},\overbrace{0,\cdots,0}^{s}).
\end{equation*}
Since $\eta_{0}$ corresponds to $\eta^{'}$ in the joint harmonics $\mathscr{H}_{p,q,n}$, then $\eta_{s}$ corresponds to $\eta^{'}$ in the joint harmonics $\mathscr{H}_{p+s,q+s,n}$ by Lemma \ref{Correspondence in Joint Harmonics}. By Theorem \ref{Algorithm1} (also see the proof of Theorem 5.8 in \cite{LPTZ03Correspondence}) and Remark \ref{Multiplicity}, $\eta_{s}$ is the lowest $Sp(p+s)\times Sp(q+s)$-type of the standard module $X(r+s,\lambda,\Psi,\mu_{s},\nu_{s})$ with multiplicity one. By Theorem \ref{Joint Harmonics}, $\pi^{'}$ corresponds to $\pi(r+s,\lambda,\Psi,\mu_{s},\nu_{s})$, the unique constituent of $X(r+s,\lambda,\Psi,\mu_{s},\nu_{s})$ containing $\eta_{s}$. By Theorem \ref{Langlands Classification}, we finish the proof.
\end{pf}

\section{Explicit Theta Correspondence for $(Sp(p,q),O^{*}(2))$}

In this chapter, we briefly state the explicit theta correspondence for reductive dual pairs $(Sp(p,q),O^{*}(2))$. This work is done by R. Howe and J.-S.  Li. Their strategy is suggestive for calculating explicit theta correspondence for reductive dual pairs $(Sp(p,q),O^{*}(4))$.

By Lemma \ref{Contragredient Representation} ,Theorem \ref{Stable Ranger Theorem} and Theorem \ref{Going-Up}, the problem is reduced to determination of $\theta_{p,q}$ for $p\geq q$ and $q\leq 1$. Irreducible admissible representations of $O^{*}(2)=U(1)$ are $\chi_{k}:t\mapsto t^{k}$ with $k\in \mathbb{Z}$.

The following is a list of systems of positive roots which will be mentioned in this chapter:
\begin{enumerate}[(1)]
  \item For each integer $p\geq 0$, we denote by $\Psi_{1}$ the system of positive roots for the compact group $Sp(p)$ such that $(p,p-1,\cdots,1)$ is dominant with respect to $\Psi_{1}$;
  \item For each integer $p\geq 1$, we denote by $\Psi_{2}$ the system of positive roots for the group $Sp(p,1)$ such that $(p,p-1,\cdots,1;p+1)$ is dominant with respect to $\Psi_{2}$;
  \item For each integer $p\geq 1$, we denote by $\Psi_{3}$ the system of positive roots for the group $Sp(p,1)$ such that $(p+1,p,\cdots,2;1)$ is dominant with respect to $\Psi_{3}$.
\end{enumerate}

For $p\geq 1$, we denote by $\mathbf{A}_{p,1}$ the set of elements in $\mathcal{R}(Sp(p,1))$ whose lowest $Sp(p)\times Sp(1)$-types are all of the form $(a_{1},0,\cdots,0;b_{1})$. In the appendix, we list all infinitesimal equivalence classes in $\mathbf{A}_{p,1}$ in terms of Langlands parameters.

\begin{thm}
The following is a list of all explicit theta $(p,0)$-lifts on $\mathcal{R}(O^{*}(2))$ for $p\geq 0$:
\begin{enumerate}[(1)]
  \item If $p=0$, then
\begin{equation*}
  \theta_{0,0}(\chi_{k})=\left\{
     \begin{array}{ll}
       \pi(0,\emptyset), & \hbox{if $k=0$,} \\
       0, & \hbox{if $k\neq 0$.}
     \end{array}
   \right.
\end{equation*}
  \item If $p\geq 1$, then
\begin{equation*}
  \theta_{p,0}(\chi_{k})=\left\{
                            \begin{array}{ll}
                              0, & \hbox{if $k<p$,} \\
                              \pi((k,p-1,p-2,\cdots,1),\Psi_{1}), & \hbox{if $k\geq p$.}
                            \end{array}
                          \right.
\end{equation*}
\end{enumerate}
\end{thm}
\begin{pf}
First assume $p=0$. Then $Sp(0,0)$ is the trivial group and the explicit theta correspondence $\theta_{0,0}$ is obvious. The trivial representation of $Sp(0,0)$ in terms of Langlands parameters is $\pi(0,\emptyset)$. Next assume $p>0$. By Lemma \ref{Correspondence in Joint Harmonics}, the $U(1)$-type $(k)$ does not occur in the Fock space $\mathscr{F}_{p,0,1}$ if $k<p$. Then $\theta_{p,0}(\chi_{k})=0$ if $k<p$. If $k\geq p$, by Proposition 6.1 and Theorem 6.2 of \cite{Li90ThetaLifting} (also see Proposition 6.1 of \cite{VZ84Unitary}), $\chi_{k}$ can be theta lifted to a unitary representation with nonzero cohomology. By Theorem \ref{Joint Harmonics}, Lemma \ref{Correspondence in Joint Harmonics}, Theorem \ref{Correspondence of Infinitesimal Characters} and Theorem \ref{Lowest and Lowest Degree K-types}, $\theta_{p,0}(\chi_{k})$ has a lowest $Sp(p)$-type $(k-p,0,\cdots,0)$ and its infinitesimal character is $(k,p-1,p-2,\cdots,1)$. It is easy to see that
\begin{equation*}
  \theta_{p,0}=\pi((k,p-1,p-2,\cdots,1),\Psi_{1})
\end{equation*}
if $k\geq p$.
\end{pf}

\begin{thm}\label{Explicit Theta $(p,1)$-Lifts of Characters}
The following is a list of all explicit theta $(p,1)$-lifts on $\mathcal{R}(O^{*}(2))$ for $p\geq 1$:
\begin{enumerate}[(1)]
  \item If $p=1$, then
\begin{equation*}
  \theta_{1,1}(\chi_{k})=\left\{
                            \begin{array}{ll}
                              \pi((k;1),\Psi_{3}),  & \hbox{if $k>0$,} \\
                              \pi((1;-k),\Psi_{2}), & \hbox{if $k<0$,} \\
                              \pi(1,0,\emptyset,1,1), & \hbox{if $k=0$.}
                            \end{array}
                          \right.
\end{equation*}
  \item If $p>1$ and $k\geq p-1$, then
\begin{equation*}
  \theta_{p,1}(\chi_{k})= \pi(1,(k,p-2,p-3,\cdots,1),\Psi_{1},1,2p-1).
\end{equation*}
  \item If $p>1$ and $-p<k<p-1$, then
\begin{equation*}
  \theta_{p,1}(\chi_{k})=\pi(1,(p-1,p-2,\cdots,1),\Psi_{1},p-k,p+k).
\end{equation*}
  \item If $p>1$ and $k\leq -p$, then
\begin{equation*}
  \theta_{p,1}(\chi_{k})=\pi((p,p-1,\cdots,1;-k),\Psi_{2}).
\end{equation*}
\end{enumerate}
\end{thm}
\begin{pf}
By Theorem \ref{Stable Ranger Theorem}, $\theta_{p,1}(\chi_{k})$ is nonzero for each $k\in \mathbb{Z}$. If $k\geq p-1$, by Theorem \ref{Joint Harmonics}, Lemma \ref{Correspondence in Joint Harmonics}, Theorem \ref{Correspondence of Infinitesimal Characters} and Theorem \ref{Lowest and Lowest Degree K-types}, we know that $\theta_{p,1}(\chi_{k})$ has a unique lowest $Sp(p)\times Sp(1)$-type $(k-p+1,0,\cdots,0;0)$ and its infinitesimal character is $(k,p,p-1,\cdots,1)$. By the appendix, there is a unique element in $\mathbf{A}_{p,1}$ satisfying these conditions. Similarly, if $k<p-1$, $\theta_{p,1}(\chi_{k})$ has a unique lowest $Sp(p)\times Sp(1)$-type $(0,0,\cdots,0;p-1-k)$ and its infinitesimal character is $(k,p,p-1,\cdots,1)$. By the appendix, there is a unique element in $\mathbf{A}_{p,1}$ satisfying these conditions. Then all desired theta $(p,1)$-lifts in this theorem are determined in terms of Langlands parameters.
\end{pf}

\section{Explicit Theta Correspondence For $(Sp(p,q),O^{*}(4))$}

In this chapter, we determine theta correspondence for reductive dual pairs $(Sp(p,q),O^{*}(4))$ explicitly in terms of Langlands parameters. By Lemma \ref{Contragredient Representation}, Theorem \ref{Stable Ranger Theorem} and Theorem \ref{Going-Up}, the problem is reduced to determination of $\theta_{p,q}$ explicitly for $p\geq q$ and $q\leq 2$. Here we list all infinitesimal equivalence classes in $\mathcal{R}(O^{*}(4))$ as in Theorem \ref{Admissible Dual2}:
\begin{enumerate}[(1)]
  \item $P_{\lambda_{1},\lambda_{2}}$ with $\lambda_{1},\lambda_{2}\notin \mathbb{Z}$, $\mathrm{Re}(\lambda_{1}+\lambda_{2})\geq 0$ and $\lambda_{1}-\lambda_{2}\in \mathbb{Z}_{\geq 1}$;
  \item $D_{\lambda_{1},\lambda_{2}}$ with $\lambda_{1},\lambda_{2}\in \mathbb{Z}$, $\lambda_{1}+\lambda_{2}\in \mathbb{Z}_{\geq 0}$ and  $\lambda_{1}-\lambda_{2}\in \mathbb{Z}_{\geq 1}$;
  \item $\overline{D}_{\lambda_{1},\lambda_{2}}$ with $\lambda_{1},\lambda_{2}\in \mathbb{Z}$, $\lambda_{1}+\lambda_{2}\in \mathbb{Z}_{\leq 0}$ and  $\lambda_{1}-\lambda_{2}\in \mathbb{Z}_{\geq 1}$;
  \item $F_{\lambda_{1},\lambda_{2}}$ with $\lambda_{1},\lambda_{2}\in \mathbb{Z}$, $\lambda_{1}+\lambda_{2}\in \mathbb{Z}_{\geq 1}$ and  $\lambda_{1}-\lambda_{2}\in \mathbb{Z}_{\geq 1}$.
\end{enumerate}

Let condition (A) and (B) be those of $\S$2.1. Let $\Psi_{1}$, $\Psi_{2}$ and $\Psi_{3}$ be those in $\S$5. The following is a complementary list of the systems of positive roots which will be mentioned in this chapter:
\begin{enumerate}[(1)]
  \item For $p\geq 1$, we denote by $\Psi_{4}$ the system of positive roots of the group $Sp(p,1)$ such that $(p+1,p-1,p-2,\cdots,1;p)$ is dominant with respect to $\Psi_{4}$;
  \item For $p\geq 2$, we denote by $\Psi_{5}$ the system of positive roots of the group $Sp(p,2)$ such that $(p,p-1,\cdots,1;p+2,p+1)$ is dominant with respect to $\Psi_{5}$;
  \item For $p\geq 2$, we denote by $\Psi_{6}$ the system of positive roots of the group $Sp(p,2)$ such that $(p+2,p+1,\cdots,3;2,1)$ is dominant with respect to $\Psi_{6}$.
\end{enumerate}

\subsection{Explicit Theta $(p,0)$-Lifting for $O^{*}(4)$}

\begin{thm}\label{Explicit Theta (p,0)-Lifts}
The following is a list of all nonzero theta $(p,0)$-lifts on $\mathcal{R}(O^{*}(4))$ for $p\geq 0$:
\begin{enumerate}[(1)]
  \item If $p=0$, then
\begin{equation*}
  \theta_{0,0}(F_{1,0})=\pi(0,\emptyset).
\end{equation*}
  \item If $p=1$, then
\begin{equation*}
  \theta_{1,0}(D_{\lambda_{1},0})=\pi((\lambda_{1}),\Psi_{1}).
\end{equation*}
  \item If $p\geq 2$ and $\lambda_{2}\geq p-1$, then
\begin{equation*}
  \theta_{p,0}(D_{\lambda_{1},\lambda_{2}})=\pi((\lambda_{1},\lambda_{2},p-2,p-3,\cdots,1),\Psi_{1}).
\end{equation*}
\end{enumerate}
\end{thm}
\begin{pf}
The theta correspondence for $(Sp(0,0),O^{*}(4))$ is trivial. The trivial representation of $Sp(0,0)$ is $\pi(0,\emptyset)$ and the trivial representation of $O^{*}(4)$ is $F_{1,0}$. By Theorem 5.1 of \cite{LPTZ03Correspondence}, we determine $\theta_{1,0}$ directly. Assume $p\geq 2$. For $\pi^{'}\in \mathcal{R}(O^{*}(4))$, $\theta_{p,0}(\pi^{'})$ is nonzero only if all $U(2)$-types in $\mathcal{K}(\pi^{'})$ occur in $\mathscr{F}_{p,0,2}$. By Theorem \ref{Admissible Dual2} and Lemma \ref{Correspondence in Joint Harmonics}, $\theta_{p,0}(\pi^{'})$ is nonzero only if $\pi^{'}$ is of the form $D_{\lambda_{1},\lambda_{2}}$ with $\lambda_{2}\geq p-1$. If $\lambda_{2}\geq p-1$, by Proposition 6.1 and Theorem 6.2 of \cite{Li90ThetaLifting},
\begin{equation*}
  \theta_{p,0}(D_{\lambda_{1},\lambda_{2}})=\pi((\lambda_{1},\lambda_{2},p-2,p-3,\cdots,1),\Psi_{1}).
\end{equation*}
\end{pf}

\subsection{Explicit Theta $(p,1)$-Lifting for $O^{*}(4)$}

Assume $p\geq 1$. Let $G=Sp(p,1)$ and let $G^{'}=O^{*}(4)$. As in $\S$2.4, $G^{'}=(G_{1}^{'}\times G_{2}^{'})/\{\pm I\}$ with $G_{1}^{'}=SU(2)$ and $G_{2}^{'}=SL(2,\mathbb{R})$. Let $P_{2}^{'}$ be the upper triangle matrix subgroups of $G_{2}^{'}$ and let $P^{'}=(G_{1}^{'}\times P_{2}^{'})/\{\pm I\}$. Let $P$ be a proper cuspidal parabolic subgroup of $G$. The Levi component $L$ of $P$ is isomorphic to $Sp(p-1)\times GL(1,\mathbb{H})$. For $\mu\in \mathbb{Z}_{\geq 1}$ and $\nu\in \mathbb{C}$, by Theorem \ref{Parabolic Induction Principle1} and Proposition \ref{Type II2}, we have a nonzero $Sp(p,1)\times O^{*}(4)$-projection (on the level of Harish-Chandra module)
\begin{equation}\label{Projection2}
  \Phi: \omega_{p,1,2}\mapsto Ind_{P}^{Sp(p,1)}(\tau(\mu,-\nu)\otimes 1{\kern-0.36em} 1)\otimes Ind_{P^{'}}^{O^{*}(4)}(\tau(\mu,\nu)\otimes 1{\kern-0.36em} 1),
\end{equation}
where $\tau(\mu,\nu)$ is the irreducible representation of $GL(1,\mathbb{H})$ as in $\S$2.1. For convenience, we denote the first standard module by $I(\mu,-\nu)$ and the second standard module by $I^{'}(\mu,\nu)$. We calculate the explicit theta $(p,1)$-lifts of irreducible principal series first.

\begin{thm}
The following is a list of explicit theta $(p,1)$-lifts of irreducible principal series:
\begin{equation*}
  \theta_{p,1}(P_{\lambda_{1},\lambda_{2}})=\pi(1,(p-1,p-2,\cdots,1),\Psi_{1},\lambda_{1}-\lambda_{2},\lambda_{1}+\lambda_{2}).
\end{equation*}
\end{thm}
\begin{pf}
By Theorem \ref{Admissible Dual2}, the representation $P_{\lambda_{1},\lambda_{2}}$ is the standard module $I^{'}(\lambda_{1}-\lambda_{2},\lambda_{1}+\lambda_{2})$. By (\ref{Projection2}), there is a nonzero $Sp(p,1)\times O^{*}(4)$-projection $\Phi$ from $\omega_{p,1,2}$ to $I(\lambda_{1}-\lambda_{2},-(\lambda_{1}+\lambda_{2}))\otimes I^{'}(\lambda_{1}-\lambda_{2},\lambda_{1}+\lambda_{2})$.
Then we know that $\theta_{p,1}(P_{\lambda_{1},\lambda_{2}})$ is a constituent of $I(\lambda_{1}-\lambda_{2},-(\lambda_{1}+\lambda_{2}))$. By Theorem 1.1 of \cite{SV80Reducibility}, $I(\lambda_{1}-\lambda_{2},-(\lambda_{1}+\lambda_{2}))$ is irreducible (also see Theorem 5.2.2 of \cite{LTZ01TensorProduct}) and $\theta_{p,1}(P_{\lambda_{1},\lambda_{2}})$ is just the standard module $I(\lambda_{1}-\lambda_{2},-(\lambda_{1}+\lambda_{2}))$. By Theorem \ref{Langlands Classification}, $I(\lambda_{1}-\lambda_{2},-(\lambda_{1}+\lambda_{2}))$ is just $\pi(1,(p-1,p-2,\cdots,1),\Psi_{1},\lambda_{1}-\lambda_{2},\lambda_{1}+\lambda_{2})$ .
\end{pf}

Next we calculate the explicit theta $(p,1)$-lifts of irreducible finite dimensional representations $F_{\lambda_{1},\lambda_{2}}$.

\begin{thm}\label{Explicit Theta Lifts of Finite Dimensional Representations}
The following is a list of explicit theta $(p,1)$-lifts of irreducible finite dimensional representations of $O^{*}(4)$ for $p\geq 1$:
\begin{enumerate}[(1)]
  \item If $\lambda_{1}\geq p$, then
\begin{equation*}
  \theta_{p,1}(F_{\lambda_{1},\lambda_{2}})=\pi(1,(p-1,p-2,\cdots,1),\Psi_{1},\lambda_{1}-\lambda_{2},\lambda_{1}+\lambda_{2}).
\end{equation*}
  \item If $\lambda_{1}<p$, then
\begin{equation*}
  \theta_{p,1}(F_{\lambda_{1},\lambda_{2}})=0.
\end{equation*}
\end{enumerate}
\end{thm}
\begin{pf}
First assume $\lambda_{1}<p$. By Theorem \ref{Admissible Dual2} and Lemma \ref{Correspondence in Joint Harmonics}, every $U(2)$-type in $\mathcal{K}(F_{\lambda_{1},\lambda_{2}})$ does not occur in the joint harmonics $\mathscr{H}_{p,1,2}$. Then $\theta_{p,1}(F_{\lambda_{1},\lambda_{2}})$ is zero by Theorem \ref{Joint Harmonics}.

Next assume $\lambda_{1}\geq p$. Let $\Phi$ be the nonzero projection from $\omega_{p,1,2}$  to $I(\lambda_{1}-\lambda_{2},-(\lambda_{1}+\lambda_{2}))\otimes I^{'}(\lambda_{1}-\lambda_{2},\lambda_{1}+\lambda_{2})$ as in (\ref{Projection2}).
Denote the first standard module by $I$ and the second standard module by $I^{'}$. If $\lambda_{2}\leq 1-p$, by Theorem \ref{Admissible Dual2} and Lemma \ref{Correspondence in Joint Harmonics}, we know that the $U(2)$-type $\eta^{'}=(-\lambda_{2},1-\lambda_{1})$ occurs in $F_{\lambda_{1},\lambda_{2}}$ (also in $I^{'}$) with multiplicity one and is of the minimal degree $\lambda_{1}-\lambda_{2}-1$. The restriction of $\eta^{'}$ to $SU(2)$ is the $SU(2)$-type $(\lambda_{1}-\lambda_{2}-1)$ and has degree $\lambda_{1}-\lambda_{2}-1$ by Proposition \ref{Type II1}. Then $\eta^{'}$ satisfy the condition (1), (2) and (3) of Theorem \ref{Parabolic Induction Principle2}, and the condition (4) is also satisfied by Theorem 1.1 of \cite{SV80Reducibility}. By Theorem \ref{Parabolic Induction Principle2}, we know that $\eta\otimes \eta^{'}$ occurs in $\mathrm{Im}\Phi$, where $\eta$ corresponds to $\eta^{'}$ in $\mathscr{H}_{p,1,2}$. By Theorem \ref{Admissible Dual1}, $F_{\lambda_{1},\lambda_{2}}$ can be regarded as a quotient of $I^{'}$ such that $F_{\lambda_{1},\lambda_{2}}=I^{'}/M^{'}$. Denote $\mathscr{M}=I\otimes M^{'}$. Then $\mathrm{Im}\Phi/(\mathrm{Im}\Phi\cap \mathscr{M})$ is a nonzero submodule of $I\otimes F_{\lambda_{1},\lambda_{2}}$ ($\eta\otimes \eta^{'}$ does not lie in $\mathscr{M}$.). Since $\mathrm{Im}\Phi/(\mathrm{Im}\Phi\cap \mathscr{M})$ is a $Sp(p,1)\times O^{*}(4)$-quotient of $\omega_{p,1,2}$, we know that $\theta_{p,1}(F_{\lambda_{1},\lambda_{2}})$ is a constituent of $I$. If $\lambda_{2}>1-p$, take $\eta^{'}=(p-1,p-\lambda_{1}+\lambda_{2})$. Similarly, we prove that $\theta_{p,1}(F_{\lambda_{1},\lambda_{2}})$ is a constituent of $I$.

By Theorem \ref{Algorithm1}, Theorem \ref{Admissible Dual2} and Lemma \ref{Correspondence in Joint Harmonics}, we can check the lowest $Sp(p)\times Sp(1)$-types of $I$ correspond to some lowest degree $U(2)$-types of $F_{\lambda_{1},\lambda_{2}}$ in $\mathscr{H}_{p,1,2}$ as follows.

If $\lambda_{1}-\lambda_{2}\leq 2p-1$, by Theorem \ref{Algorithm1} and Lemma \ref{Correspondence in Joint Harmonics}, the unique lowest $Sp(p)\times Sp(1)$-type of $I$ is $(0,\cdots,0;\lambda_{1}-\lambda_{2}-1)$, which corresponds to the $U(2)$-type $(p-1,p-\lambda_{1}+\lambda_{2})$ in $\mathscr{H}_{p,1,2}$. Since $\lambda_{1}\geq p$, $\lambda_{1}+\lambda_{2}\geq 1$ and $\lambda_{1}-\lambda_{2}\leq 2p-1$, we check that $-\lambda_{2}\leq p-1\leq \lambda_{1}-1$. By Theorem \ref{Admissible Dual2} and Lemma \ref{Correspondence in Joint Harmonics}, $(p-1,p-\lambda_{1}+\lambda_{2})$ is a lowest degree $U(2)$-type of $F_{\lambda_{1},\lambda_{2}}$.

If $\lambda_{1}-\lambda_{2}>2p-1$, by Theorem \ref{Algorithm1}, the lowest $Sp(p)\times Sp(1)$-types of $I$ are
\begin{equation*}
  (\frac{\lambda_{1}-\lambda_{2}+1+x}{2}-p,0,\cdots,0;\frac{\lambda_{1}-\lambda_{2}-3-x}{2}+p), \quad x=0,\pm 1.
\end{equation*}
Here $x$ is chosen so that $\lambda_{1}-\lambda_{2}+1+x$ is even. By Lemma \ref{Correspondence in Joint Harmonics}, the lowest $Sp(p)\times Sp(1)$-types correspond to
\begin{equation}\label{Form of K-types2}
  (\frac{\lambda_{1}-\lambda_{2}-1+x}{2},\frac{1+\lambda_{2}-\lambda_{1}+x}{2})
\end{equation}
in $\mathscr{H}_{p,1,2}$. Since $\lambda_{1}-\lambda_{2}>2p-1$, we know that $\frac{1+\lambda_{2}-\lambda_{1}+x}{2}\leq p-1\leq \frac{\lambda_{1}-\lambda_{2}-1+x}{2}$. By Lemma \ref{Correspondence in Joint Harmonics}, the $U(2)$-types in (\ref{Form of K-types2}) are lowest degree $U(2)$-types in $I^{'}$ of degree $\lambda_{1}-\lambda_{2}-1$. Since $\lambda_{1}+\lambda_{2}\geq 1$, we check that $-\lambda_{2}\leq \frac{\lambda_{1}-\lambda_{2}-1+x}{2}\leq \lambda_{1}-1$. By Theorem \ref{Admissible Dual2}, the $U(2)$-types in (\ref{Form of K-types2}) are lowest degree $U(2)$-types of $F_{\lambda_{1},\lambda_{2}}$.

 By Theorem \ref{Joint Harmonics}, we know that $\theta_{p,1}(F_{\lambda_{1},\lambda_{2}})$ is the constituent of $I$ containing all lowest $Sp(p)\times Sp(1)$-types. We finish the proof.
\end{pf}

Next we calculate the explicit theta $(p,1)$-lifts of limits of highest weight discrete series $\overline{D}_{\lambda_{1},\lambda_{2}}$.

\begin{thm}\label{Explicit Theta $(p,1)$-lifts of Highest Weight Discrete Series}
The following is a list of explicit theta $(p,1)$-lifts of limits of highest weight discrete series for $p\geq 1$:
\begin{enumerate}[(1)]
  \item If $\lambda_{1}\geq p$, then
\begin{equation*}
  \theta_{p,1}(\overline{D}_{\lambda_{1},\lambda_{2}})=\pi((\lambda_{1},p-1,p-2,\cdots,1;-\lambda_{2}),\Psi_{2}).
\end{equation*}
  \item If $\lambda_{1}<p$, then
\begin{equation*}
  \theta_{p,1}(\overline{D}_{\lambda_{1},\lambda_{2}})=0.
\end{equation*}
\end{enumerate}
\end{thm}
\begin{pf}
First assume $\lambda_{1}<p$. By Theorem \ref{Admissible Dual2} and Lemma \ref{Correspondence in Joint Harmonics}, every $U(2)$-type in $\mathcal{K}(\overline{D}_{\lambda_{1},\lambda_{2}})$ does not occur in the joint harmonics $\mathscr{H}_{p,1,2}$. Then $\theta_{p,1}(\overline{D}_{\lambda_{1},\lambda_{2}})$ is zero by Theorem \ref{Joint Harmonics}.

Next assume $\lambda_{1}\geq p$. Let $\Phi$ be the nonzero projection from $\omega_{p,1,2}$ to $I(\lambda_{1}-\lambda_{2},-(\lambda_{1}+\lambda_{2}))\otimes I^{'}(\lambda_{1}-\lambda_{2},\lambda_{1}+\lambda_{2})$ as in (\ref{Projection2}).
Denote the first standard module by $I$ and the second standard module by $I^{'}$. By Theorem \ref{Admissible Dual1},  $\overline{D}_{\lambda_{1},\lambda_{2}}$ can be realized as a quotient of $I^{'}$. By Theorem \ref{Admissible Dual2} and Lemma \ref{Correspondence in Joint Harmonics}, we know that $\eta^{'}=(\lambda_{1}-1,\lambda_{2})$ is a lowest degree $U(2)$-type of $\overline{D}_{\lambda_{1},\lambda_{2}}$ of degree $\lambda_{1}-\lambda_{2}-1$. As in the proof of Theorem \ref{Explicit Theta Lifts of Finite Dimensional Representations}, we prove that $\theta_{p,1}(\overline{D}_{\lambda_{1},\lambda_{2}})$ is a constituent of $I$. It will be proven in Proposition \ref{Conclusion2} that $\eta=(\lambda_{1}-p,0,\cdots,0;p-1-\lambda_{2})$ occurs in $\theta_{p,1}(\overline{D}_{\lambda_{1},\lambda_{2}})$ as a lowest $Sp(p)\times Sp(1)$-type. By Theorem \ref{Correspondence of Infinitesimal Characters}, the infinitesimal character of $\overline{D}_{\lambda_{1},\lambda_{2}}$ is $(\lambda_{1},\lambda_{2},p-1,p-2,\cdots,1)$. By the appendix, the unique element in $\mathcal{R}(Sp(p,1))$ satisfying these conditions is $\pi((\lambda_{1},p-1,p-2,\cdots,1;-\lambda_{2}),\Psi_{2})$. We finish the proof.
\end{pf}

Finally, we calculate the explicit theta $(p,1)$-lifts of limits of lowest weight discrete series $D_{\lambda_{1},\lambda_{2}}$.

\begin{thm}\label{Explicit Theta $(p,1)$-Lifts of Lowest Weight Discrete Series}
The following is a list of explicit theta $(p,1)$-lifts of limits of lowest weight discrete series for $p\geq 1$:
\begin{enumerate}[(1)]
  \item If $\lambda_{1}\leq p-1$ and $\lambda_{1}+\lambda_{2}\neq 0$, then
\begin{equation*}
  \theta_{p,1}(D_{\lambda_{1},\lambda_{2}})=\pi(1,(p-1,p-2,\cdots,1),\Psi_{1},\lambda_{1}-\lambda_{2},\lambda_{1}+\lambda_{2}).
\end{equation*}
  \item If $\lambda_{1}\leq p-1$ and $\lambda_{1}+\lambda_{2}=0$, then
\begin{equation*}
  \theta_{p,1}(D_{\lambda_{1},-\lambda_{1}})=\pi((p-1,\cdots,\lambda_{1}+1,\lambda_{1},\lambda_{1},\lambda_{1}-1,\cdots,1;\lambda_{1}),\Psi) \end{equation*}
with $\Psi$ uniquely determined by the condition $(A)$ of $\S$2.1.
  \item If $\lambda_{1}\geq p$ and $\lambda_{2}\leq -p$, then
\begin{equation*}
  \theta_{p,1}(D_{\lambda_{1},\lambda_{2}})=\pi((\lambda_{1},p-1,p-2,\cdots,1;-\lambda_{2}),\Psi_{4}).
\end{equation*}
  \item If $p\geq 2$, $\lambda_{1}\geq p$ and $1-p<\lambda_{2}\leq p-2$, then
\begin{equation*}
  \theta_{p,1}(D_{\lambda_{1},\lambda_{2}})=\pi(1,(\lambda_{1},p-2,p-3,\cdots,1),\Psi_{1},p-1-\lambda_{2},p-1+\lambda_{2}).
\end{equation*}
  \item If $p\geq 2$, $\lambda_{1}\geq p$ and $\lambda_{2}=1-p$, then
\begin{equation*}
  \theta_{p,1}(D_{\lambda_{1},1-p})=\pi((\lambda_{1},p-1,p-2,\cdots,1;p-1),\Psi_{4}).
\end{equation*}
  \item If $p\geq 3$ and $\lambda_{2}\geq p-1$, then
\begin{equation*}
  \theta_{p,1}(D_{\lambda_{1},\lambda_{2}})=\pi(1,(\lambda_{1},\lambda_{2},p-3,p-4,\cdots,1),\Psi_{1},1,2p-3).
\end{equation*}
  \item If $p=2$ and $\lambda_{2}\geq 1$, then
\begin{equation*}
  \theta_{2,1}(D_{\lambda_{1},\lambda_{2}})=\pi((\lambda_{1},\lambda_{2};1),\Psi_{3}).
\end{equation*}
  \item If $p=1$ and $\lambda_{2}\geq 0$, then
\begin{equation*}
  \theta_{1,1}(D_{\lambda_{1},\lambda_{2}})=0.
\end{equation*}
\end{enumerate}
\end{thm}
\begin{pf}
Let $\Phi$ be the nonzero projection from $\omega_{p,1,2}$ to $I(\lambda_{1}-\lambda_{2},\lambda_{1}+\lambda_{2})\otimes I^{'}(\lambda_{1}-\lambda_{2},-(\lambda_{1}+\lambda_{2}))$ as in (\ref{Projection2}).
Denote the first standard module by $I$ and the second standard module by $I^{'}$. By Theorem \ref{Admissible Dual1},  $D_{\lambda_{1},\lambda_{2}}$ can be realized as a quotient of $I^{'}$.

If $\lambda_{1}\leq p-1$, by Theorem \ref{Admissible Dual2} and Lemma \ref{Correspondence in Joint Harmonics}, we know that  $\eta^{'}=(p-1,p-\lambda_{1}+\lambda_{2})$ is a lowest degree $U(2)$-type of $D_{\lambda_{1},\lambda_{2}}$ of degree $\lambda_{1}-\lambda_{2}-1$. As in the proof of Theorem \ref{Explicit Theta Lifts of Finite Dimensional Representations}, we prove that $\theta_{p,1}(D_{\lambda_{1},\lambda_{2}})$ is a constituent of $I$.
If $\lambda_{1}+\lambda_{2}\neq 0$, by Theorem \ref{Algorithm1} and Lemma \ref{Correspondence in Joint Harmonics}, we know that  $\eta=(0,\cdots,0;\lambda_{1}-\lambda_{2}-1)$ is a lowest $Sp(p)\times Sp(1)$-type of $I$ and corresponds to $\eta^{'}$ in $\mathscr{H}_{p,1,2}$. By Theorem \ref{Joint Harmonics},  $\theta_{p,1}(D_{\lambda_{1},\lambda_{2}})$ is the unique constituent of $I$ containing all lowest $Sp(p)\times Sp(1)$-types. Then
\begin{equation*}
  \theta_{p,1}(D_{\lambda_{1},\lambda_{2}})=\pi(1,(p-1,p-2,\cdots,1),\Psi_{1},\lambda_{1}-\lambda_{2},\lambda_{1}+\lambda_{2}).
\end{equation*}
If $\lambda_{1}+\lambda_{2}=0$, the standard module $I$ is just the limit of discrete series (see Lemma 3.2.6 of \cite{Paul1998Unitary})
\begin{center}
  $\pi((p-1,\cdots,\lambda_{1}+1,\lambda_{1},\lambda_{1},\lambda_{1}-1,\cdots,1;\lambda_{1}),\Psi)$
\end{center}
with $\Psi$ uniquely determined by condition (A) since the Harish-Chandra parameter $(p-1,\cdots,\lambda_{1}+1,\lambda_{1},\lambda_{1},\lambda_{1}-1,\cdots,1;\lambda_{1})$ satisfies condition (B).

If $\lambda_{1}\geq p$ and $\lambda_{2}\leq p-2$, by Theorem \ref{Admissible Dual2} and Lemma \ref{Correspondence in Joint Harmonics}, we know that $\eta^{'}=(\lambda_{1},\lambda_{2}+1)$ is a lowest degree $U(2)$-type of $D_{\lambda_{1},\lambda_{2}}$ of degree $\lambda_{1}-\lambda_{2}-1$. As in the proof of Theorem 6.3, we prove that $\theta_{p,1}(D_{\lambda_{1},\lambda_{2}})$ is a constituent of $I$. It will be proven in Proposition \ref{Conclusion1} that $\theta_{p,1}(D_{\lambda_{1},\lambda_{2}})$ has a lowest $Sp(p)\times Sp(1)$-type $(\lambda_{1}-p+1,0,\cdots,0;p-2-\lambda_{2})$. By Theorem \ref{Correspondence of Infinitesimal Characters}, the infinitesimal character of $\theta_{p,1}(D_{\lambda_{1},\lambda_{2}})$ is $(\lambda_{1},\lambda_{2},p-1,p-2,\cdots,1)$. By the appendix, there is a unique element in $\mathcal{R}(Sp(p,1))$ satisfying these conditions. The explicit theta $(p,1)$-lifts are listed above.

If $p\geq 3$ and $\lambda_{2}\geq p-1$, by Theorem \ref{Explicit Theta (p,0)-Lifts} and Theorem \ref{Going-Up},
\begin{equation*}
  \theta_{p,1}(D_{\lambda_{1},\lambda_{2}})=\pi(1,(\lambda_{1},\lambda_{2},p-3,p-4,\cdots,1),\Psi_{1},1,2p-3).
\end{equation*}

If $p=1$ and $\lambda_{2}\geq 0$, by Theorem \ref{Admissible Dual2} and Lemma \ref{Correspondence in Joint Harmonics}, every $U(2)$-type in $\mathcal{K}(D_{\lambda_{1},\lambda_{2}})$ does not occur in the joint harmonics $\mathscr{H}_{1,1,2}$. Then
\begin{equation*}
  \theta_{1,1}(D_{\lambda_{1},\lambda_{2}})=0.
\end{equation*}

If $p=2$ and $\lambda_{2}\geq 2$, by Proposition 6.1 and Theorem 6.2 of \cite{Li90ThetaLifting},
\begin{center}
  $\theta_{2,1}(D_{\lambda_{1},\lambda_{2}})=\pi((\lambda_{1},\lambda_{2};1),\Psi_{3})$.
\end{center}
It will be proven in Proposition \ref{Conclusion3} that $\theta_{2,1}(D_{\lambda_{1},1})$ is nonzero. By Theorem \ref{Admissible Dual2}, Lemma \ref{Correspondence in Joint Harmonics}, Theorem \ref{Correspondence of Infinitesimal Characters} and Theorem \ref{Lowest and Lowest Degree K-types}, $\theta_{2,1}(D_{\lambda_{1},1})$ has a unique lowest $Sp(2)\times Sp(1)$-type $(\lambda_{1}-1,1;0)$ and its infinitesimal character is $(\lambda_{1},1,1)$. It is not difficult to check the unique element in $\mathcal{R}(Sp(p,1))$ satisfying these conditions is $\pi((\lambda_{1},1;1),\Psi_{3})$. We finish the proof.
\end{pf}

Here we give the proofs of three conclusions which remain unproved until now.

\begin{prop}\label{Conclusion1}
If $\lambda_{1}\geq p$, $\lambda_{2}\leq p-2$ and $\lambda_{1}+\lambda_{2}\geq 0$, then $\theta_{p,1}(D_{\lambda_{1},\lambda_{2}})$ has a unique lowest $Sp(p)\times Sp(1)$-type $(\lambda_{1}-p+1,0,\cdots,0;p-2-\lambda_{2})$.
\end{prop}
\begin{pf}
By Theorem \ref{Admissible Dual2} and Lemma \ref{Correspondence in Joint Harmonics}, we know that
\begin{equation*}
  \mathcal{A}_{p,1,2}(D_{\lambda_{1},\lambda_{2}})=\{(k+\lambda_{1}-\lambda_{2}-1,k)|\lambda_{2}+1\leq k\leq p-1\}.
\end{equation*}
By Theorem \ref{Joint Harmonics}, Lemma \ref{Correspondence in Joint Harmonics} and Theorem \ref{Lowest and Lowest Degree K-types},  $\mathcal{A}(\theta_{p,1}(D_{\lambda_{1},\lambda_{2}}))$ is contained in
\begin{equation*}
  \mathcal{D}=\{(k-p+\lambda_{1}-\lambda_{2},0,\cdots,0;p-1-k)|\lambda_{2}+1\leq k\leq p-1\}.
\end{equation*}
Then the lowest $Sp(p)\times Sp(1)$-types are those of minimal norm (see Definition 5.1 of \cite{Vogan79AlgebraicStructure}) in $\mathcal{D}$. Since $2\rho_{c}=(2p,2(p-1),\cdots,2;2)$, we just need to find out the smallest integer of
\begin{center}
  $\{(k+p+\lambda_{1}-\lambda_{2})^{2}+(p+1-k)^{2}|\lambda_{2}+1\leq k\leq p-1\}$.
\end{center}
Since $\lambda_{1}+\lambda_{2}\geq 0$, the smallest integer is obtained when $k=\lambda_{2}+1$. This means that the unique lowest $Sp(p)\times Sp(1)$-type of $\theta_{p,1}(D_{\lambda_{1},\lambda_{2}})$ is $(\lambda_{1}-p+1,0,\cdots,0;p-2-\lambda_{2})$.
\end{pf}

\begin{prop}\label{Conclusion2}
If $\lambda_{1}\geq p$, $\lambda_{2}\leq -p$ and $\lambda_{1}+\lambda_{2}\leq 0$, then $\theta_{p,1}(\overline{D}_{\lambda_{1},\lambda_{2}})$ has a unique lowest $Sp(p)\times Sp(1)$-type $(\lambda_{1}-p,0,\cdots,0;p-1-\lambda_{2})$.
\end{prop}
\begin{pf}
The proof is similar to the proof of Proposition \ref{Conclusion1}.
\end{pf}

To prove the occurrence of $D_{\lambda_{1},1}$ for $(Sp(2,1),O^{*}(4))$, we need to introduce Witt Towers and conservation relations. We define an equivalence relation on the set of right $\mathbb{H}$ Hermitian spaces. We say a right $\mathbb{H}$ Hermitian space $V_{1}$ with signature $(p_{1},q_{1})$ is equivalent to a right $\mathbb{H}$ Hermitian space $V_{2}$ with signature $(p_{2},q_{2})$ if and only if $p_{1}-q_{1}=p_{2}-q_{2}$. A Witt Tower $\mathrm{T}$ is defined to be an equivalence class in the set of right $\mathbb{H}$ Hermitian spaces. For non-negative integers $p$ and $q$, let $\mathbb{H}^{p,q}$ be the right column vector space $\mathbb{H}^{p+q}$ with Hermitian form (\ref{Hermitian form}). For $\pi^{'}\in \mathcal{R}(O^{*}(4))$, we define
\begin{equation*}
  \mathrm{n}_{\mathrm{T}}(\pi^{'})=\mathrm{min}\{p+q|\theta_{p,q}(\pi^{'})\neq 0, \mathbb{H}^{p,q}\in \mathrm{T}\}.
\end{equation*}
For two Witt towers $\mathrm{T}_{1}$ and $\mathrm{T}_{2}$, we define the distance between $\mathrm{T}_{1}$ and $\mathrm{T}_{2}$ to be
\begin{equation*}
  \mathrm{dist}(\mathrm{T}_{1},\mathrm{T}_{2})=|p_{1}+q_{2}-q_{1}-p_{2}|
\end{equation*}
if $\mathbb{H}^{p_{1},q_{1}}\in \mathrm{T}_{1}$ and $\mathbb{H}^{p_{2},q_{2}}\in \mathrm{T}_{2}$. It is easy to check that the distance is well-defined. We have the following Theorem (see Theorem 7.6 of \cite{SZ2012conservation}).

\begin{thm}\label{Conservation Relations}
For $\pi^{'}\in \mathcal{R}(O^{*}(4))$, there are two different Witt towers $\mathrm{T}_{1}$ and $\mathrm{T}_{2}$ such that
\begin{center}
$\mathrm{n}_{\mathrm{T}_{1}}(\pi^{'})+\mathrm{n}_{\mathrm{T}_{2}}(\pi^{'})=5$.
\end{center}
 For two different Witt towers $\mathrm{T}_{3}$ and $\mathrm{T}_{4}$,
\begin{center}
$\mathrm{n}_{\mathrm{T}_{3}}(\pi^{'})+\mathrm{n}_{\mathrm{T}_{4}}(\pi^{'})\geq 4+\mathrm{dist}(\mathrm{T}_{3},\mathrm{T}_{4})$.
\end{center}
\end{thm}

\begin{prop}\label{Conclusion3}
 Let $\pi^{'}=D_{\lambda_{1},1}$. Then $\theta_{2,1}(\pi^{'})$ is nonzero.
\end{prop}
\begin{pf}
We denote by $\mathrm{T}_{1}$ the Witt tower containing $\mathbb{H}^{2,0}$ and by $\mathrm{T}_{2}$ the Witt tower containing $\mathbb{H}^{2,1}$. By Theorem 5.1 of \cite{LPTZ03Correspondence} and Theorem \ref{Explicit Theta (p,0)-Lifts}, $\theta_{p,q}(\pi^{'})=0$ for all $p+q\leq 2$ except $(2,0)$. Then $\mathrm{n}_{\mathrm{T}}(\pi^{'})\geq 2$ for any Witt tower $\mathrm{T}$ and the equality holds if and only if $\mathrm{T}=\mathrm{T}_{1}$. By Theorem \ref{Conservation Relations}, there exists a Witt tower $\mathrm{T}^{'}$ such that $\mathrm{n}_{\mathrm{T}^{'}}(\pi^{'})=3$. By Theorem \ref{Explicit Theta (p,0)-Lifts} and Theorem \ref{Conservation Relations}, $\theta_{3,0}(\pi^{'})=0$ and $\mathrm{T}_{2}$ is the unique Witt Tower $\mathrm{T}^{'}$ such that $\mathrm{n}_{\mathrm{T}^{'}}(\pi^{'})=3$. This means that $\theta_{2,1}(\pi^{'})$ is nonzero.
\end{pf}

\subsection{Explicit Theta $(p,2)$-Lifting for $O^{*}(4)$}

Assume $p\geq 2$. By results in $\S$6.2, the following is a list of infinitesimal equivalence classes in $\mathcal{R}(O^{*}(4))$ whose theta $(p-1,1)$-lifts are zero:
\begin{enumerate}[(1)]
  \item Limits of highest weight discrete series $\overline{D}_{\lambda_{1},\lambda_{2}}$, $\lambda_{1}<p-1$;
  \item Irreducible finite dimensional representations $F_{\lambda_{1},\lambda_{2}}$, $\lambda_{1}<p-1$;
  \item $p=2$, Lowest weight discrete series $D_{\lambda_{1},\lambda_{2}}$, $\lambda_{2}\geq 0$.
\end{enumerate}

By Theorem \ref{Stable Ranger Theorem} and Theorem \ref{Going-Up}, the problem is reduced to calculating the theta $(p,2)$-lifts of the elements listed above. We denote by $\mathbf{G}_{p,2}$ the set $\{\overline{D}_{\lambda_{1},\lambda_{2}},F_{\lambda_{1},\lambda_{2}}|\lambda_{1}<p-1\}$. If $\lambda_{1}<p-1$, by Theorem \ref{Admissible Dual2} and Lemma \ref{Correspondence in Joint Harmonics}, the unique lowest degree $U(2)$-type of $\overline{D}_{\lambda_{1},\lambda_{2}}$ (resp. $F_{\lambda_{1},\lambda_{2}}$) is $(\lambda_{1}-1,\lambda_{2})$. By Theorem \ref{Correspondence in Joint Harmonics}, Lemma \ref{Correspondence in Joint Harmonics} and Theorem \ref{Lowest and Lowest Degree K-types}, we know that the unique lowest $Sp(p)\times Sp(2)$-type of $\theta_{p,2}(\overline{D}_{\lambda_{1},\lambda_{2}})$ (resp. $\theta_{p,2}(F_{\lambda_{1},\lambda_{2}})$) is of the form $(0,\cdots,0;b_{1},b_{2})$. We denote by $\theta_{p,2}(\mathbf{G}_{p,2})$ the set of theta $(p,2)$-lifts of elements in $\mathbf{G}_{p,2}$ and denote by $\mathbf{B}_{p,2}$ the set of infinitesimal equivalence classes in $\mathcal{R}(Sp(p,2))$ whose lowest $Sp(p)\times Sp(2)$-types are all of the form $(0,\cdots,0;b_{1},b_{2})$. By the appendix, we have the following theorem.

\begin{thm}\label{List for Theta (p,2)-lifts}
The following is a list of infinitesimal equivalence classes in $\mathbf{B}_{p,2}$:
\begin{enumerate}[(1)]
  \item $\pi((p,p-1,\cdots,1;\delta_{1},\delta_{2}),\Psi_{5})$, $\delta_{1}>\delta_{2}\geq p$;
  \item $\pi((p-1,\cdots,\delta_{2}+1,\delta_{2},\delta_{2},\delta_{2}-1,\cdots,1;\delta_{1},\delta_{2}),\Psi)$, $\delta_{1}\geq p-1\geq \delta_{2}$. Here $\Psi$ is the unique system of positive roots such that condition $(A)$ is satisfied and $-e_{1}+f_{1}, e_{1}-f_{2}\in \Psi$;
  \item $\pi((\lambda_{1},\lambda_{2},\cdots,\lambda_{p};\delta_{1},\delta_{2}),\Psi)$, $\delta_{1}\leq \lambda_{1}=p-2$. Here $\Psi$ is the unique system of positive roots such that condition $(A)$ is satisfied;
  \item $\pi(1,(p-1,p-2,\cdots,1;\delta),\Psi_{2},\mu,\nu)$, $\delta\geq p-1$. Here $\mu\leq 2p-1$ and the equality holds only if $\delta>p-1$;
  \item $\pi(1,(p-2,\cdots,\delta+1,\delta,\delta,\delta-1,\cdots,1;\delta),\Psi,\mu,\nu)$, $\delta\leq p-2$, $\mu\leq 2p-3$. Here $\Psi$ is uniquely determined by condition $(A)$;
  \item $\pi(2,(p-2,p-3,\cdots,1),\Psi_{1},(\mu_{1},\mu_{2}),(\nu_{1},\,\nu_{2}))$, $\mu_{1}\geq \mu_{2}$ and $\mu_{1}\leq 2p-3$.
\end{enumerate}
\end{thm}

We denote by $\mathbf{B}_{p,2}^{(1)}$ the set of infinitesimal equivalence classes in Case (1) of Theorem \ref{List for Theta (p,2)-lifts}. Similarly, we define the sets $\mathbf{B}_{p,2}^{(2)}$, $\mathbf{B}_{p,2}^{(3)}$, $\mathbf{B}_{p,2}^{(4)}$, $\mathbf{B}_{p,2}^{(5)}$ and $\mathbf{B}_{p,2}^{(6)}$. Denote $\mathbf{B}_{p,2}^{+}=\mathbf{B}_{p,2}^{(1)}\cup \mathbf{B}_{p,2}^{(4)}\cup \mathbf{B}_{p,2}^{(6)}$ and $\mathbf{B}_{p,2}^{-}=\mathbf{B}_{p,2}^{(2)}\cup \mathbf{B}_{p,2}^{(3)}\cup \mathbf{B}_{p,2}^{(5)}$. We have the following theorem.

\begin{thm}\label{Exclusion}
The set $\theta_{p,2}(\mathbf{G}_{p,2})$ is a subset of $\mathbf{B}_{p,2}^{+}$.
\end{thm}
\begin{pf}
Since $\theta_{p,2}(\mathbf{G}_{p,2})$ is a subset of $\mathbf{B}_{p,2}$ and $\mathbf{B}_{p,2}=\mathbf{B}_{p,2}^{+}\cup \mathbf{B}_{p,2}^{-}$, we only need to prove that $\theta_{p,2}(\mathbf{G}_{p,2})\cap \mathbf{B}_{p,2}^{-}=\emptyset$.

Suppose $\pi$ is an element in $\mathbf{B}_{p,2}^{(3)}$. The absolute values of all entries in the infinitesimal character of $\pi$ are less than $p$. On the other hand, by Theorem \ref{Correspondence of Infinitesimal Characters}, the infinitesimal character of $\theta_{p,2}(\overline{D}_{\lambda_{1},\lambda_{2}})$ (resp.  $\theta_{p,2}(F_{\lambda_{1},\lambda_{2}})$) is $(\lambda_{1},\lambda_{2},p,p-1,\cdots,1)$. Then $\pi$ is not the theta $(p,2)$-lift of an element in $\mathbf{G}_{p,2}$. Then we prove that $\theta_{p,2}(\mathbf{G}_{p,2})\cap \mathbf{B}_{p,2}^{(3)}=\emptyset$.

Suppose $\pi$ is an element in $\mathbf{B}_{p,2}^{(2)}$ as in Theorem \ref{List for Theta (p,2)-lifts}. By Theorem \ref{Correspondence of Infinitesimal Characters}, the infinitesimal character of $\pi$ corresponds to the infinitesimal character of an element in $\mathbf{G}_{p,2}$ only if $\delta_{1}=p$. If $\delta_{1}=p$, the unique element in $\mathbf{G}_{p,2}$ whose infinitesimal character corresponds to the infinitesimal character of $\pi$ is $\overline{D}_{\delta_{2},-\delta_{2}}$. By Theorem \ref{Algorithm1}, the lowest $Sp(p)\times Sp(2)$-type of $\pi$ is $(0,\cdots,0;2p-2,2\delta_{2}-1)$ while the lowest degree $U(2)$-type of $\overline{D}_{\delta_{2},-\delta_{2}}$ is $(\delta_{2}-1,-\delta_{2})$. By Lemma \ref{Correspondence in Joint Harmonics}, the $Sp(p)\times Sp(2)$-type $(0,\cdots,0;2p-2,2\delta_{2}-1)$ does not correspond to the $U(2)$-type $(\delta_{2}-1,-\delta_{2})$ in $\mathscr{H}_{p,2,2}$. By Theorem \ref{Joint Harmonics} and Theorem \ref{Lowest and Lowest Degree K-types}, $\pi$ is not the theta $(p,2)$-lifts of $\overline{D}_{\delta_{2},-\delta_{2}}$. We prove that $\theta_{p,2}(\mathbf{G}_{p,2})\cap \mathbf{B}_{p,2}^{(2)}=\emptyset$.

Suppose $\pi$ is an element in $\mathbf{B}_{p,2}^{(5)}$ as in Theorem \ref{List for Theta (p,2)-lifts}. By Theorem \ref{Correspondence of Infinitesimal Characters}, the infinitesimal character of $\pi$ corresponds to the infinitesimal character of an element in $\mathbf{G}_{p,2}$ only if $(\mu,\nu)=(1,2p-1)$. If $(\mu,\nu)=(1,2p-1)$, the unique element in $\mathbf{G}_{p,2}$ whose infinitesimal character corresponds to the infinitesimal character of $\pi$ is $\overline{D}_{\delta,-\delta}$. By Theorem \ref{Algorithm1}, the lowest $Sp(p)\times Sp(2)$-type of $\pi$ is $(0,\cdots,0;2\delta-1,0)$ while the lowest degree $U(2)$-type of $\overline{D}_{\delta,-\delta}$ is $(\delta-1,-\delta)$. By Lemma \ref{Correspondence in Joint Harmonics}, the $Sp(p)\times Sp(2)$-type $(0,\cdots,0;2\delta-1,0)$ does not correspond to the $U(2)$-type $(\delta-1,-\delta)$ in $\mathscr{H}_{p,2,2}$. By Theorem \ref{Joint Harmonics} and Theorem \ref{Lowest and Lowest Degree K-types}, $\pi$ is not the theta $(p,2)$-lifts of $\overline{D}_{\delta,-\delta}$. Then we prove that $\theta_{p,2}(\mathbf{G}_{p,2})\cap \mathbf{B}_{p,2}^{(5)}=\emptyset$.

In summary, we prove that $\theta_{p,2}(\mathbf{G}_{p,2})\cap \mathbf{B}_{p,2}^{-}=\emptyset$.
\end{pf}

\begin{thm}
The following is a list of theta $(p,2)$-lifts of elements in $\mathbf{G}_{p,2}$:
\begin{enumerate}[(1)]
  \item If $\lambda_{1},\lambda_{2}\leq -p$, then
\begin{equation*}
  \theta_{p,2}(\overline{D}_{\lambda_{1},\lambda_{2}})=\pi((p,p-1,\cdots,1;-\lambda_{2},-\lambda_{1}),\Psi_{5}).
\end{equation*}
  \item If $\lambda_{2}\leq -p+1\leq \lambda_{1}$, then
\begin{equation*}
  \theta_{p,2}(\overline{D}_{\lambda_{1},\lambda_{2}})=\pi(1,(p-1,p-2,\cdots,1;-\lambda_{2}),\Psi_{2},p-\lambda_{1},p+\lambda_{1}).
\end{equation*}
  \item If $\lambda_{2}\geq -p+2$, then
\begin{equation*}
  \theta_{p,2}(\overline{D}_{\lambda_{1},\lambda_{2}})=\pi(2,(p-2,p-3,\cdots,1),\Psi_{1},(\mu_{1},\mu_{2}),(\nu_{1},\nu_{2}))
\end{equation*}
with $(\mu_{1},\mu_{2})=(p-1-\lambda_{2},p-\lambda_{1})$ and $(\nu_{1},\nu_{2})=(p-1+\lambda_{2},p+\lambda_{1}))$.
  \item If $\lambda_{1}<p-1$, then
\begin{equation*}
\theta_{p,2}(F_{\lambda_{1},\lambda_{2}})=\pi(2,(p-2,p-3,\cdots,1),\Psi_{1},(\mu_{1},\mu_{2}),(\nu_{1},\nu_{2}))
\end{equation*}
with $(\mu_{1},\mu_{2})=(p-1-\lambda_{2},p-\lambda_{1})$ and $(\nu_{1},\nu_{2})=(p-1+\lambda_{2},p+\lambda_{1}))$.
\end{enumerate}
\end{thm}
\begin{pf}
By Theorem \ref{Admissible Dual2}, Lemma \ref{Correspondence in Joint Harmonics}, Theorem \ref{Correspondence of Infinitesimal Characters} and Theorem \ref{Lowest and Lowest Degree K-types},  we know that $\theta_{p,2}(\overline{D}_{\lambda_{1},\lambda_{2}})$ (resp. $\theta_{p,2}(F_{\lambda_{1},\lambda_{2}})$) has a unique lowest $Sp(p)\times Sp(2)$-type $(0,\cdots,0;p-2-\lambda_{2},p-1-\lambda_{1})$ and its infinitesimal character is $(\lambda_{1},\lambda_{2},p,p-1,\cdots,1)$. By Theorem \ref{Exclusion}, there is a unique element in $\mathbf{B}_{p,2}^{+}$ satisfying these condition. We list the explicit theta $(p,2)$-lifts above.
\end{pf}

We denote by $\mathbf{S}_{2,2}$ the set $\{D_{\lambda_{1},\lambda_{2}}|\lambda_{2}\geq 0\}$, and by $\mathbf{C}_{2,2}$ the set of elements in $\mathcal{R}(Sp(2,2))$ whose lowest $Sp(2)\times Sp(2)$-types are all of the form $(a_{1},a_{2};0,0)$. In the appendix, we list all infinitesimal equivalence classes in $\mathbf{C}_{2,2}$ in terms of Langlands parameters.

\begin{thm}\label{Theta $(2,2)$-Lifts}
The following is a list of explicit theta $(2,2)$-lifts of elements in $\mathbf{S}_{2,2}$:
\begin{enumerate}[(1)]
  \item If $\lambda_{2}=0$, then
\begin{equation*}
  \theta_{2,2}(D_{\lambda_{1},0})=\pi(1,(\lambda_{1};1),\Psi_{3},2,2).
\end{equation*}
  \item If $\lambda_{2}=1$, then
\begin{equation*}
  \theta_{2,2}(D_{\lambda_{1},1})=\pi(1,(\lambda_{1};1),\Psi_{3},3,1).
\end{equation*}
  \item If $\lambda_{2}\geq 2$, then
\begin{equation*}
  \theta_{2,2}(D_{\lambda_{1},\lambda_{2}})=\pi((\lambda_{1},\lambda_{2};2,1),\Psi_{6}).
\end{equation*}
\end{enumerate}
\end{thm}
\begin{pf}
By Theorem \ref{Admissible Dual2}, Lemma \ref{Correspondence in Joint Harmonics}, Theorem \ref{Correspondence of Infinitesimal Characters} and Theorem \ref{Lowest and Lowest Degree K-types},  we know that $\theta_{2,2}(D_{\lambda_{1},\lambda_{2}})$ has a unique lowest $Sp(2)\times Sp(2)$-type $(\lambda_{1},\lambda_{2}+1;0,0)$ and its infinitesimal character is $(\lambda_{1},\lambda_{2},2,1)$. By the appendix, there is a unique element in $\mathbf{C}_{2,2}$ satisfying these conditions. We list the explicit theta $(2,2)$-lifts above.
\end{pf}

\appendix

\section*{Appendix}

In this appendix, we list infinitesimal equivalence classes in $\mathbf{A}_{p,1}$, $\mathbf{B}_{p,2}$ and $\mathbf{C}_{2,2}$ in terms of Langlands parameters for $p\geq 2$. Let $\Psi_{1}$, $\Psi_{2}$, $\Psi_{3}$, $\Psi_{4}$, $\Psi_{5}$ and $\Psi_{6}$ be those mentioned in $\S$5 and $\S$6. Let condition (A) and (B) be those of $\S$2.1.

First, we list infinitesimal equivalence classes in $\mathbf{A}_{p,1}$ in terms of Langlands parameters. It is obvious that $\mathbf{A}_{1,1}=\mathcal{R}(Sp(1,1))$. Assume $p\geq 2$. By Theorem \ref{Algorithm1}, the following is a list of infinitesimal equivalence classes in $\mathbf{A}_{p,1}$ for $p\geq 2$:
\begin{enumerate}[(1)]
  \item $\pi((\lambda_{1},p-1,p-2,\cdots,1;\delta_{1}),\Psi_{2})$, $\delta_{1}\geq \lambda_{1}>p-1$;
  \item $\pi((\lambda_{1},p-1,p-2,\cdots,1;\delta_{1}),\Psi_{4})$, $\lambda_{1}\geq \delta_{1}>p-1$;
  \item $\pi((\lambda_{1},p-1,p-2,\cdots,1;p-1),\Psi_{4})$, $\lambda_{1}\geq p-1$;
  \item $\pi((\lambda_{1},p-2,\cdots,\delta_{1}+1,\delta_{1},\delta_{1},\delta_{1}-1,\cdots,1;\delta_{1}),\Psi)$, $\delta_{1}\leq p-2$. Here $\Psi$ is uniquely determined by condition (A);
  \item $\pi(1,(p-1,p-2,\cdots,1),\Psi_{1},\mu,\nu)$;
  \item $\pi(1,(\lambda_{1},p-2,p-3,\cdots,1),\Psi_{1},\mu,\nu)$, $\lambda_{1}>p-1$ and $\mu\leq 2p-3$.
\end{enumerate}

The desired theta $(p,1)$-lifts in Theorem \ref{Explicit Theta $(p,1)$-Lifts of Characters}, Theorem \ref{Explicit Theta $(p,1)$-lifts of Highest Weight Discrete Series} and Theorem \ref{Explicit Theta $(p,1)$-Lifts of Lowest Weight Discrete Series} are elements in $\mathbf{A}_{p,1}$. We determined these theta $(p,1)$-lifts in terms of Langlands parameters by their lowest $Sp(p)\times Sp(1)$-types and infinitesimal characters.

Next we list infinitesimal equivalence classes in $\mathbf{B}_{p,2}$ in terms of Langlans parameters. Assume $p\geq 2$. For $\pi\in \mathcal{R}(Sp(p,2))$, we denote by $\mathcal{A}^{1}(\pi)$ the set $\{a_{1}|(a_{1},a_{2},\cdots,a_{p};b_{1},b_{2})\in \mathcal{A}(\pi)\}$. We know that $\pi\in \mathbf{B}_{p,2}$ if and only if $\mathcal{A}^{1}(\pi)=\{0\}$.

Suppose $\pi=\pi((\lambda_{1},\lambda_{2},\cdots,\lambda_{p};\delta_{1},\delta_{2}),\Psi)$.  By Theorem \ref{Algorithm1}, the set $\mathcal{A}^{1}(\pi)$ is:
\begin{equation*}
\begin{split}
 \left\{
   \begin{array}{ll}
     \{\lambda_{1}-p\} & \hbox{if $\lambda_{1}<\delta_{2}$,} \\
     \{\lambda_{1}-p+1\} & \hbox{if $\delta_{2}<\lambda_{1}<\delta_{1}$,} \\
     \{\lambda_{1}-p+2\} & \hbox{if $\delta_{1}<\lambda_{1}$,} \\
     \{\lambda_{1}-p\}      & \hbox{if $\delta_{1}>\delta_{2}=\lambda_{1}>\lambda_{2}$ and $-e_{1}+f_{2}\in \Psi$,} \\
     \{\lambda_{1}-p+1\}      & \hbox{if $\delta_{1}>\delta_{2}=\lambda_{1}>\lambda_{2}$ and $e_{1}-f_{2}\in \Psi$,} \\
     \{\lambda_{1}-p+1\}      & \hbox{if $\delta_{1}=\lambda_{1}>\lambda_{2},\delta_{2}$ and $-e_{1}+f_{1}\in \Psi$,} \\
     \{\lambda_{1}-p+2\}      & \hbox{if $\delta_{1}=\lambda_{1}>\lambda_{2},\delta_{2}$ and $e_{1}-f_{1}\in \Psi$,} \\
     \{\lambda_{1}-p+2\} & \hbox{if $\delta_{1}=\lambda_{1}=\lambda_{2}>\delta_{2}$,} \\
     \{\lambda_{1}-p+1\} & \hbox{if $\delta_{2}=\lambda_{1}=\lambda_{2}<\delta_{1}$,} \\
     \{\lambda_{1}-p+1\} & \hbox{if $\lambda_{1}=\delta_{1}=\delta_{2}>\lambda_{2}$,} \\
     \{\lambda_{1}-p+2\} & \hbox{if $\lambda_{1}=\lambda_{2}=\delta_{1}=\delta_{2}$ and $e_{1}-f_{1}\in \Psi$,} \\
     \{\lambda_{1}-p+1\} & \hbox{if $\lambda_{1}=\lambda_{2}=\delta_{1}=\delta_{2}$ and $-e_{1}+f_{1}\in \Psi$.}
   \end{array}
 \right.
\end{split}
\end{equation*}

Suppose  $\pi=\pi(1,(\lambda_{1},\lambda_{2},\cdots,\lambda_{p-1};\delta_{1}),\Psi,\mu_{1},\nu_{1})$. By Theorem \ref{Algorithm1}, the set $\mathcal{A}^{1}(\pi)$ is:
\begin{equation*}
\begin{split}
  \left\{
    \begin{array}{ll}
      \{\lambda_{1}-p+2\} & \hbox{if $\lambda_{1}>\delta_{1}$ and $\lambda_{1}>\frac{\mu_{1}}{2}$,} \\
      \{\lambda_{1}-p+1\}, & \hbox{if $\frac{\mu_{1}}{2}<\lambda_{1}<\delta_{1}$,} \\
      \{\frac{\mu_{1}+3}{2}-p\} & \hbox{if $\delta_{1},\lambda_{1}<\frac{\mu_{1}}{2}$ and $\mu_{1}$ is odd,} \\
      \{\frac{\mu_{1}+1}{2}-p\} & \hbox{if $\delta_{1}>\frac{\mu_{1}}{2}>\lambda_{1}$ and $\mu_{1}$ is odd,} \\
      \{\frac{\mu_{1}}{2}-p+1,\frac{\mu_{1}}{2}-p+2\} & \hbox{if $\delta_{1},\lambda_{1}<\frac{\mu_{1}}{2}$ and $\mu_{1}$ is even,} \\
      \{\frac{\mu_{1}}{2}-p,\frac{\mu_{1}}{2}-p+1\} & \hbox{if $\delta_{1}>\frac{\mu_{1}}{2}>\lambda_{1}$ and $\mu_{1}$ is even,} \\
      \{\lambda_{1}-p+2\} & \hbox{if $\lambda_{1}=\delta_{1}\geq \frac{\mu_{1}}{2}$ and $e_{1}-f_{1}\in \Psi$,} \\
     \{\lambda_{1}-p+1\} & \hbox{if $\lambda_{1}=\delta_{1}\geq \frac{\mu_{1}}{2}$ and $-e_{1}+f_{1}\in \Psi$,} \\
      \{\lambda_{1}-p+2\} & \hbox{if $\lambda_{1}=\frac{\mu_{1}}{2}>\delta_{1}$ ,} \\
      \{\lambda_{1}-p+1\} & \hbox{if $\lambda_{1}=\frac{\mu_{1}}{2}<\delta_{1}$ ,} \\
      \{\frac{\mu_{1}}{2}-p+1\} & \hbox{if $\delta_{1}=\frac{\mu_{1}}{2}>\lambda_{1}$ .}
\end{array}
  \right.
\end{split}
\end{equation*}

Suppose $p\geq 3$ and $\pi=\pi(2,(\lambda_{1},\lambda_{2},\cdots,\lambda_{p-2}),\Psi,(\mu_{1},\mu_{2}),(\nu_{1},\nu_{2}))$ with $\mu_{1}\geq \mu_{2}$. By Theorem \ref{Algorithm1}, the set $\mathcal{A}^{1}(\pi)$ is :
\begin{equation*}
\begin{split}
 \left\{
    \begin{array}{ll}
      \{\lambda_{1}-p+2\} & \hbox{if $\lambda_{1}>\frac{\mu_{1}}{2}$,} \\
      \{\frac{\mu_{1}+3}{2}-p\} & \hbox{if $\lambda_{1}<\frac{\mu_{1}}{2}$ and $\mu_{1}$ is odd,} \\
        \{\frac{\mu_{1}}{2}-p+1,\frac{\mu_{1}}{2}-p+2\} & \hbox{if $\lambda_{1}<\frac{\mu_{1}}{2}$ and $\mu_{1}$ is even,} \\
        \{\lambda_{1}-p+2\} & \hbox{if $\lambda_{1}=\frac{\mu_{1}}{2}$,}
    \end{array}
  \right.
\end{split}
\end{equation*}

Suppose $p=2$ and $\pi=\pi(2,0,\emptyset,(\mu_{1},\mu_{2}),(\nu_{1},\nu_{2}))$ with $\mu_{1}\geq \mu_{2}$. By Theorem \ref{Algorithm1}, the set $\mathcal{A}^{1}(\pi)$ is:
\begin{equation*}
  \left\{
     \begin{array}{ll}
       \{\frac{\mu_{1}-1}{2}\} & \hbox{if $\mu_{1}$ is odd,} \\
       \{\frac{\mu_{1}}{2}-1,\frac{\mu_{1}}{2}\} & \hbox{if $\mu_{1}$ is even.}
     \end{array}
   \right.
\end{equation*}

In summary, the following is a list of infinitesimal equivalence classes in $\mathbf{B}_{p,2}$:
\begin{enumerate}[(1)]
  \item $\pi((p,p-1,\cdots,1;\delta_{1},\delta_{2}),\Psi_{5})$, $\delta_{2}\geq p$;
  \item $\pi((p-1,\cdots,\delta_{2}+1,\delta_{2},\delta_{2},\delta_{2}-1,\cdots,1;\delta_{1},\delta_{2}),\Psi)$, $\delta_{2}\leq p-1\leq \delta_{1}$. Here $\Psi$ is the unique system of positive roots such that condition (A) is satisfied and $-e_{1}+f_{1},e_{1}-f_{2}\in \Psi$;
  \item $\pi((p-2,\lambda_{2},\lambda_{3},\cdots,\lambda_{p};\delta_{1},\delta_{2}),\Psi)$, $\delta_{1}\leq p-2$. Here $\Psi$ the unique system of positive roots such that condition (A) is satisfied;
  \item $\pi(1,(p-1,p-2,\cdots,1;\delta_{1}),\Psi_{2},\mu_{1},\nu_{1})$, $\delta_{1}>p-1$ and $\frac{\mu_{1}}{2}\leq p-\frac{1}{2}$;
  \item $\pi(1,(p-1,p-2,\cdots,1;p-1),\Psi_{2},\mu_{1},\nu_{1})$, $\frac{\mu_{1}}{2}\leq p-1$;
  \item $\pi(1,(p-2,\cdots,\delta_{1}+1,\delta_{1},\delta_{1},\delta_{1}-1,\cdots,1;\delta_{1}),\Psi,\mu_{1},\nu_{1})$, $\delta_{1}\leq p-2$, $\frac{\mu_{1}}{2}\leq p-\frac{3}{2}$. Here $\Psi$ is uniquely determined by condition (A);
  \item $\pi(2,(p-2,p-3,\cdots,1),\Psi_{1},(\mu_{1},\mu_{2}),(\nu_{1},\nu_{2}))$, $\mu_{2}\leq \mu_{1}\leq 2p-3$.
\end{enumerate}

The infinitesimal equivalence classes listed above are just those listed in Theorem \ref{List for Theta (p,2)-lifts}.

Finally, we list infinitesimal equivalence classes in $\mathbf{C}_{2,2}$ in terms of Langlands parameters. For $\pi\in \mathcal{R}(Sp(2,2))$, we denote by $\mathcal{A}^{3}(\pi)$ the set $\{b_{1}|(a_{1},a_{2};b_{1},b_{2})\in \mathcal{A}(\pi)\}$. We know that $\pi\in \mathbf{C}_{2,2}$ if and only if $\mathcal{A}^{3}(\pi)=\{0\}$. Let $\Psi_{7}$ be a system of positive roots for the group $Sp(2,2)$ such that $(4,2;3,1)$ is dominant with respect to $\Psi_{7}$.

Suppose $\pi=\pi((\lambda_{1},\lambda_{2};\delta_{1},\delta_{2}),\Psi)$. By Theorem \ref{Algorithm1}, the set $\mathcal{A}^{3}(\pi)$ is
\begin{equation*}
\left\{
  \begin{array}{ll}
    \{\delta_{1}-2\}, & \hbox{if $\delta_{1}<\lambda_{2}$;} \\
    \{\delta_{1}-1\}, & \hbox{if $\lambda_{2}<\delta_{1}<\lambda_{1}$;} \\
    \{\delta_{1}\}, & \hbox{if $\delta_{1}>\lambda_{1}$;} \\
    \{\delta_{1}-2\}, & \hbox{if $\lambda_{1}>\lambda_{2}=\delta_{1}>\delta_{2}$ and $e_{2}-f_{1}\in \Psi$;} \\
    \{\delta_{1}-1\}, & \hbox{if $\lambda_{1}>\lambda_{2}=\delta_{1}>\delta_{2}$ and $-e_{2}+f_{1}\in \Psi$;} \\
   \{\delta_{1}-1\}, & \hbox{if $\lambda_{2},\delta_{2}<\delta_{1}=\lambda_{1}$ and $e_{1}-f_{1}\in \Psi$;} \\
   \{\delta_{1}\}, & \hbox{if $\lambda_{2},\delta_{2}<\delta_{1}=\lambda_{1}$ and $-e_{1}+f_{1}\in \Psi$;} \\
   \{\delta_{1}-1\}, & \hbox{if $\lambda_{1}=\lambda_{2}=\delta_{1}>\delta_{2}$;} \\
  \{\delta_{1}-1\} , & \hbox{if $\lambda_{1}>\lambda_{2}=\delta_{1}=\delta_{2}$;} \\
  \{\delta_{1}\}, & \hbox{if $\lambda_{1}=\delta_{1}=\delta_{2}>\lambda_{2}$;} \\
  \{\delta_{1}-1\}, & \hbox{if $\lambda_{1}=\lambda_{2}=\delta_{1}=\delta_{2}$ and $e_{1}-f_{1}\in \Psi$;} \\
   \{\delta_{1}\}, & \hbox{if $\lambda_{1}=\lambda_{2}=\delta_{1}=\delta_{2}$ and $-e_{1}+f_{1}\in \Psi$.}
  \end{array}
\right.
\end{equation*}

Suppose $\pi=\pi(1,(\lambda_{1};\delta_{1}),\Psi,\mu,\nu)$. By Theorem \ref{Algorithm1}, the set $\mathcal{A}^{3}(\pi)$ is
\begin{equation*}
  \left\{
    \begin{array}{ll}
      \{\frac{\mu-1}{2}\}, & \hbox{if $\lambda_{1},\delta_{1}<\frac{\mu}{2}$ and $\mu$ is odd;} \\
      \{\frac{\mu}{2},\frac{\mu}{2}-1\}, & \hbox{if $\lambda_{1},\delta_{1}<\frac{\mu}{2}$ and $\mu$ is even;} \\
      \{\delta_{1}-1\}, & \hbox{if $\lambda_{1}>\delta_{1}>\frac{\mu}{2}$;} \\
      \{\delta_{1}\}, & \hbox{if $\delta_{1}>\lambda_{1}$ and $\delta_{1}\geq \frac{\mu}{2}$;} \\
      \{\frac{\mu-3}{2}\}, & \hbox{if $\lambda_{1}>\frac{\mu}{2}>\delta_{1}$ and $\mu$ is odd;} \\
      \{\frac{\mu}{2}-1,\frac{\mu}{2}-2\}, & \hbox{if $\lambda_{1}>\frac{\mu}{2}>\delta_{1}$ and $\mu$ is even;} \\
      \{\frac{\mu}{2}-1\}, & \hbox{if $\frac{\mu}{2}=\lambda_{1}>\delta_{1}$;} \\
      \{\delta_{1}-1\}, & \hbox{if $\lambda_{1}>\delta_{1}=\frac{\mu}{2}$;} \\
      \{\delta_{1}-1\}, & \hbox{if $\lambda_{1}=\delta_{1}\geq \frac{\mu}{2}$ and $e_{1}-f_{1}\in \Psi$;} \\
      \{\delta_{1}\}, & \hbox{if $\lambda_{1}=\delta_{1}\geq \frac{\mu}{2}$ and $-e_{1}+f_{1}\in \Psi$.}
    \end{array}
  \right.
\end{equation*}

Suppose $\pi=\pi(2,0,\emptyset,(\mu_{1},\mu_{2}),(\nu_{1},\nu_{2}))$ with $\mu_{1}\geq \mu_{2}$. By Theorem \ref{Algorithm1}, the set $\mathcal{A}^{3}(\pi)$ is
\begin{equation*}
  \left\{
    \begin{array}{ll}
      \{\frac{\mu_{1}-1}{2}\}, & \hbox{if $\mu_{1}$ is odd;} \\
      \{\frac{\mu_{1}}{2},\frac{\mu_{1}}{2}-1\}, & \hbox{if $\mu_{1}$ is even.}
    \end{array}
  \right.
\end{equation*}

In summary, the following is a list of infinitesimal equivalence classes in $\mathbf{C}_{2,2}$:
\begin{enumerate}[(1)]
  \item $\pi((\lambda_{1},\lambda_{2};2,1),\Psi_{6})$;
  \item $\pi((\lambda_{1},1;1,1),\Psi_{7})$;
  \item $\pi(1,(\lambda_{1};1),\Psi_{3},3,\nu)$, $\lambda_{1}>1$;
  \item $\pi(1,(\lambda_{1};1),\Psi_{3},\mu,\nu)$, $\mu\leq 2$;
  \item $\pi(2,0,\emptyset,(1,1),(\nu_{1},\nu_{2}))$.
\end{enumerate}

\section*{Acknowledgement}
We would like to thank R. Howe and J.-S. Li for their suggestive strategy. We also would like to thank B. Sun for his key idea on conservation relations and Xiang Fan for heuristic discussions.








\bibliographystyle{model1b-num-names}
\bibliography{MyRef1}







\end{document}